%% file: PickloRyanNonUniSIACMRA.tex
\documentclass[11pt]{article}

 \usepackage{geometry}
\usepackage[utf8]{inputenc}
\usepackage{tabularx}
\usepackage{tikz}
\usepackage{amsmath,amsthm,amssymb,amsfonts}
\usetikzlibrary{patterns,patterns.meta}
\usetikzlibrary{arrows}
\usetikzlibrary{arrows.meta}

\usepackage{algorithm}
\usepackage{algpseudocode}
\usepackage{mathtools}
\usepackage{multirow}

\usepackage{mathrsfs}
\usepackage{graphicx}
\usepackage{adjustbox}
\usepackage{pstricks-add}
\usepackage{pgf}
\usepackage{lineno}
\usepackage[toc,page]{appendix}
\usepackage{bbm}
\usepackage{bbold}
\usepackage{lscape,rotating,pdflscape}
\usepackage{enumerate}
\usepackage{caption}
\usepackage{subcaption,wrapfig}
\usepackage{float,multirow}
\usepackage{comment}
\usepackage{array}
 \usepackage{relsize}
\usepackage{lmodern,babel,booktabs,multirow}
\usepackage{setspace}
\usepackage{pdflscape}
\usepackage{colortbl}
\usepackage{mathdots}
\usepackage{xcolor}
 \usepackage{blkarray,stackengine,tikz}
 \usepackage{trimclip}
 \usepackage{footnote}
 \usepackage[bottom]{footmisc}
 \usepackage[normalem]{ulem}
 \usepackage{hyperref}

 \newcommand{\norm}[1]{\left\lVert#1\right\rVert}
 
 \usetikzlibrary{arrows}
\usetikzlibrary{decorations.markings}
\definecolor{DarkRed}{rgb}{0.65,0.00,0.05}
\definecolor{structure}{rgb}{0.2,0.2,0.7}

\newcolumntype{M}[1]{>{\centering\arraybackslash}m{#1}}
\newcolumntype{N}{@{}m{0pt}@{}}

\providecommand{\keywords}[1]
{
  \small	
  \textbf{\textit{Keywords---}} #1
}

\providecommand{\AMS}[1]
{
  \small	
  \textbf{\textit{AMS---}} #1
}

\title{LSIAC-MRA for Nonuniform Meshes and Applications to Mesh Adaptivity}

\author{Matthew J. Picklo{\footnotemark[1] \footnotemark[2] } \and
Jennifer K. Ryan{ \footnotemark[2] \footnotemark[3] \footnotemark[4]} 
}

\footnotetext[1]{Department of Applied Mathematics and Statistics, Colorado School of Mines, Golden, CO 80401, USA.}
\footnotetext[2]{Research supported by Air Force Office of Scientific Research (AFOSR), Computational Mathematics Program (Program Manager: Dr. Fariba Fahroo), under grant number FA9550-20-1-0166.}

\footnotetext[3]{{\it Corresponding Author.} {Email:} jryan@kth.se}
\footnotetext[4]{Department of Mathematics and Linn\'{e} Flow Centre, KTH Royal Institute of Technology, 114 28 Stockholm, Sweden.}

\date{\today}

\begin{document}

\maketitle

\begin{abstract}
In this article we consider the extension of the (L)SIAC-MRA enhancement procedure to nonuniform meshes. We demonstrate that error reduction can be obtained on perturbed quadrilateral and Delaunay meshes, and investigate the effect of limited resolution and its impact on the procedure for various function types. We show that utilizing mesh-based localized kernel scalings, which were shown to reduce approximation errors for LSIAC filters, improve the performance of the LSIAC-MRA enhancement procedure. Lastly, we demonstrate the usefulness of enhanced approximations generated by (L)SIAC-MRA in mesh adaptivity applications, and show that SIAC reconstruction can be used in identification of regions of high error in steady-state DG approximations.

\end{abstract}


\AMS{65D60, 65D10, 65M50, 65M60} 

\input{sections/introduction.tex}

\input{sections/enhancement_procedure.tex}

\input{sections/Mesh_Transition}

\input{sections/adaptive_schemes}

\input{sections/conclusions.tex}

\section{Acknowledgements}
We would like to thank Dr. Ayaboe Edoh for his valuable comments in the presentation of this paper.  This work is done in connection with the Linné Flow Centre at KTH. Research of both authors was supported by the Air Force Office of Scientific Research (AFOSR) Computational Mathematics Program (Program Manager: Dr. Fariba Fahroo) under Grant number FA9550-20-1-0166.

\bibliographystyle{siamplain.bst}
\bibliography{PickloRyanNonUniSIACMRA.bib}
\end{document}

%% file: sections/introduction.tex
\section{Introduction}\label{sec:intro}

The ability to utilize multiple levels of resolution, including subgrid scales, is useful in a variety of applications such as trouble cell detection \cite{Vuik2014,Vuik2017}, multiresolution numerical methods \cite{Hovhannisyan2014,Muller,Gerhard2022}, and error estimation for $h/p$ adaptivity \cite{Bautista,bautista_thesis,Bautista_paper}. Through application of carefully selected Smoothness-Increasing Accuracy-Conserving (SIAC) kernels, the authors of \cite{SIACMRA,picklo} were able to demonstrate error reduction under multiple mesh transitions for a nested sequence of uniform hierarchical meshes. The error reduction was obtained for quadrilateral Cartesian meshes through the use of a SIAC kernel post-processor prior to projection onto a refined mesh. Furthermore, the authors constructed explicit transition matrices for performing the procedure efficiently. However, the approach only considered uniform meshes.  Extending this procedure to nonuniform meshes and adapted meshes presents many challenges because repeating patterns in the mesh and kernel breaks do not exist globally. This unfortunately limits the use of the procedure for mesh adaptivity purposes. In this work, we return to the fundamental filtering and projection procedures and extend these techniques to triangular and nonuniform meshes. Furthermore, we demonstrate the ability of (Line)-SIAC Multiresolution Analysis ((L)SIAC-MRA) to construct accurate fine mesh representations of coarse mesh data that can be applied within the context of adaptive Discontinuous Galerkin (DG) methods. 

\newpage
The (L)SIAC-MRA procedure views the application of SIAC filtering and projection within a multiresolution framework, where fine mesh data is represented as a coarse mesh projection plus additional multiwavelet details. In (L)SIAC-MRA, coarse data is convolved with a carefully chosen SIAC kernel and projected onto a refined mesh. This enhanced reconstruction represents a fine mesh approximation which is equivalent to a coarse mesh approximation with additional multiwavelet details:
\begin{equation}
\tilde{u}_h^{(n+1)}(\mathbf{x})=\tilde{u}^{(n)}_h(\mathbf{x})+\sum_{\tau\in\mathcal{T}^{(n)}}\sum_{k\in\mathcal{I}}\bar{d}^{n}_{k,\tau}\psi^{n}_{\tau,k}(\mathbf{x}).
\end{equation}
Here the $n+1$ superscript is used to indicate that the approximation is defined on a finer mesh than the coarse level-$n$ initial data.

This representation is typically much more accurate than simply projecting the approximation on the coarse mesh alone. The procedure was originally developed for uniform meshes using tensor-product kernels in \cite{SIACMRA}. It was then later extended to higher dimension via a Line SIAC filter in \cite{picklo}. In both cases, error reduction under iterated refinement was observed numerically. The latter paper demonstrated that this procedure is also useful for Finite Volume type approximations. In this article, we introduce a refinement procedure that enables the same error reduction for unstructured meshes and applies that capability for refinement indication in adaptive DG schemes.

The SIAC filters themselves were originally developed in \cite{Bramble1977} for application to finite element data arising from the solution of elliptic PDEs. Their effectiveness was later demonstrated for the hyperbolic case in \cite{Cockburn2002}. To improve the computational efficiency in higher dimensions, the less costly Line-SIAC variants were introduced in \cite{LSIAC}. The crucial characteristic enabling this improved efficiency was the substitution of a one-dimensional kernel for a multidimensional tensor-product kernel.  The tensor-product kernel requires spatial quadrature of the same dimension as the data, whereas, with a line kernel, the support is single parameterized line segment only requiring one-dimensional quadrature. The orientation of the line segment is chosen to allow the filter to simultaneously traverse multiple ordinates, which with the correct orientation enables the same superconvergence phenomena associated with the more costly tensor-product construction. Additionally, the line filter has a reduction in the error constant.

The L-SIAC filter has since been augmented by an adaptive kernel scaling in \cite{Jallepalli2019B}. The adaptively scaled kernel uses a scaling choice that is based on the local average resolution of the mesh about the filtering location. Previously, the maximum mesh edge length has been chosen as the kernel scaling for nonuniform meshes. Unfortunately, when the mesh resolution varies significantly over the domain, exact evaluation of the convolution is particularly costly. By choosing the maximum edge length, the kernel support is fixed for the entire domain regardless of the local resolution.  As a consequence, in high resolution regions, many separate quadrature sums are required to ensure exactness. This contrasts with the low resolution regions which would typically require fewer quadrature sums because the kernel scaling is closer to the local characteristic length scale of the mesh. In the uniform case, the choice of kernel scaling directly influences the amount of dissipation applied by the filter, and so by choosing an adaptively scaled kernel, the magnitude of attenuation of data is also localized, preventing excessive dissipation in high resolution regions. Using the adaptive scaling, the authors of \cite{Jallepalli2019B} showed not only improved computational efficiency but also greater error reductions for unstructured meshes where the mesh resolution varied significantly. 

One direct application of (L)SIAC-MRA is mesh adaptivity and by extension approximation refinement. We show how incorporating this enhancement procedure into a multiwavelet-based framework can be useful for mesh adaptivity. Adaptive multiresolution schemes have previously been considered for Finite-Volume methods \cite{Muller,Harten1993b,Harten1994}, as well as in the DG-context \cite{Calle,Gerhard2015,Gerhard2016,Gerhard2022,NilsThesis2017,Shelton2008,Hovhannisyan2014}. These approaches have mainly concentrated on utilizing the multiresolution analysis framework for computational efficiency, specifically the truncation of insignificant detail coefficients, which can be thought of as extraneous degrees of freedom (DOF). The present approach is motivated by ideas similar to \cite{Bautista_paper,bautista_thesis}. Our (L)SIAC-MRA reconstruction is used to generate additional detail coefficients, which effectively defines a finer mesh approximation. These modified detail coefficients are then used to indicate regions where adaptivity is necessary. This work differs from those previously mentioned in that the SIAC enhancement methodologies are employed to increase resolution rather than projections of least-squares polynomials. 

We begin with a review of the enhancement procedure in one- and two-dimensions in Section \ref{sec:enhancement_procedure}, which includes an overview of multiresolution analysis as it relates to this work. We then demonstrate the  effectiveness of LSIAC-MRA for error reduction when the enhancement procedure applied to unstructured meshes in Section \ref{sec:Mesh_Transition}. As part of this, the role of an adaptive kernel scaling is investigated. Lastly, a demonstration of the usefulness of the enhancement technique in the context of mesh refinement schemes is provided in Section \ref{sec:adaptive_DG}.

%% file: sections/enhancement_procedure.tex
\section{Background} \label{sec:enhancement_procedure}

In this section we describe the relevant background for the enhancement procedure for the one- and two-dimensional (L)SIAC-MRA as originally formulated in \cite{SIACMRA, picklo}. We first provide an overview of multiresolution analysis as it pertains to this work. Next we detail the enhancement procedure, formulating the SIAC filters and projection components used in the enhancement. Briefly discussing the issue of adaptive scaling, this section concludes with instruction for computation of the multiwavelet detail coefficient to be used later for mesh adaptivity.

\subsection{Data format}
We first note that, for the purposes of this work, it is assumed that the data $u_h^{(n)}(\mathbf{x})$ is defined on a mesh $\mathcal{T}^{(n)}$ composed of elements $\tau$ that partition some domain $\Omega$. We further assume that on each element, $\tau$, the data can be expressed in terms of some polynomial basis $\{\phi_{k,\tau}\}_{k\in \mathcal{I}}$ with $\mathcal{I}$ being an index set for the basis. Hence, we write the approximation on an element $\tau$ as:
\begin{equation}\label{eq:initialization}
    u_h^{(n)}\Big|_{\tau}(\mathbf{x})=\sum_{k\in \mathcal{I}}\hat{u}^k_{\tau}\phi_{k,\tau}(\mathbf{x}).
\end{equation}
Throughout this article the initializations considered to obtain the $\hat{u}$ coefficients are determined by $L^2-$projection, that is by requiring that
\begin{equation}
\label{eq:L2}
    \langle u,\phi_{m,\tau}\rangle_{\tau}=\langle u_h^{(n)},\phi_{m,\tau}\rangle_{\tau},
\end{equation}
for each $m\in \mathcal{I}$ and $\tau\in\mathcal{T}^{(n)}$.

\subsection{Multiresolution Analysis}

In order to understand the enhancement procedure, we first discuss the multiresolution analysis framework (MRA).  It is important to note that this discussion is for illustrative purposes as the enhancement procedure does not actually make explicit use of the multiwavelet details.  For unstructured meshes, the multiwavelet bases may be expensive to compute (see Section \ref{sec:Mesh_Transition}).  These details are important in the adaptive solution of the PDEs considered in section \ref{sec:adaptive_DG} where the multiwavelet detail indicator of \cite{Bautista_paper} is used. For simplicity, we consider multiresolution analysis specifically for the case of hierarchies of dyadically uniformly subdivided meshes, and, in the multidimensional case, we will restrict ourselves to quadrilateral elements. Note that a multiresolution analyses naturally exists and can be constructed for other hierarchies (see Yu et al. for triangular mesh-based multiwavelets, and Gerhard et al. for wavelet-free methods for nonuniform meshes \cite{Yu1997,Gerhard2022}). For more in-depth treatments of MRA we refer the reader to \cite{Alp93,Alp02BGV,Wavelets2,Yu1997,Keinert2004}.

For clarity, we first consider the one-dimensional case and follow the presentation of \cite{Vuik2017,Bautista}. Define an initial mesh $\mathcal{T}^{(0)}=[-1,1]$ to consist of a single element. Let $\mathcal{T}^{(n)}$ denote the mesh obtained by uniformly subdividing every element in $\mathcal{T}^{(n-1)}$, i.e.
\begin{equation}
\mathcal{T}^{(n)}=\cup_{j=0}^{2^n-1}I^n_j,
\end{equation}
where $I^n_j=(-1+2^{1-n}j,-1+2^{1-n}(j+1))$. Over each mesh, we define the approximation space
\begin{equation}
V^p_n=\{v\;:v\big|_{I^n_j}\in\mathbb{P}^p(I^n_j),j=0,\hdots,2^n-1 \}.
\end{equation}
Owing to the dyadic hierarchy of the meshes, the approximation spaces are nested:
\begin{equation}
V^p_{0}\subset V^p_{1}\subset\hdots\subset V^p_n \subset \hdots. 
\end{equation}
Denote by $\{\phi_k(x)\}_{k=0}^p=\{\phi^0_{k,0}(x)\}_{k=0}^p$ an orthonormal basis for the coarsest mesh consisting of one element. These basis functions are called \textit{scaling functions}. We can obtain bases for any $V^p_n$ by scaling and translating this initial basis of scaling functions. The approximation space $V^p_n$ can also be defined as the span of $\{\phi^n_{k,j}\}_{k,j}$ where 
\begin{equation}
\phi^n_{k,j}=2^{n/2}\phi_{k}\big(2^n(x+1)-2j-1\big),\;\;\;k=0,\hdots,p,\;\;\;j=0,\hdots,2^n-1.
\end{equation}
As such, a function $u$ can be approximated in $V^p_n$ by
\begin{equation}
\label{eq:coarsefine}
u_h^{(n)}(x)=\left(\mathcal{P}^n u\right)(x)=\sum_{j=0}^{2^n-1}\sum_{k=0}^p s^n_{k,j}\phi^n_{k,j}(x),
\end{equation}
where
\begin{equation}
s^n_{k,j}=\langle u,\phi^n_{k,j} \rangle_{I^n_{j}},\;\;\; k=0,\hdots,p,\;\;\;j=0,\hdots,2^n-1.
\end{equation}
To transition from a coarse approximation space  $V_{n}^p$ to a finer approximation space $V_{n+1}^p$ {additional information is required}. This information is contained within the multiwavelet space $W_{n}^p$ defined to be the orthogonal compliment of $V_{n}^p$ in $V_{n+1}^p$:
\begin{equation}\label{eq:nest} 
V_{n+1}^p=V_{n}^p\oplus W_{n}^p,\hspace{1cm} W_n^p\subset V^p_{n+1},\hspace{1cm} W^p_n\perp V_n^p.
\end{equation}
Inductively, this hierarchy allows the expression of the fine mesh approximation space as the direct sum of the coarsest mesh approximation spaces and the intermediary wavelet spaces:
\[
     V_n^{p} = V_0^{p} \oplus W_0^{p} \oplus W_1^{p} \oplus \cdots \oplus W_{n-1}^{p}.
\]

An orthonormal basis for these multiwavelet spaces can be obtained via a Gram-Schmidt orthogonalization procedure and obtain uniqueness by imposing additional conditions \cite{Keinert2004}. Here the special case of Alpert's multiwavelets (see \cite{Alp93}) is considered, where Legendre polynomials are chosen as scaling functions,
\begin{equation}
    \{\phi_k(x)\}_{k=0}^p=\{\phi^0_{k,0}\}_{k=0}^p=
    \begin{cases}
    \sqrt{\frac{2k+1}{2}}P_k(x),&\quad x\in [-1,1]\\
    0,&\quad \text{Else}
    \end{cases}
\end{equation}
and the multiwavelet bases for the spaces $W^p_n$ satisfy the properties given below. Denoting the multiwavelets on the coarsest scale by 
\[\{\psi_k(x)\}_{k=0}^p=\{\psi^0_{k,0}\}_{k=0}^p,\]
 multiwavelet bases $\{\psi^{n}_{k,j}\}_{k,j}$ are obtained for each of the multiwavelet spaces $W^p_n$ by dilation and scaling where
\begin{equation}
\psi^{n}_{k,j}
=2^{n/2}\psi_k(2^n(x+1)-2j-1),\;\;\;k=0,\hdots,p,\;\;\;\quad j=0,\hdots,2^n-1. 
\end{equation}
These multiwavelets satisfy the following properties relevant to our analysis:
\begin{enumerate}
    \item Orthogonality conditions: $\langle \phi^n_{k,i},\psi^n_{\ell,j} \rangle_{I^n_j}=0$ and $\langle \psi^n_{k,i},\psi^n_{\ell,j} \rangle_{I^n_j}=\delta_{i,j}\delta_{k,
    \ell}$.
    \item Support conditions: $\text{supp}(\psi^n_{k,j})=\text{supp}(\phi^n_{k,j})=I^n_j$.
    \item Moment conditions: $\langle x^m,\psi^{n}_{k,j} \rangle_{I^n_j}=0$ for $m\leq k+p+1$. 
\end{enumerate}
We refer the interested reader to the original work of Alpert \cite{Alp93} for details as to their construction and additional properties satisfied by these piecewise polynomial multiwavelets, and to the thesis of Vuik \cite{Vuik2017} for explicit formulas for $\psi$ for the first few polynomials degrees. We note that by Equation \eqref{eq:nest} and the orthongonality conditions, the detail necessary for transition from a coarse mesh approximation space to a fine mesh approximation space
\begin{equation}
\sum_{j=0}^{2^n-1}d^n_{k,j}\psi^n_{k,j}(x)=\left(\mathcal{P}^{n+1}u\right)(x)-\left(\mathcal{P}^{n}u \right)(x),
\end{equation}
are obtained simply by projection of the initial data onto the multiwavelet space:
\[d^n_{k,j}=\langle u,\psi^n_{k,j}\rangle_{I^n_j}.\]
Hence, a fine mesh representation is simply a coarse mesh representation with additional multiwavelet details:
\begin{equation}
u_h^{(n+1)}(x)=\sum_{j=0}^{2^n-1}\sum_{k=0}^p s^n_{k,j}\phi^n_{k,j}(x)+\sum_{j=0}^{2^n-1}\sum_{k=0}^p d^n_{k,j}\psi^n_{k,j}(x).
\end{equation}
When considering the enhancement technique for adaptive PDEs in later sections, it is beneficial to compute the coarse mesh scaling and detail coefficients from the fine mesh scaling coefficients. To that end we define the quadrature mirror filters \cite{Vuik2017} useful in performing this exact task. Denote by $\mathbf{s}^n_j=[s^n_{0,j},\hdots,s^n_{p,j}]^T$ and $\mathbf{d}^n_j=[d^n_{0,j},\hdots,d^n_{p,j}]^T$, then 
\begin{align}\label{eq:scaling}
    \mathbf{s}^{n-1}_j&=H^{(0)}\mathbf{s}^{n}_{2j}+H^{(1)}\mathbf{s}^{n}_{2j+1}\\
    \mathbf{d}^{n-1}_j&=G^{(0)}\mathbf{s}^{n}_{2j}+G^{(1)}\mathbf{s}^{n}_{2j+1},\label{eq:details}
\end{align}
where $H^{(\;)}$ and $G^{(\;)}$ are $(p+1)\times(p+1)$ matrices with entries given by
\begin{align}
    H^{(0)}(\ell,r)=\langle \phi^0_{\ell,0},\phi^1_{r,0}\rangle_{I^1_0}, \;\;\;\;\;\;\;&\;\;\;\;\;\;\;H^{(1)}(\ell,r)=\langle \phi^0_{\ell,0},\phi^1_{r,1}\rangle_{I^1_1},\\
    G^{(0)}(\ell,r)=\langle \psi^0_{\ell,0},\phi^1_{r,0}\rangle_{I^1_0}, \;\;\;\;\;\;\;&\;\;\;\;\; \;\;G^{(1)}(\ell,r)=\langle \psi^0_{\ell,0},\phi^1_{r,1}\rangle_{I^1_1}.
\end{align}
One benefit of restricted uniform hierarchical meshes of the present analysis is that these matrices for transitioning detail and scaling coefficients between scales are fixed and only need be computed once. For more general meshes these quadrature mirror filters may be spatially varying, which increases the computational cost of the decomposition into multiwavelet and scaling coefficients.

\subsection{Extension to two dimensions}

The same techniques generalize to multiple dimensions easily via tensor-products of the approximation spaces and their associated bases. Here we follow the presentation of \cite{Bautista_paper}. To simplify notation, denote element indices on level $n$ by $\mathbf{j}=(j_x,j_y)$ for $j_x,j_y\in \{0,\hdots,2^{n-1}\}$ and modal indices by $\mathbf{k}=(k_x,k_y)$ for $k_x,k_y\in \{0,\hdots,p\}$. Let $\mathcal{T}^{(0)}=[-1,1]^2$ be the coarsest level mesh consisting of a single element, and $\mathcal{T}^{(n)}$ the mesh obtained by uniform subdivision into quadrants of each element of $\mathcal{T}^{(n-1)}$. The elements on mesh level $n$ are given by $Q^n_{\mathbf{j}}=I^n_{j_x}\otimes I^n_{j_y}$. Define the piecewise polynomial approximation space over mesh level $n$ to be
\begin{equation}
\mathbf{V}^p_n=\{v\;:v\big|_{I^n_j}\in\mathbb{Q}^p(Q^n_{\mathbf{j}}),\;Q^n_{\mathbf{j}} \in \mathcal{T}^{(n)}\}.
\end{equation}
A basis for this approximation space is given by $\{\mathbf{\phi}^n_{\mathbf{k},\mathbf{j}}\;:\;\norm{\mathbf{k}}_{\ell^{\infty}}\leq p,\;Q^n_{\mathbf{j}} \in \mathcal{T}^{(n)}\},$ where 
\begin{equation}
\mathbf{\phi}^n_{\mathbf{k},\mathbf{j}}(\mathbf{x})=\phi^n_{k_x,j_x}(x)\phi^n_{k_y,j_y}(y).
\end{equation}
The bases for the multiwavelet subspaces are obtained similar to the one-dimensional case except now there are three multiwavelet spaces denoted by $\mathbf{W}^p_{n,\alpha},\mathbf{W}^p_{n,\beta}$, and $\mathbf{W}^p_{n,\gamma}$. The bases for these subspaces are given by 
\begin{align}
    \{\mathbf{\psi}^{n,\alpha}_{\mathbf{k},\mathbf{j}}\;:\;\norm{\mathbf{k}}_{\ell^{\infty}}\leq p,\;Q^n_{\mathbf{j}} \in \mathcal{T}^{(n)}  \},\\
    \{\mathbf{\psi}^{n,\beta}_{\mathbf{k},\mathbf{j}}\;:\;\norm{\mathbf{k}}_{\ell^{\infty}}\leq p,\;Q^n_{\mathbf{j}} \in \mathcal{T}^{(n)}  \},\\
    \intertext{and}
    \{\mathbf{\psi}^{n,\gamma}_{\mathbf{k},\mathbf{j}}\;:\;\norm{\mathbf{k}}_{\ell^{\infty}}\leq p,\;Q^n_{\mathbf{j}} \in \mathcal{T}^{(n)}  \},
    \intertext{respectively, where}  \mathbf{\psi}^{n,\alpha}_{\mathbf{k},\mathbf{j}}(\mathbf{x})=\psi^n_{k_x,j_x}(x)\phi^n_{k_y,j_y}(y),\\
 \mathbf{\psi}^{n,\beta}_{\mathbf{k},\mathbf{j}}(\mathbf{x})=\phi^n_{k_x,j_x}(x)\psi^n_{k_y,j_y}(y),\\
    \intertext{and}   \mathbf{\psi}^{n,\gamma}_{\mathbf{k},\mathbf{j}}(\mathbf{x})=\psi^n_{k_x,j_x}(x)\psi^n_{k_y,j_y}(y).
\end{align}
Note that $W^p_{n,\alpha}$ which has multiwavelet functions for the x-component of the tensor-product basis can be thought of as corresponing to the x-direction, similarly, $W^p_{n,\beta}$ and $W^p_{n,\gamma}$ correspond to the $y$- and $xy$-directions respectively \cite{Vuik2017}. 

Similar to the one-dimensional case, the approximation in $\mathbf{V}^p_{n+1}$ is expressible as a sum of approximations in $\mathbf{V}^p_{n},\mathbf{W}^p_{n,\alpha},\mathbf{W}^p_{n,\beta}$, and $\mathbf{W}^p_{n,\gamma}$:
\begin{align}
\mathcal{P}^{n+1}u&=\sum_{Q^{n+1}_{\mathbf{j}} \in \mathcal{T}^{(n+1)}}\sum_{\norm{\mathbf{k}}_{\ell^{\infty}}\leq p} s^{n+1}_{\mathbf{k},\mathbf{j}}\mathbf{\phi}^{n+1}_{\mathbf{k},\mathbf{j}}\\
&=
\sum_{Q^n_{\mathbf{j}} \in \mathcal{T}^{(n)}}\sum_{\norm{\mathbf{k}}_{\ell^{\infty}}\leq p} s^{n}_{\mathbf{k},\mathbf{j}}\mathbf{\phi}^n_{\mathbf{k},\mathbf{j}}
+\\&\;\;\;\;
\sum_{Q^n_{\mathbf{j}} \in \mathcal{T}^{(n)}}\sum_{\norm{\mathbf{k}}_{\ell^{\infty}}\leq p} d^{n,\alpha}_{\mathbf{k},\mathbf{j}}\mathbf{\psi}^{n,\alpha}_{\mathbf{k},\mathbf{j}}+d^{n,\beta}_{\mathbf{k},\mathbf{j}}\mathbf{\psi}^{n,\beta}_{\mathbf{k},\mathbf{j}}+d^{n,\gamma}_{\mathbf{k},\mathbf{j}}\mathbf{\psi}^{n,\gamma}_{\mathbf{k},\mathbf{j}}.\nonumber
\end{align}
The quadrature mirror filters previously defined can also be used in higher dimensions to compute the coarse mesh scaling and multiwavelet coefficients from the fine mesh scaling coefficients \cite{Vuik2014}:
\begin{align}
s^{n}_{\mathbf{k},\mathbf{j}}=\sum_{\tilde{j}_x,\tilde{j}_y=0}^1\sum_{r_x,r_y=0}^pH^{(\tilde{j_x})}(k_x,r_x)H^{(\tilde{j_y})}(k_y,r_y)s^{n+1}_{\mathbf{r},2\mathbf{j}+\tilde{\mathbf{j}}},\\
d^{n,\alpha}_{\mathbf{k},\mathbf{j}}=\sum_{\tilde{j}_x,\tilde{j}_y=0}^1\sum_{r_x,r_y=0}^pG^{(\tilde{j}_x)}(k_x,r_x)H^{(\tilde{j}_y)}(k_y,r_y)s^{n+1}_{\mathbf{r},2\mathbf{j}+\tilde{\mathbf{j}}},\\
d^{n,\beta}_{\mathbf{k},\mathbf{j}}=\sum_{\tilde{j}_x,\tilde{j}_y=0}^1\sum_{r_x,r_y=0}^pH^{(\tilde{j}_x)}(k_x,r_x)G^{(\tilde{j}_y)}(k_y,r_y)s^{n+1}_{\mathbf{r},2\mathbf{j}+\tilde{\mathbf{j}}},\\
d^{n,\gamma}_{\mathbf{k},\mathbf{j}}=\sum_{\tilde{j}_x,\tilde{j}_y=0}^1\sum_{r_x,r_y=0}^pG^{(\tilde{j}_x)}(k_x,r_x)G^{(\tilde{j}_y)}(k_y,r_y)s^{n+1}_{\mathbf{r},2\mathbf{j}+\tilde{\mathbf{j}}}.
\end{align}

\subsection{Local vs. Global MRA}

We note for use in the adaptive mesh DG methods described later that the enhancement procedure does not require that we start from a mesh resultant from dyadic splitting. As a result, rather than start a global MRA described above, i.e. the mesh for $n=0$ is $[-1,1]^d$, we can instead choose to compute local multiresolution analyses (see \cite{Bautista_paper} for details) where we initialize a level $n=0$ multiresolution analysis for each interval or quadrilateral element $\Omega_e$ of some initial coarse mesh. We describe the idea in one dimension.  Let $\xi^{(e)}=\frac{2}{|\Omega_e|}(x-x_e)$, where $x_e$ is the center of the interval $\Omega_e$, denote the reference mapping $\xi^{(e)}(\Omega_e)\rightarrow[-1,1]$. Choose scaling functions $\tilde{\phi}^0_{k,0}=\sqrt{\frac{2}{|\Omega_e|}}\phi_k(\xi^{(e)}(x))$ and in turn multiwavelets  $\tilde{\psi}^0_{k,0}=\sqrt{\frac{2}{|\Omega_e|}}\psi_k(\xi^{(e)}(x))$. The multiresolution analyses generated by these choices of scaling and multiwavelet functions satisfy the same orthogonality properties described previously except $I^n_j$ now represents the $j$th interval in the $n$-times repeated dyadic subdivision of $\Omega^{(e)}$. We use the tilde notation $\tilde{d}^n_{k,j}$ and $\tilde{s}^n_{k,j}$ to denote the detail and scaling coefficients resultant from this local subdivision. For more details, see \cite{Bautista_paper}.

\subsection{Enhancement procedure}

To enhance the MRA procedure for the construction of a fine mesh approximation $\tilde{u}^{(n+1)}_h$, the coarse mesh data,  $u^{(n)}_h$, is filtered by convolving it with a specially chosen kernel function, $K_H$, before projecting onto a refined mesh. Mathematically, we write the enhanced fine mesh approximation as
\begin{equation}\label{eq:recon} 
\tilde{u}^{(n+1)}_h=\mathcal{P}^{n+1}(K_H\star u^{(n)}_h). 
\end{equation}
 where $\mathcal{P}^{n+1}$ represents the projection operator projecting the filtered approximation onto the fine mesh approximation space.  This allows for computing the multiwavelet detail coefficients from the enhanced fine mesh approximation as shown in \cite{SIACMRA,picklo}.  In the following, where it does not impact clarity, we will drop the mesh level indicator $(n)$ in the superscripts for readability.  This is an alternative to the work of Bautista et al., who considered a similar refinement.  In \cite{Bautista}, the fine mesh approximation came from the projection of a least-squares polynomial formed from the coarse mesh approximation. In section \ref{sec:adaptive_DG} we show how fine mesh data can be leveraged in adaptive DG schemes, but first  we provide a review of fundamental concepts used in the enhancement procedure.

\subsection{SIAC filter}

Our enhancement procedure makes use of the Smoothness-Increasing Accuracy-Conserving (SIAC) filters that were originally developed in \cite{Bramble1977} for achieving superconvergence in finite element data. SIAC extensions include position-dependent kernels for non-periodic boundaries \cite{ryan2003,SRV,XLiOne}, Hexagonal filters for hexagonal meshes \cite{HexSIAC}, and Line-SIAC (L-SIAC) filters for computationally efficient filtering of higher dimensional data \cite{LSIAC}. More recently, SIAC filters have been used as an enhancement tool for enriching data in a multiresolution framework \cite{picklo,SIACMRA}. For the purposes of this paper the focus is on these last two ideas. 

The one-dimensional SIAC filtered approximation is given by
\begin{equation}
u^{\star}_h(x)=\int_{-\infty}^{\infty}K^{(r+1,\ell)}_H(y-x)u_h(y)\;dy,
\end{equation}
 where $u_h(x)$ is the original approximation and the scaled SIAC kernel itself is defined by the relation
\begin{equation}\label{eq:kernel}
K^{(r+1,\ell)}_H(t)=\frac{1}{H}\sum_{\gamma=-r}^rc_{\gamma}B^{\ell}\Big(\frac{t}{H}-x_\gamma\Big),
\end{equation}
where $B^{\ell}$ are central B-splines of order $\ell$ defined by the recursion relation (see de Boor \cite{deBoor})
\begin{align}
    B^1(t)&=\chi_{[-1/2,1/2)}(t)\\
    B^{\ell+1}(t)&=\frac{1}{\ell}\Big[\Big(\frac{\ell+1}{2}+t\Big)B^{\ell}(t+1/2)+\Big(\frac{\ell+1}{2}-t\Big)B^{\ell}(t-1/2)\Big],\;\;\;\ell>1.
\end{align}
The coefficients $c_{\gamma}$ are chosen s.t. the kernel reproduces polynomials up through degree $r$ under convolution:
\begin{equation}\label{eq:coeffs}
K^{(r+1,\ell)}\star x^{m}=x^m,\;\;\;m=0,\hdots,r.
\end{equation}
In the context of multiresolution analysis, we consider symmetric kernels with $x_\gamma = \gamma$. A study of the linear system resulting from Equation \ref{eq:coeffs} has been conducted for kernels whose splines have uniform knot sequences in \cite{Mirzargar2016}, and more general knot sequences formulations in \cite{peters2015,peters2016}. For the purposes of our study, we choose to use $2p+1$ B-splines of order $1$ in constructing the filter, where $p$ is the polynomial degree of the underlying data on a given element. It was observed in \cite{SIACMRA} that the choice of $r+1=2p+1$ and $\ell=1$ led to the greatest error reduction under the enhancement procedure, though it did not lead to the superconvergent phenomena associated with a single filter step using $2p+1$ B-splines of order $p+1$.

In two dimensions we consider the L-SIAC filter which maintains the one-dimensional support of the original kernel, and thereby avoids the costly multi-dimensional quadrature associated with tensor-product formulations. Given initial data $u_h(\mathbf{x})$, we post-process the data by convolution along a parameterized line $\Gamma(t)$, obtaining a filtered approximation:
\begin{equation}
u^{\star}_h(\mathbf{x})=\int_{-\infty}^{\infty}K^{(r+1,\ell)}_H(t)u_h(\Gamma(t))\;dt,
\end{equation}
where 
\begin{equation}
\Gamma(t)=\mathbf{x}+t(\cos(\theta),\sin(\theta)).
\end{equation}
Here the kernel function is given by Equation \eqref{eq:kernel}, where $\theta$ denotes the orientation of the line over which the kernel is supported, and $H$ denotes the kernel scaling, which is usually associate with the mesh size. For uniform quadrilateral meshes the choices  $\theta=\pi/4 = \arctan\left(\frac{\Delta y}{\Delta x}\right)$ and $H=\sqrt{2}h$ are typically made to ensure error reduction. For  unstructured meshes, no general procedure for choosing the orientation has been introduced.  Frequently,  the maximum edge length of the mesh is chosen for the scaling, though adaptive scalings have also been considered \cite{Jallepalli2019B}. In this paper, the orientation is taken to be $\theta=\pi/4$ and we investigate setting the scaling, $H$, to be the maximum edge length as well as the adaptive scaling (see section \ref{sec:adap_scaling}). 

\subsection{$L^2$-projection} 
\label{sec:L-2proj}

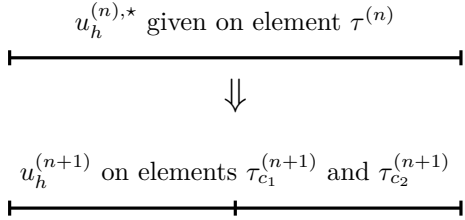
\begin{wrapfigure}[]{l}{0.4\textwidth}
\begin{tikzpicture}
\draw[very thick] (-3,0)--(3,0);
\draw[very thick]  (-3,-0.1)--(-3,0.1);
\draw[very thick]  (3,-0.1)--(3,0.1);
\draw (0,0.5) node{ $u^{(n),\star}_h$ given on element $\tau^{(n)}$};

\draw[very thick]  (0,-0.5) node{{\Large $\Downarrow$}};

\draw[very thick]  (-3,-2)--(3,-2);
\draw[very thick]  (-3,-2.1)--(-3,-1.9);
\draw[very thick]  (3,-2.1)--(3,-1.9);
\draw[very thick]  (0,-2.1)--(0,-1.9);
\draw (0,-1.5) node{ $u^{(n+1)}_h$ on elements $\tau^{(n+1)}_{c_1}$ and $\tau^{(n+1)}_{c_2}$};

\end{tikzpicture}
    \caption{Illustration of enhancement}
    \label{fig:1D procedure}
\end{wrapfigure}
For the purpose of refinement in the enhancement procedure, the filtered approximation is projected onto a refined mesh. This procedure is similar to that in the initialization Equation \eqref{eq:L2}. Assuming $\tau^{(n)}$ is the coarse mesh element over which $u^{(n),\star}_h$ is defined and $\tau^{(n+1)}_{c}$ is a sub-element or ``child" element of $\tau^{(n)}$ which exists on the finer mesh, the enhanced approximation on that sub-element is computed by requiring:
\begin{equation}
    \langle u^{(n),\star}_h,\phi_{m,\tau^{(n+1)}_{c}}\rangle_{\tau^{(n+1)}_{c}}=\langle u_h^{(n+1)},\phi_{m,\tau^{(n+1)}_{c}}\rangle_{\tau^{(n+1)}_{c}},
\end{equation}
for each $m\in \mathcal{I}$ (see Figure \ref{fig:1D procedure}). Repeating this procedure for each $\tau^{(n)}\in\mathcal{T}^{(n)}$ and $\tau^{(n+1)}_c\subset \tau^{(n)}$ in $\mathcal{T}^{(n+1)}$ a polynomial representation is obtained on the refined mesh. Note that we drop the star superscript from the enhanced approximation as the approximation now belongs to the fine mesh approximation space over $\mathcal{T}^{(n+1)}$.

\subsection{Adaptive kernel scaling}\label{sec:adap_scaling}

In considering meshes where elements sizes vary significantly we investigate the impact of the adaptive kernel scaling introduced by Jallepalli \cite{Jallepalli2019B} on the enhancement procedure. The computation of this adaptive scaling consists of a pre-processing step assigning a characteristic length to each vertex in the mesh, and a post-processing step where based off the location of the filtering point in an element, a weighting of the characteristic lengths of the vertices in that element produces a characteristic length for that filtering point. Here we focus on the two-dimensional case.

Let $v_i$ be the $i$-th vertex in the mesh, $\tau_j$ an element with vertex $v_i$, and $L_{\tau_j}$ the length of the edge of $\tau_j$ with maximum length. Then a characteristic length $H(v_i)$ is assigned to vertex $v_i$ via the expression:
\begin{equation}
H(v_i)=\frac{\sum_{j}L_{\tau_j}\cdot Area(\tau_j)}{\sum_{j}Area(\tau_j)},
\end{equation}
where $Area(\tau_j)$ is the area of element $\tau_j$. Computing these scaling for each vertex in the mesh constitutes the pre-processing step. Then, given a filtering point $x^{\star}$, the kernel scaling used at that location is computed via the relation
\begin{equation}
H(x^{\star})=\sum_{\ell=0}^N\lambda_{\ell}H(v_{\ell})
\end{equation}
where $v_{\ell}$ are the vertices of the element containing $x^{\star}$, $N+1$ is the number of vertices of that element, and $\lambda_{\ell}$ are the barycentric coordinates of $x^{\star}$ in that element.

\subsection{Multiwavelet details from (L)SIAC-MRA}
Now that we have defined our enhancement procedure, we can define the multi-wavelet detail coefficients of the enhanced approximation. As displayed in equation \ref{eq:recon}, application of the (L)SIAC-MRA procedure provides a modal representation on a mesh one level finer than the previously given data. Therefore, in one dimension the enhanced reconstruction can be expressed as
\begin{equation}
\tilde{u}^{(n+1)}_h(x)=\sum_{j=0}^{2^{n+1}-1}\sum_{k=0}^p\tilde{s}^{n+1}_{k,j}\phi^{n+1}_{k,j}(x).
\end{equation}
Application of the quadrature mirror filters yields
\begin{equation}
\tilde{u}^{(n+1)}_h(x)=\sum_{j=0}^{2^{n}-1}\sum_{k=0}^p\tilde{s}^{n}_{k,j}\phi^{n}_{k,j}(x)+\sum_{j=0}^{2^{n}-1}\sum_{k=0}^p\tilde{d}^{n}_{k,j}\psi^{n}_{k,j}(x),
\end{equation}
where $\tilde{\mathbf{s}}^{n}$ and $\tilde{\mathbf{d}}^{n}$ are obtained by applying equations \eqref{eq:scaling} and \eqref{eq:details} to $\tilde{\mathbf{s}}^{n+1}$. It is exactly these multiwavelet details that transition {the} approximation from a coarse to fine mesh that are used for adaptivity indication in Section \ref{sec:adaptive_DG}. In two dimensions the same procedure is performed, with the only difference being the additional terms resulting from multiple multiwavelet subspaces.

%% file: sections/Mesh_Transition.tex
\section{Application to Nonuniform Mesh Transition} \label{sec:Mesh_Transition}

Now that the components necessary for enhancement on an unstructured mesh have been defined, we can investigate the error-reducing capabilities of the LSIAC-MRA procedure when applied to a variety of mesh and initial condition pairings. Considered here are Delaunay triangulations and perturbed quadrilateral meshes. For the purposes of assessing the effects of adaptive scaling on the enhanced reconstruction procedure, the following two meshes (depicted in Figure \ref{fig:adaptive_mesh_examples}) generated by using the GMSH application \cite{gmsh} are also considered. These meshes are denoted here by SVR100 and ISR2. The SVR100 mesh is smoothly varying where the ratio of the boundary to center element is 100, while the ISR2 mesh possesses a more abrupt transition between boundary elements and interior elements with an element size ratio of 2. For ease of implementation, it is enforced in all cases that the meshes contain no hanging nodes, but it is stressed that this is not a requirement for the method itself.
\input{Figures/adap_scaling/adaptive_mesh_examples}

\subsection{Mesh Refinement}
For refinement consider a uniform subdivision into $2^d$ elements. The subdivision is interpolatory in that nodes on the coarse mesh remain present after refinement. For every edge on the mesh, refine by inserting a node at the midpoint of that edge (in the quadrilateral case a node is also inserted at the average of the midpoints of the edges of each element). These new nodes are then connected with the original nodes to form four new elements on the refined mesh. A depiction of this subdivision procedure is shown in Figure \ref{fig:subdivision} and a comparison of initial and subdivided meshes in Figure \ref{fig:mesh_examples}.
\input{Figures/meshes/subdivision.tex}

\input{Figures/meshes/mesh_examples.tex}

\subsection{Numerical tests}
We consider the periodic domain $\Omega=[0,1]^2$ and meshes of $N$ elements and their refinements generated according to procedures described in Section \ref{sec:enhancement_procedure}. We note that the experiments in this section do not require the use of multiwavelet bases, which are typically expensive to compute over nonuniform meshes (see \cite{Gerhard2022}). In the following we consider $\mathbb{Q}^p$ polynomial bases for the quadrilateral meshes, and $\mathbb{P}^p$ polynomial bases for the triangular meshes. We compute the $L^2$- and $L^{\infty}$-errors via Gauss-Legendre quadrature at $6^d$ quadrature nodes on each element relative to the finest mesh scale considered. These error expressions are given by
\begin{equation}
    \norm{u-u_{h}}_{L^2(\Omega)}=\sqrt{\int_{\Omega}|u-u_{h}|^2\;d\Omega}
\qquad 
\text{ and }
\qquad
    \norm{u-u_{h}}_{\infty}=\sup_{x\in\Omega}|u-u_{h}|.
\end{equation}

\subsection{Impact of characteristic frequencies and lengths on error reduction}

Below we investigate the error reduction capabilities of the enhanced reconstruction procedure for initial data of varying characteristic frequencies. As initial data, we consider the $L^2$-projections of the tensor-products of the following functions over the domain $[0,1]^2$, i.e. $u^{(n)}_h=\mathcal{P}^{(n)}u$, $u(x,y)=f_i(x)f_i(y)$ for:
\begin{itemize}
    \item $f_1(x)=\sin(2\pi k x)$,\\
    \item $f_2(x)=\exp\Big(\frac{-(x-1/2)^2}{2(kc)^2}\Big)$,\\
    \item $f_3(x)=\frac{1}{2}\Big[\tanh\Big(\frac{x-7/20}{kw}\Big)-\tanh\Big(\frac{x-13/20}{kw}\Big)\Big].$
\end{itemize}
The first function represents a sine wave of varying frequencies $(2\pi k)$, the second a Gaussian function with varying standard deviation $(kc)$, and the last a double hyperbolic tangent function with transition scaling $(kw)$. Here $w$ and $c$ are chosen such that function $f_2$ and $f_3$ are less than $10^{-16}$ at the domain boundaries. The parameter $k$ is a scaling parameter that is varied to help understand the effects of initial data with different characteristic lengths on the LSIAC-MRA procedure. These test functions each represent a fundamental feature we might expect to encounter in physical simulations of fluid flow: sines representing a standard basis for the expression of waves, Gaussian representing vortices, and tanh functions possessing sharp gradients which are commonplace in PDEs with hyperbolic attributes. Furthermore, we note that these test functions are meant to be challenging functions for the enhancement procedure. Specifically, they either contain sharp gradients, or high-frequency oscillations relative to the scaling of the kernel. By increasing the frequency of the sines or the transition width of the sharp data, a better sense of the applicability of multi-element SIAC enhancements can be gained

Considering first the sine initial data, we observe in Figure \ref{fig:tensor sine} that errors are reduced for all polynomial degrees for $k=1$. Doubling the frequency of the initial data, we obtain error reduction up through $p=2$, with the $p=3$ approximation, the $L^2$-errors flatten for $k=2$. Similarly doubling the frequency to $k=4$ yields error reduction up through $p=1$ and for $k=8$ error reduction is obtained only for $p=0$. Similar results are obtained for both perturbed quadrilateral meshes and Delaunay triangulations. Increasing the initial mesh resolution would increase the points per wavelength (ppw) in a single direction and thereby remedy this issue for higher orders.  However, we note that, based off results from the one dimensional SIAC filter, higher frequency waves are damped more.  Hence, a tuning of kernel scalings should be considered when applying this technique and a more in-depth investigation will be performed in future work. Another interpretation is that, at least empirically, the minimum resolution for improvement increases by a factor of two with each additional polynomial degree. The benefit of the procedure as an error-reducing reconstruction on fine meshes must then be tempered by the expectation of high-frequency data.  

\input{Figures/error_plots/tensor_product_sine.tex}

Next, we consider the impact of characteristic lengths of the Gaussian and double hyperbolic tangent initial conditions. Note here that as the scaling, $k$, decreases, the Gaussian function becomes sharper and more closely resembles a point-source centered at $(1/2,1/2)$, whereas the double hyperbolic tangent function will more closely resemble the indicator function over the region $[7/20,13/20]^2$. In Figure \ref{fig:tanh_and_gaussian} we depict the {\it global} error effects of repeated application of the LSIAC-MRA reconstruction where the initial data has varying characteristic scales governed by the choice of the scaling factor $k$.  For the double hyperbolic tangent function we observe that for small $p$ and large $k$ we see global error reduction with enhancement, for $k<1/2$ we no longer see improvement, and for $p=2$ no improvement is observed. This is because the sharp gradient is highly localized, and becomes sharper with decreasing $k$. A weighted average over increasingly larger intervals as $p$ increases does not improve the ability of the method to capture this feature.  The results would also improve if an area around the discontinuous region were excluded from the error calculation.  We note that, similar to the sine case, adaptive scaling relative to the feature being captured may improve these results.

The enhancement procedure performs better for the global errors when Gaussian initial data is considered.  In this case, we observe improvement for all scalings depicted when $p=0$. For $p=1$ we observe improvement in the global $L^2-$errors in all cases, but only for the first two scalings with respect to the global $L^{\infty}-$error. When $p=2$ we only observe improvement for the largest scaling. Similar to the hyperbolic tangent case though less pronounced, for the global errors to reduce, this would require high levels of resolution to resolve without changing the kernel formulation. In this sense, a kernel scaling implicitly informed by the underlying data seems necessary and will be a study of future work.

\input{Figures/error_plots/tanh_and_gaussian}

\subsection{Adaptive versus constant scaling}
We now consider the filter scaling for non-uniform meshes and how it effects the performance of LSIAC-MRA procedure.  Specifically, we compare the traditional maximum edge length constant scaling versus the adaptive kernel scaling of Jallepalli \cite{Jallepalli2019B}. As depicted in Figure \ref{fig:adapt_vs_constant_scaling} we observe that the adaptively scaled kernel performs the same or better on variable element size meshes. Of particular interest is the $p=2$ and $p=3$ cases for the smoothly varying meshes. In this case, improvement is observed for the adaptively scaled procedure, but not in the constantly scaled case. Though we do not pursue a computational performance study here, we note that the adaptively scaled procedure is significantly faster than the constant scaling approach. This speed up of the filtering procedure by using adaptive scalings has been investigated in \cite{Jallepalli2019B}.

\input{Figures/error_plots/Adaptive_vs_constant/adaptive_vs_constant}

%% file: Figures/adap_scaling/adaptive_mesh_examples.tex
\begin{figure}
    \centering
    \begin{tabular}{c c} 
        \multicolumn{2}{c}{\textbf{Meshes and Scaling Contours}}\\       
           %\textbf{SVR10}& 
           \textbf{SVR100}&  \textbf{ISR2}  \\ 
 \includegraphics[width=0.45\linewidth]{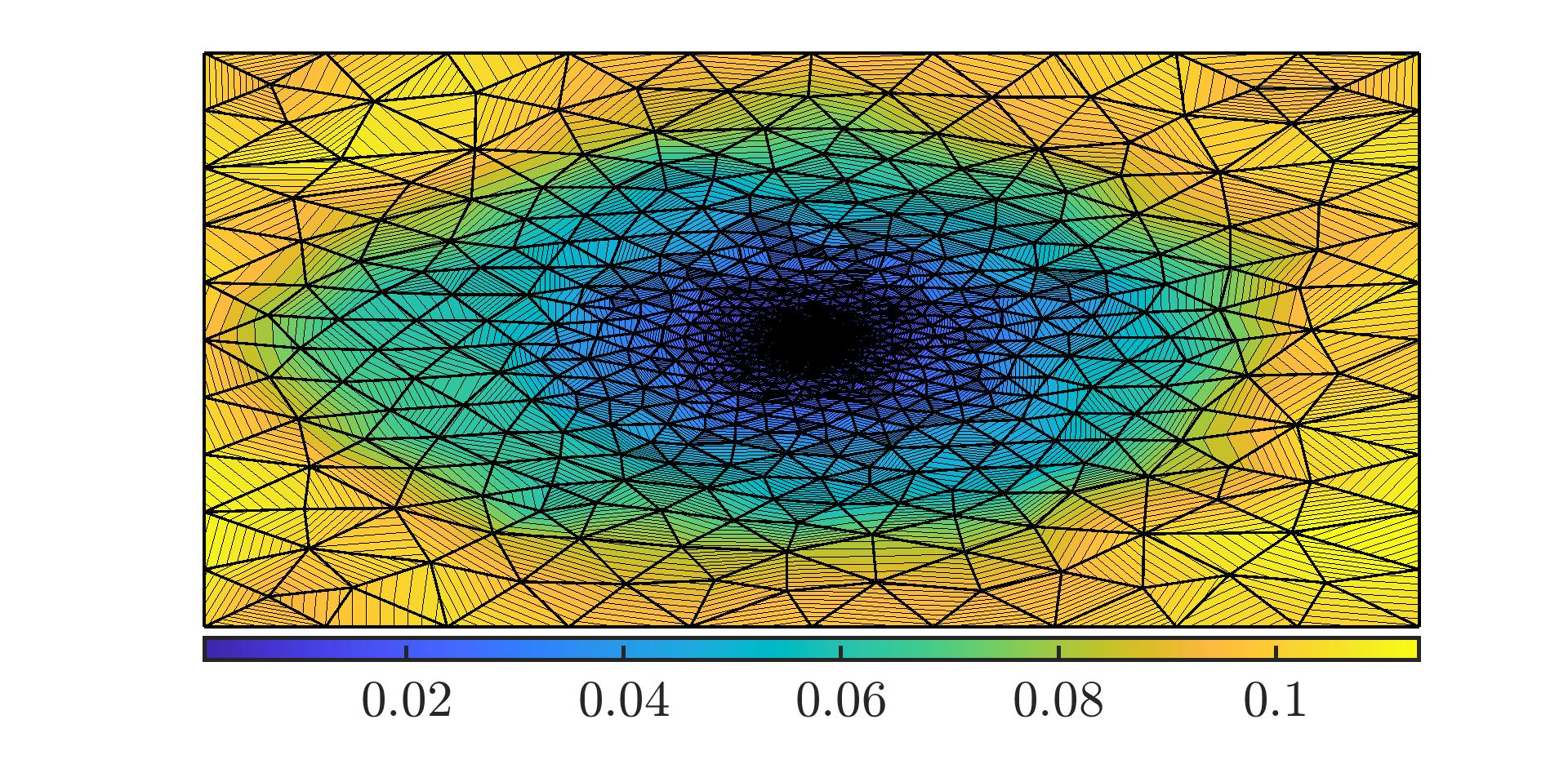} & \includegraphics[width=0.45\linewidth]{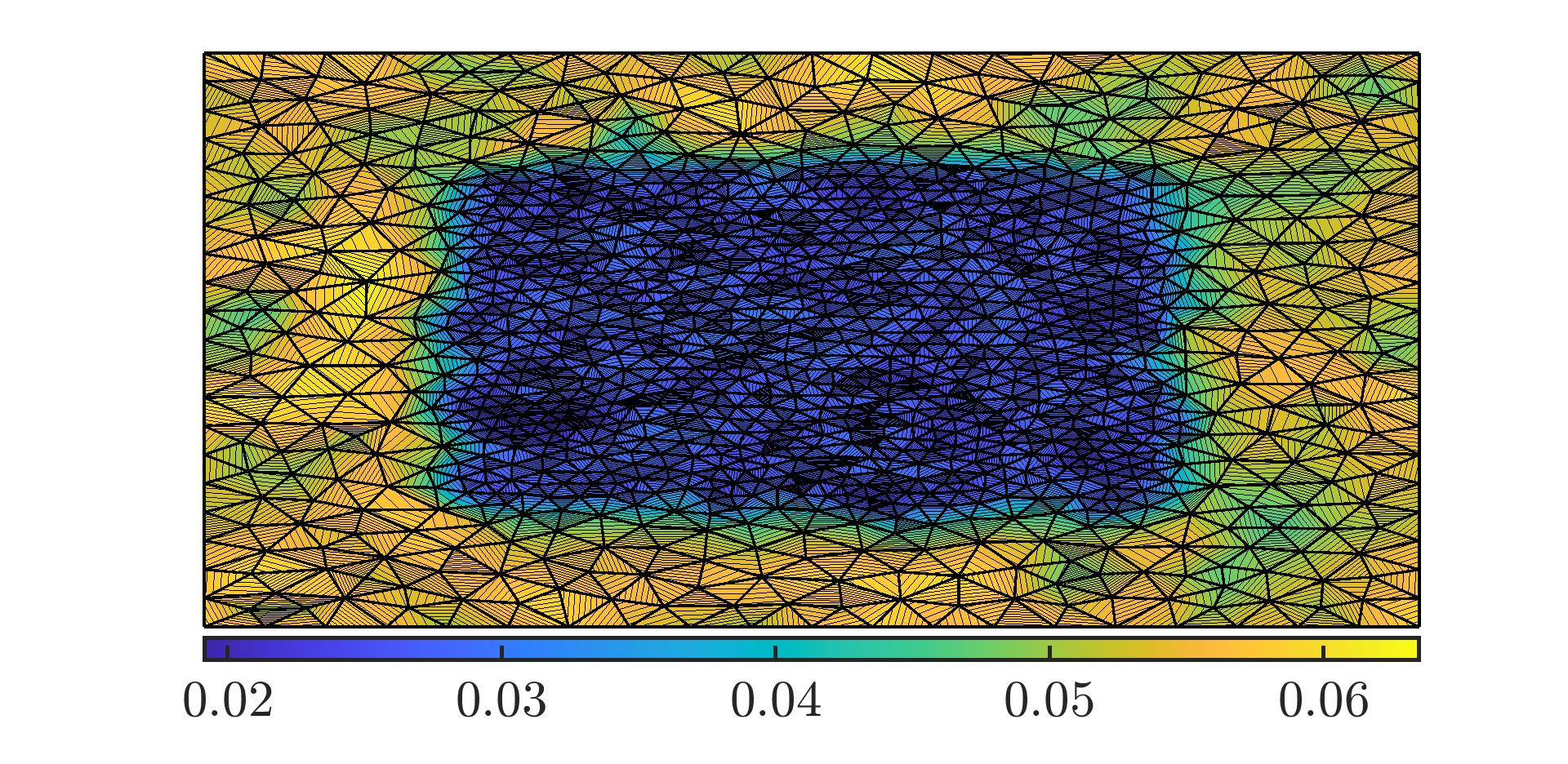}\\ 
    \end{tabular}
  \caption{Examples of nonuniform meshes with varying ratios of boundary element to interior element sizes. The color map indicates the kernel scaling used at each point in the domain.}
  \label{fig:adaptive_mesh_examples}
\end{figure}

%% file: Figures/meshes/subdivision.tex
  \begin{figure}[tp!]
  \centering
    
\begin{tabular}{c}

\begin{tikzpicture}[scale=0.75] \usetikzlibrary{calc}

 \coordinate (Origin) at (4,0);

\draw[black] (-1,-1)--(-1.1,1.25)--(0.9,1)--(1.1,-0.7)--(-1,-1);

\draw[black,fill] (-1,-1) circle (1.5pt);
\draw[black,fill] (-1.1,1.25) circle (1.5pt);
\draw[black,fill] (0.9,1) circle (1.5pt);
\draw[black,fill] (1.1,-0.7) circle (1.5pt);

\draw[black] (2,-1)--(1.9,1.25)--(3.9,1)--(4.1,-0.7)--(2,-1);

\draw[black,fill] (2,-1) circle (1.5pt);
\draw[black,fill] (1.9,1.25) circle (1.5pt);
\draw[black,fill] (3.9,1) circle (1.5pt);
\draw[black,fill] (4.1,-0.7) circle (1.5pt);

\draw[blue,fill] (1.95,0.125) circle (1.5pt);
\draw[blue,fill] (2.9,1.125) circle (1.5pt);
\draw[blue,fill] (4,0.15) circle (1.5pt);
\draw[blue,fill] (3.05,-0.85) circle (1.5pt);

\draw[blue,fill] (2.975,0.1375) circle (1.5pt);
\draw[blue,dashed] (1.95,0.125)--(4,0.15); 
\draw[blue,dashed] (3.05,-0.85)--(2.9,1.125); 
\draw [-stealth, line width=1mm] (1.2,0)--(1.75,0);

\draw[black,fill] (6,-1) circle (1.5pt);
\draw[black,fill] (6.5,1.25) circle (1.5pt);
\draw[black,fill] (8,-0.5) circle (1.5pt);

\draw[black] (6,-1)--(6.5,1.25)--(8,-0.5)--(6,-1);

\draw[black,fill] (9,-1) circle (1.5pt);
\draw[black,fill] (9.5,1.25) circle (1.5pt);
\draw[black,fill] (11,-0.5) circle (1.5pt);

\draw[black] (9,-1)--(9.5,1.25)--(11,-0.5)--(9,-1);

\draw [-stealth, line width=1mm] (8.2,0)--(9,0);

\draw[blue,fill] (9.25,0.125) circle (1.5pt);
\draw[blue,fill] (10.25,0.375) circle (1.5pt);
\draw[blue,fill] (10,-0.75) circle (1.5pt);

\draw[blue,dashed] (9.25,0.125)--(10.25,0.375); 
\draw[blue,dashed] (9.25,0.125)--(10,-0.75); 
\draw[blue,dashed] (10,-0.75)--(10.25,0.375);

 \end{tikzpicture}

  \end{tabular}

  \caption{Refinement procedure for quadrilateral and triangular elements.}

      \label{fig:subdivision}
\end{figure}
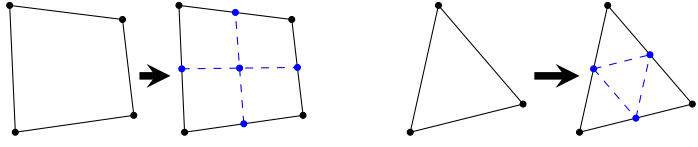

%% file: Figures/meshes/mesh_examples.tex
\begin{figure}[tp!]
    \centering

\begin{tabular}{c c}
 \multicolumn{1}{c}{\textbf{Initial Mesh}}  &  \multicolumn{1}{c}{\textbf{Refined Mesh}} \\
        \includegraphics[width=0.35\linewidth]{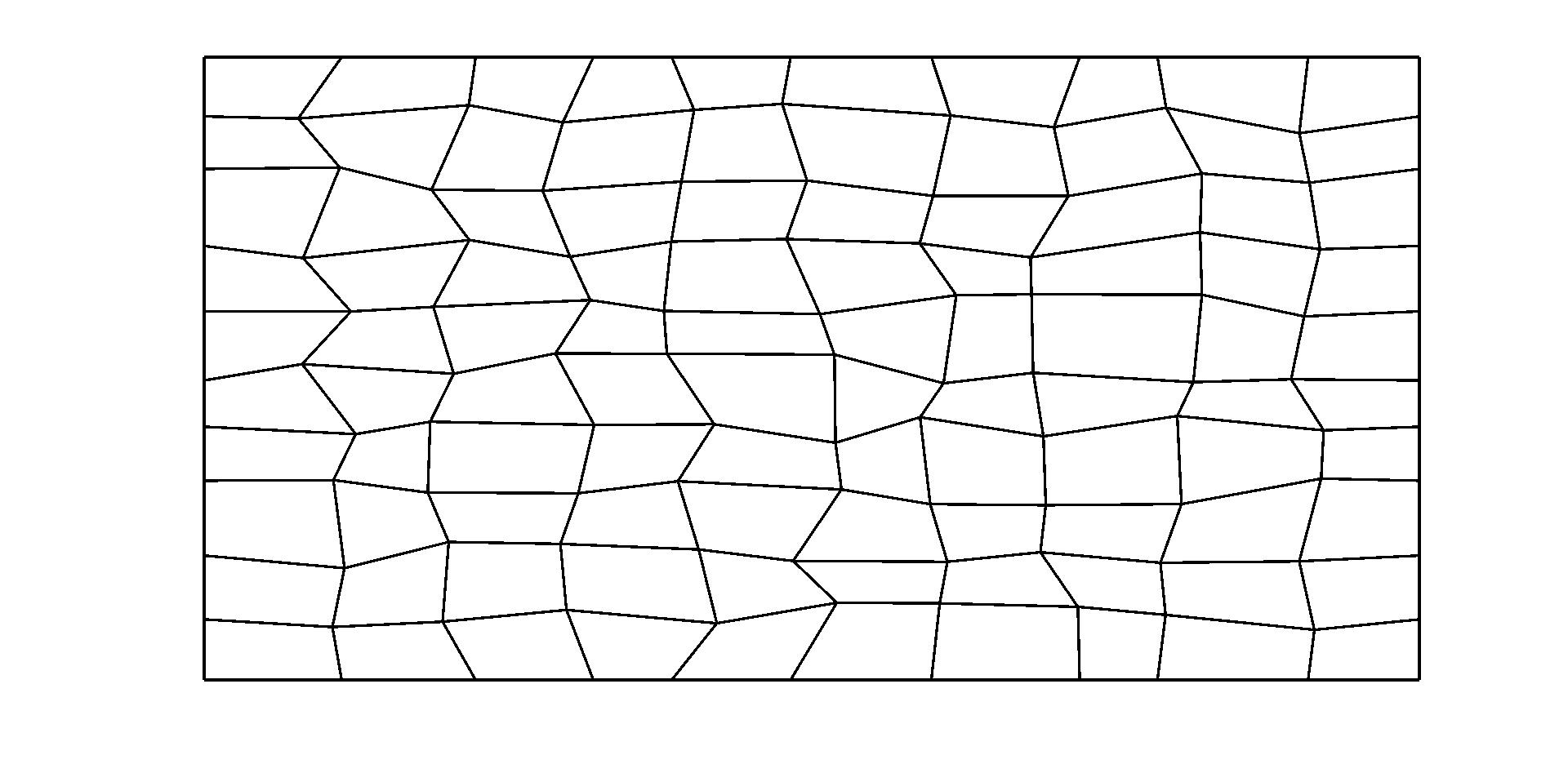}&
                \includegraphics[width=0.35\linewidth]{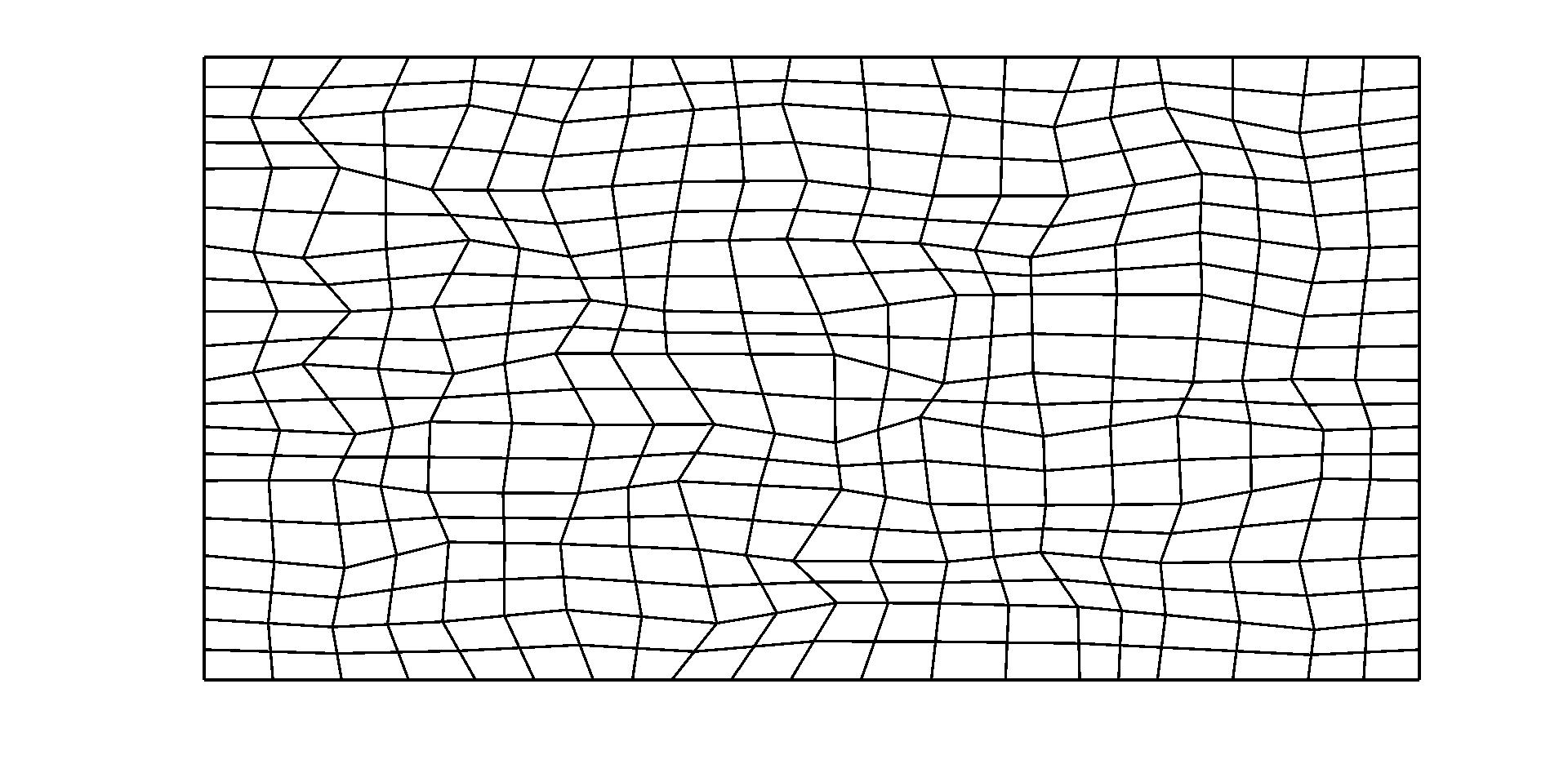} 
\\
        \includegraphics[width=0.35\linewidth]{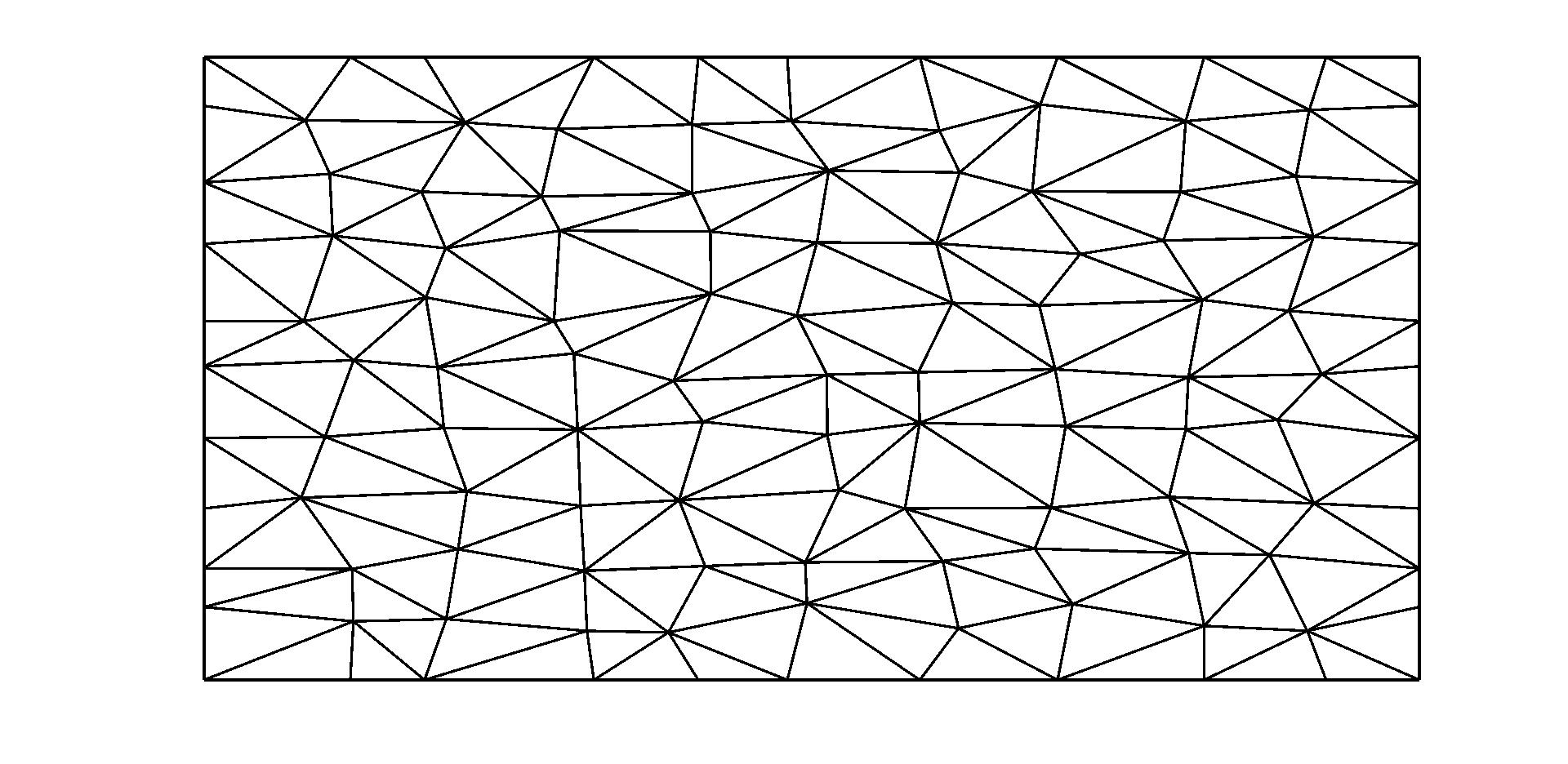}&
                \includegraphics[width=0.35\linewidth]{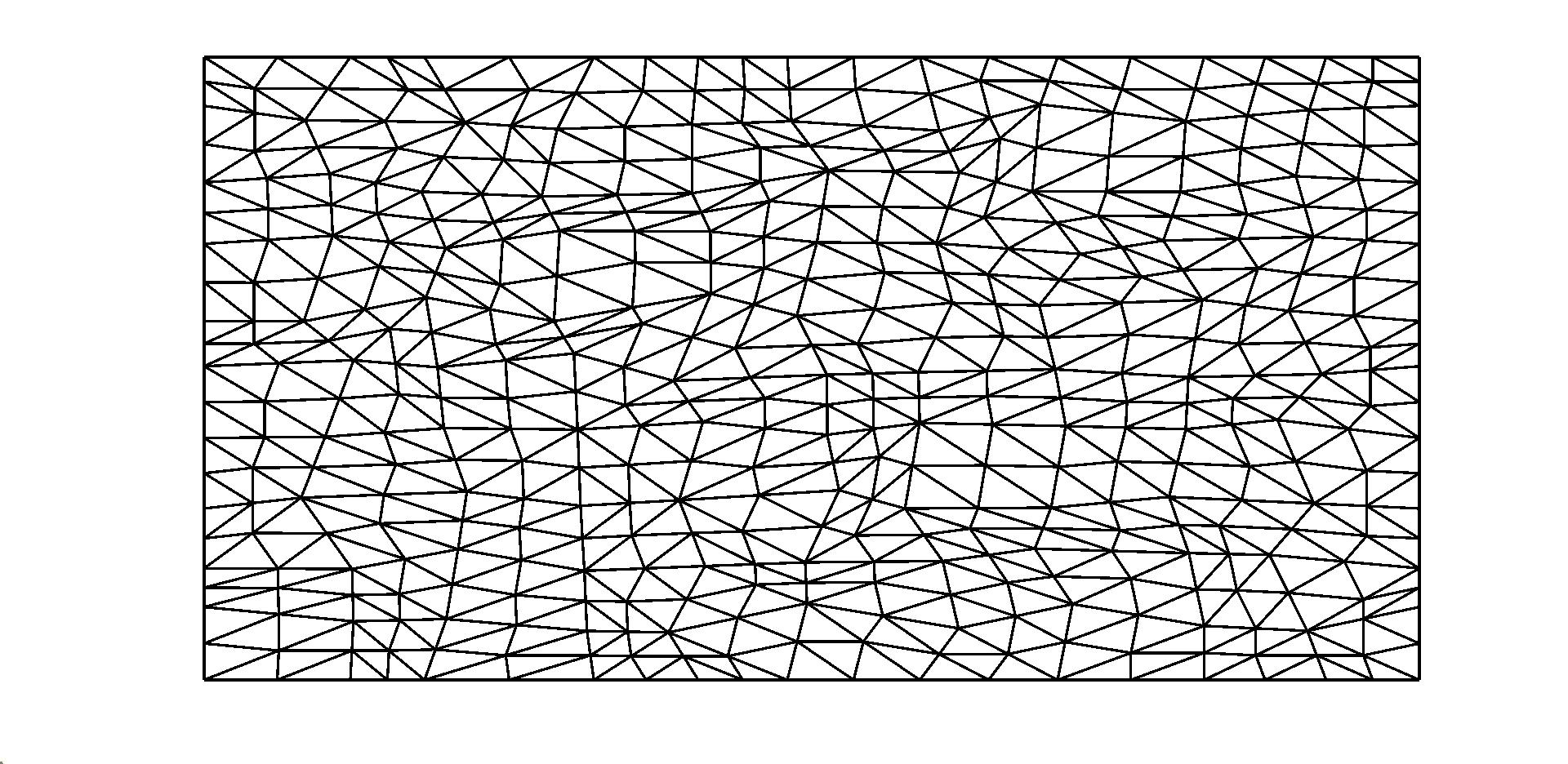} 
\end{tabular}

    \caption{Depiction of nonuniform quadrilateral and triangular meshes (left) and the refined meshes resulting from a uniform subdivision (right).}
    \label{fig:mesh_examples}
\end{figure}

%% file: Figures/error_plots/tensor_product_sine.tex
\begin{figure}[hp!]
    \centering

\begin{tabular}{c }
\multicolumn{1}{c}{$\mathbf{40\times 40}$ \textbf{Perturbed Quadrilateral Mesh}}\\
       \includegraphics[width=0.8\linewidth]{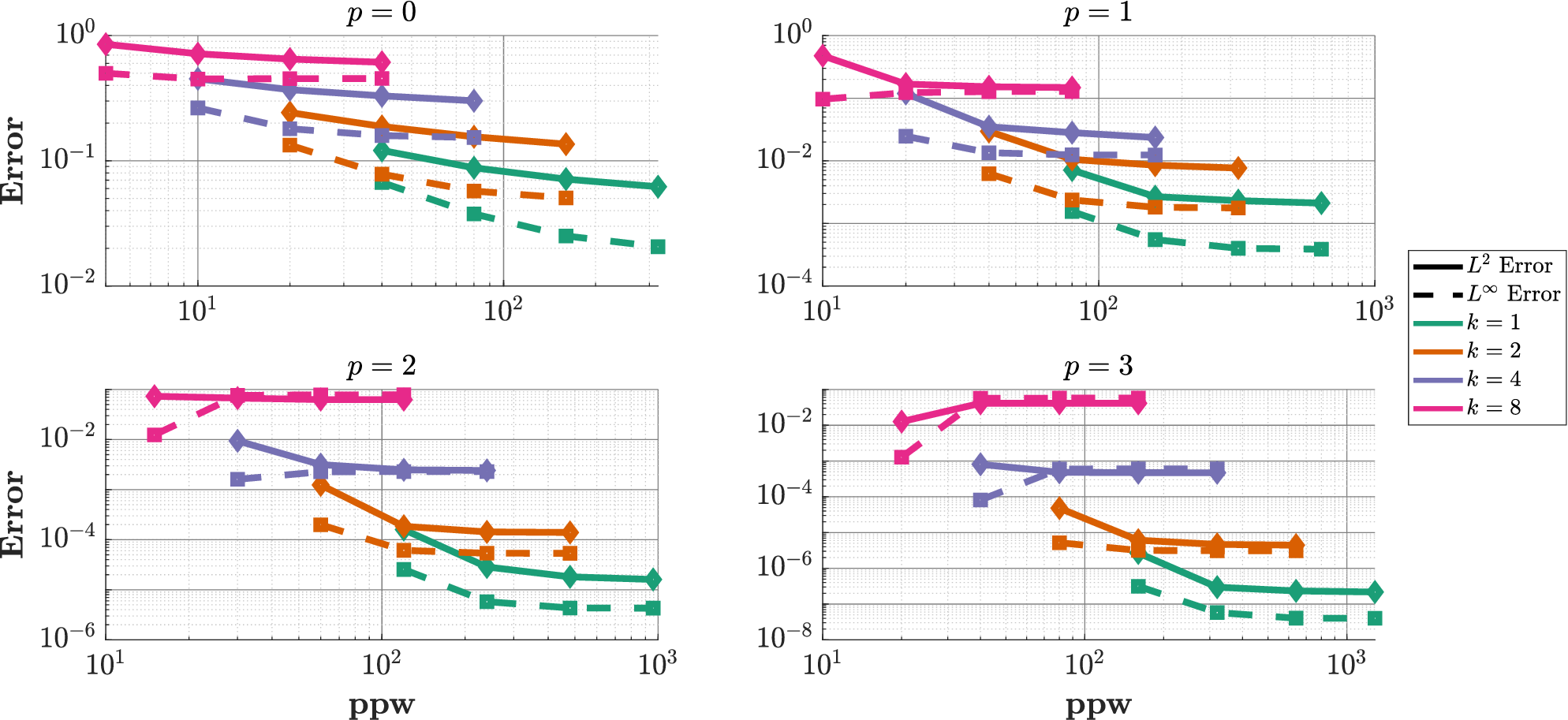}\\
       \multicolumn{1}{c}{$\mathbf{2\times40\times 40}$ \textbf{Delaunay Mesh}}\\
       \includegraphics[width=0.8\linewidth]{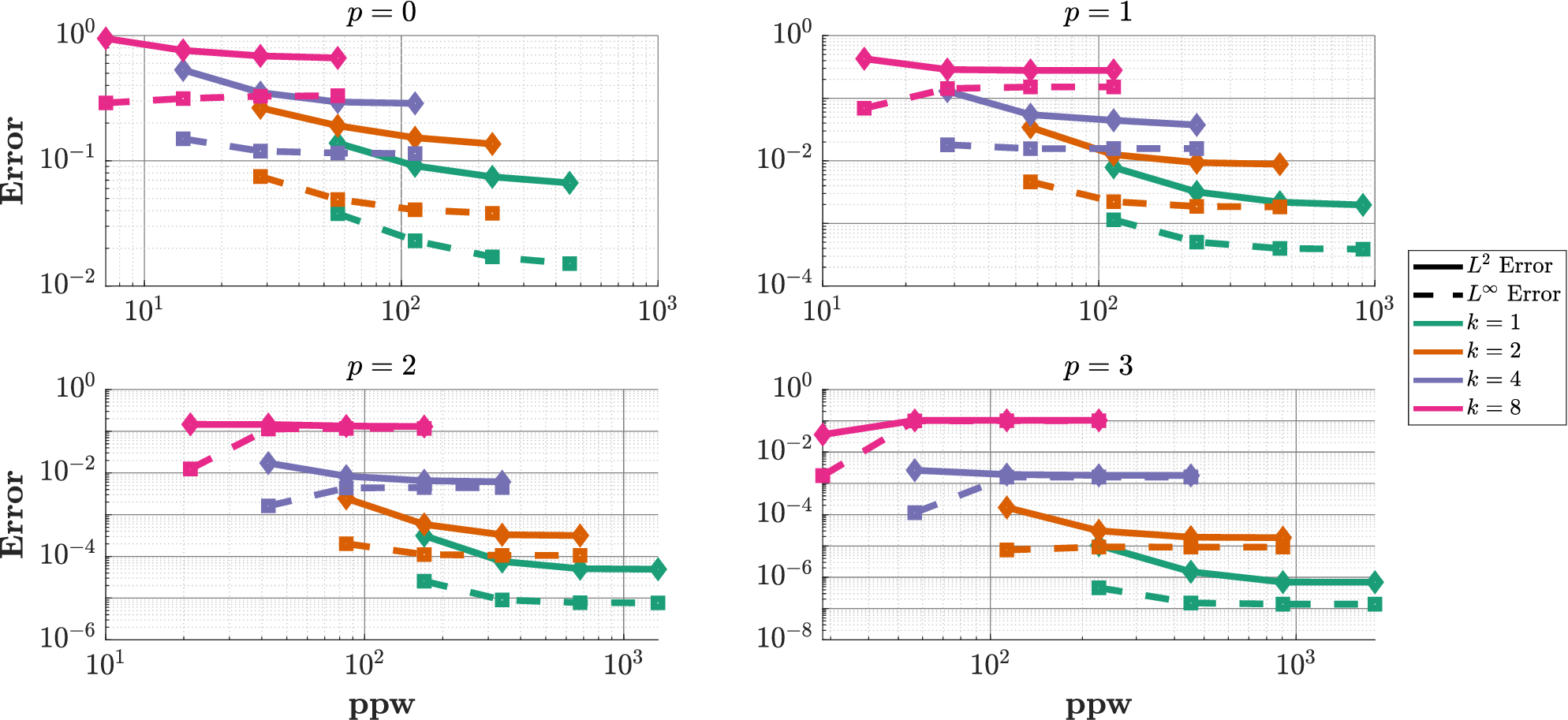}
\end{tabular}
    \caption{Log-Log Error plots for $u(x,y)=\sin(2\pi k x)\sin(2\pi k y)$ computed over the domain $[0,1]^2$. On the x-axis is plotted the points per wavelength in a single coordinate direction $\sqrt{N}(p+1)*2^{\# enh}/k$, with $\# enh$ representing the number of enhancements.}
    \label{fig:tensor sine}
\end{figure}

%% file: Figures/error_plots/tanh_and_gaussian.tex
\begin{figure}[hp!]
\begin{tabular}{c} 
\multicolumn{1}{c}{\textbf{Double Hyperbolic Tangent Data on Quadrilateral} $\mathbf{40\times 40}$ \textbf{Mesh}}\\ 
 \includegraphics[width=0.75\linewidth]{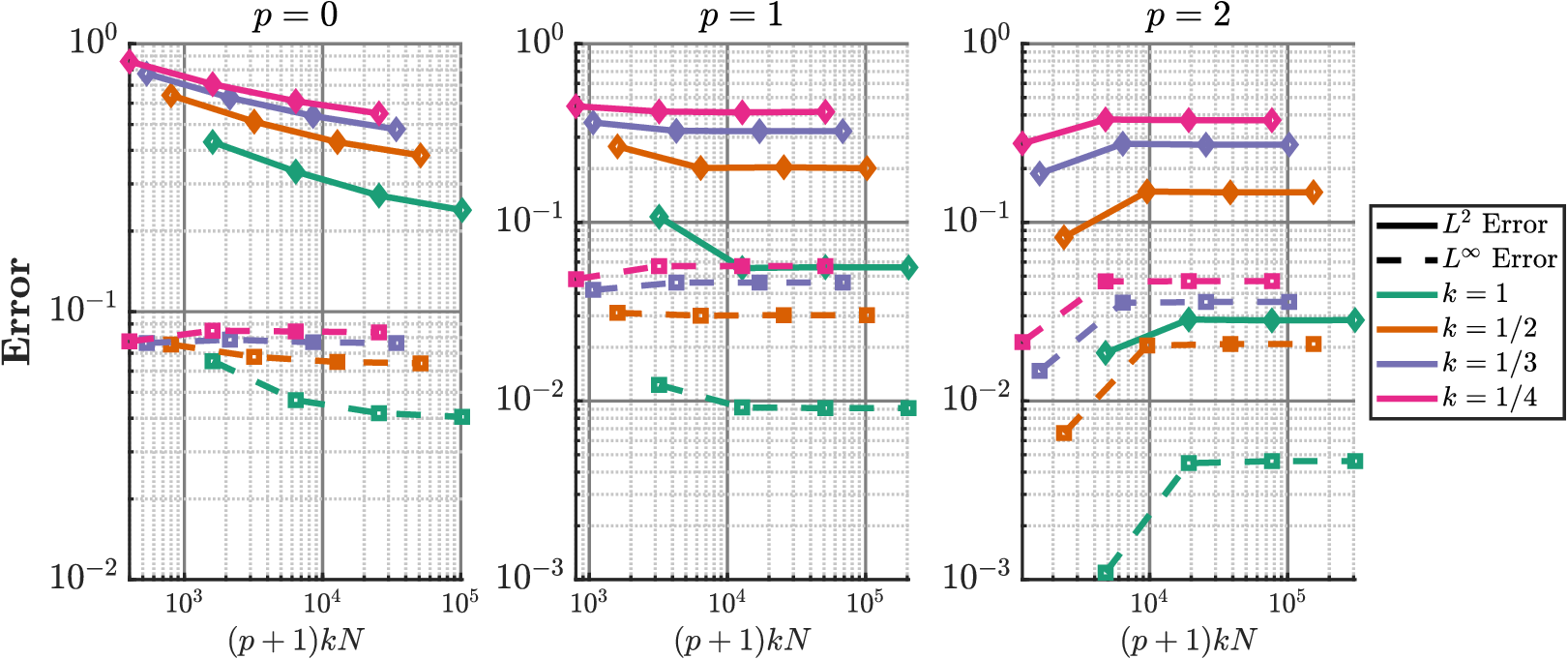}\\ 
  \multicolumn{1}{c}{\textbf{Gaussian Data on Quadrilateral} $\mathbf{40\times 40}$ \textbf{Mesh}}\\
 \includegraphics[width=0.75\linewidth]{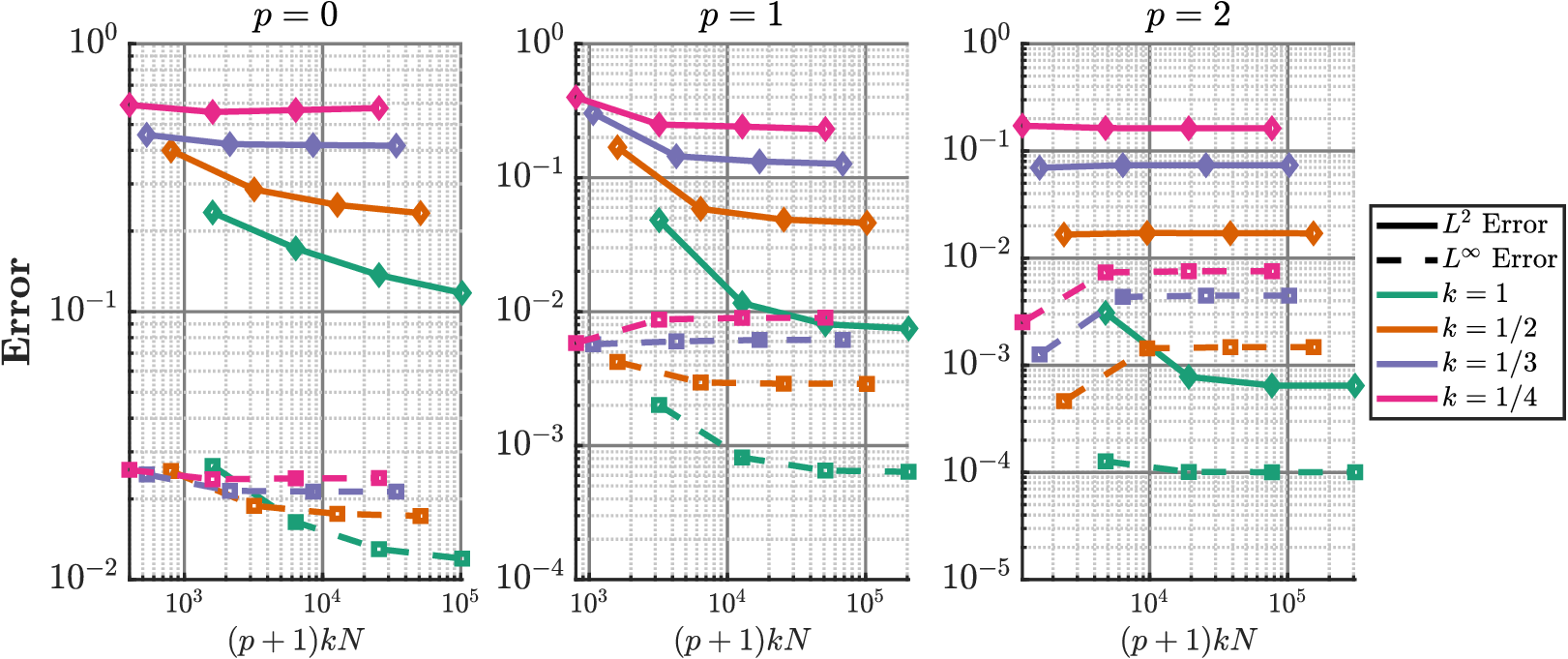} \\
\end{tabular}
    \centering
    \caption{Errors under refinement of the double hyperbolic tangent and Gaussian functions for varying scaling factors $k$. Here $N$ is the number of elements in a single coordinate direction $40*2^{\# enh}$, with $\# enh$ representing the number of enhancements.}
    \label{fig:tanh_and_gaussian}
\end{figure}

%% file: Figures/error_plots/Adaptive_vs_constant/adaptive_vs_constant.tex
\begin{figure}[tp!]
    \centering
    \begin{tabular}{c c} 
     %      \textbf{SVR10 Mesh} \\ 
 %\includegraphics[width=0.5\linewidth]{Figures/error_plots/Adaptive_vs_constant/SVR10_adaptive_vs_constant_tsine.jpg}\\ 
           \textbf{SVR100 Mesh} & \textbf{ISR2 Mesh} \\  
 \includegraphics[width=0.45\linewidth]{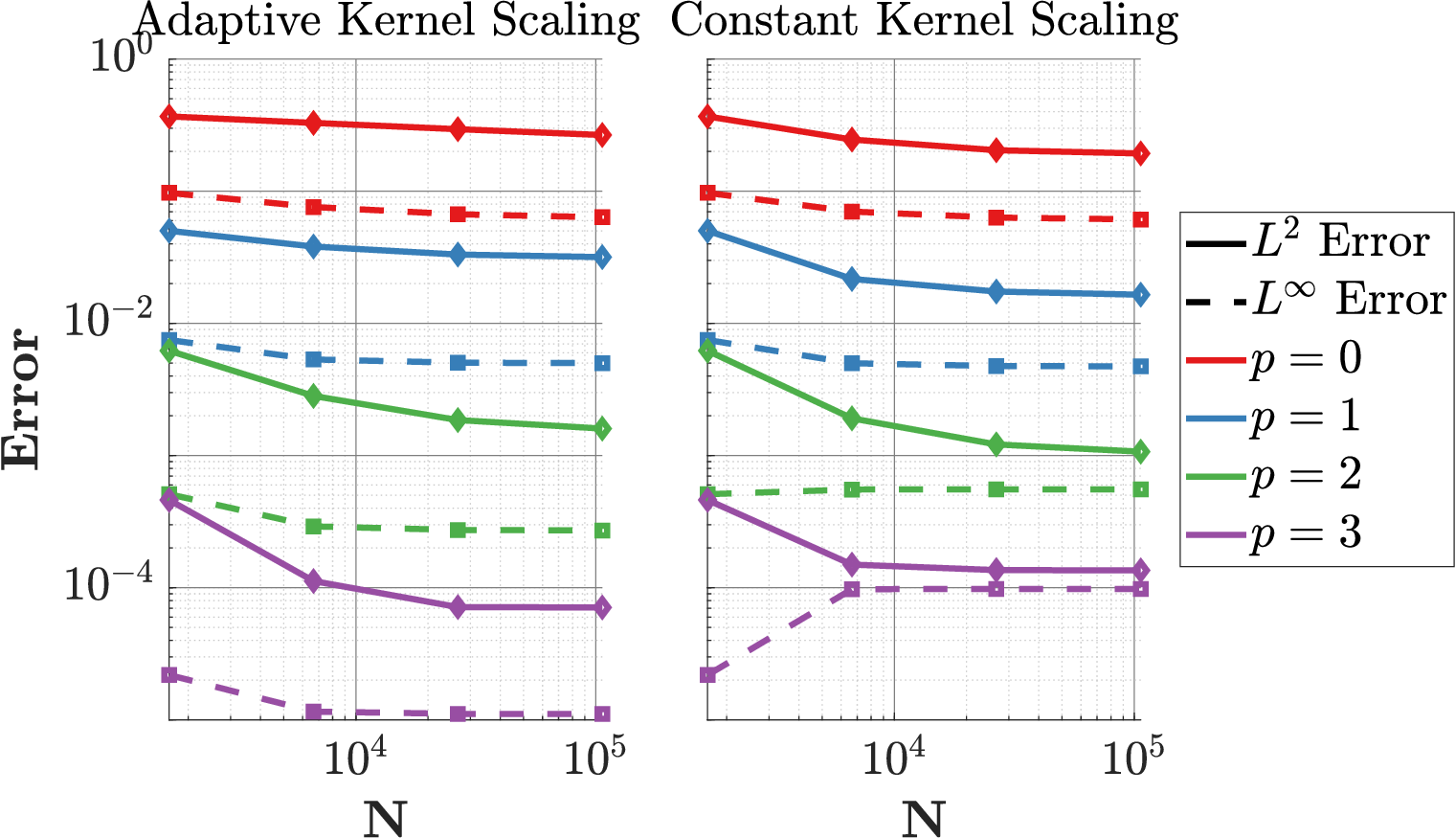}&
           
 \includegraphics[width=0.45\linewidth]{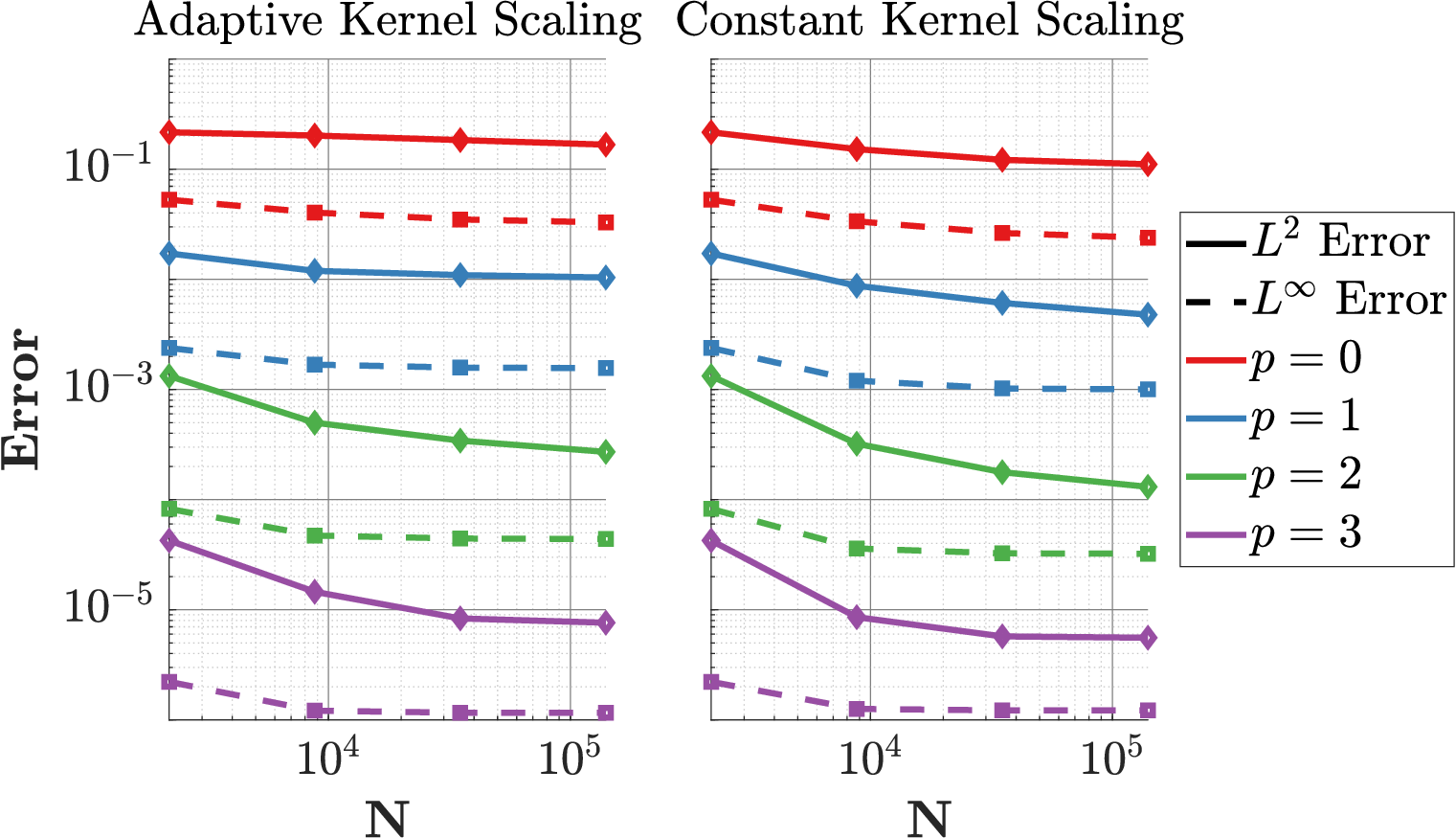}\\ 
    \end{tabular}
  \caption{Comparison of errors when using adaptive kernel scalings versus constant max edge length kernel scalings for $u(x,y)=\sin(2\pi x)\sin(2\pi y)$ on the SVR100 and ISR2 meshes over the domain $[0,1]^2$.}
    \label{fig:adapt_vs_constant_scaling}

\end{figure}

%% file: sections/adaptive_schemes.tex
\section{Application to Adaptive DG}\label{sec:adaptive_DG}

We now focus our attention on the application of LSIAC for $h$-adaptivity using multiwavelet detail indicators. Given a discrete problem, the aim is to compute an adaptive numerical approximation to $u_h$ that satisfies a local error target, $\eta_{TOL}$.  Simultaneously, we wish to reduce the number of degrees of freedom required when compared to a uniformly resolved reference approximation. To obtain such an approximation the $h$-adaptivity provided by the LSIAC-MRA methodology is employed. The work and analysis here mirrors that described in \cite{Bautista_paper} where, instead of SIAC-filtering, a reconstruction based on least-squares polynomials was employed. A benefit of the LSIAC approach is that the linear systems do not need to be solved over each sub-element -- instead a quadrature is performed.

\subsection{Adaptivity procedure}

The adaptivity procedure has the goal of computing the nonuniform mesh $\Omega^{target}_h$ for which the associated numerical solution satisfies user prescribed error targets. In general, the error of a computed numerical approximation is not known and to identify the elements that might benefit from subdivision, proxies termed error indicators are used. These are denoted here by $\eta$. To compare the quality with which a given error indicator tracks the actual error we introduce the effectivity index 
\begin{equation}
\iota_{eff}=\frac{\eta}{e_h}.
\end{equation}
Here the global error indicator $\eta$ is computed from the element-wise indicators $\eta_{\tau}$ by the relation $\eta=\sqrt{\sum_{\tau}\eta_{\tau}^2}$, and $e_h$ is defined to be the discrete approximation (computed as in the previous section) to the global $L^2$-norm of the error
\begin{equation}
e_h=\norm{u-u_h}_{L^2(\Omega)}.
\end{equation}
An effectivity index near 1 suggests that our chosen error indicator is serving as a good approximation of the actual error. For this work, we considered five different error indicators. Three are error indicators utilizing the SIAC enhancement and two are indicators appearing in the literature. They are
\begin{itemize}
\item The filtered difference indicator that simply considers the difference between the filtered solution and the approximation on a coarse grid:
\begin{equation}
\eta^{\star}_{\tau}=\norm{u_h^{\star}-u_h}_{L^2(\tau)}.
\end{equation}
\item The reconstruction difference indicator that considers the difference between the enhanced reconstruction on the finer mesh to the approximation on the coarser grid:
\begin{equation}
\eta^{REC}_{\tau}= \norm{\mathcal{P}^{n+1}u_h^{\star}-u_h}_{L^2(\tau)}. 
\end{equation}
\item The spectral decay (SD) and normalized small-scale energy decay (SSED) indicators for $p>0$ that consider the scaled effect of the highest polynomial mode, (see \cite{bautista_thesis} and the references therein):
\begin{equation}
\eta^{SD}_{\tau}=\frac{\norm{\sum^{N_p}_{\ell=1}u^{\ell}_{\tau}\phi^{\ell}_{\tau}-\sum^{N_p-1}_{\ell=1}u^{\ell}_{\tau}\phi^{\ell}_{\tau}}_{L^2({\tau})}}{\norm{\sum^{N_p}_{\ell=1}u^{\ell}_{\tau}\phi^{\ell}_{\tau}}_{L^2({\tau})}},
\end{equation}
\begin{equation}
\eta^{SSED}_{\tau}=\frac{\norm{\sum^{N_p}_{\ell=1}u^{\ell}_{\tau}\phi^{\ell}_{\tau}-\sum^{N_p-1}_{\ell=1}u^{\ell}_{\tau}\phi^{\ell}_{\tau}}_{L^2({\tau})}}{|{\tau}|^{1/2}}, 
\end{equation}
where $\ell$ is an index ordering the basis from lowest ($\ell=1$) to highest order modes ($\ell=N_p$) the highest order mode. Here $N_p=p+1$ in one dimension, and $(p+1)^2$ in two dimensions.
\item Lastly, the multiwavelet detail coefficient indicator 
\begin{align}
    \eta^{W}_{\tau}&= \norm{(\mathcal{P}^{n+1}-\mathcal{P}^{n})u_h^{\star}}_{L^2(\tau)}=\norm{\sum_{k=0}^p\tilde{d}_{k,\tau}^0\tilde{\psi}^0_{k,\tau}}_{L^2(\tau)}=\sqrt{\sum_{k=0}^p(\tilde{d}_{k,\tau}^0)^2},
    \intertext{which in two dimensions is given by}
\eta^W_{\tau}&=\sqrt{\sum_{k_x,k_y=0}^p(\tilde{d}^{0,\alpha}_{\mathbf{k},\tau})^2+(\tilde{d}^{0,\beta}_{\mathbf{k},\tau})^2+(\tilde{d}^{0,\gamma}_{\mathbf{k},\tau})^2},
\end{align}
where $\tilde{d}$ are the multiwavelet detail coefficients defined by starting a local-MRA about element $\tau$.
\end{itemize}
The adapted solution error is tracked over the course of the numerical experiment, and is denoted by
\[E_h=\norm{u-u^{adap}_h}_{L^2(\Omega)}.\]
 The refinement procedure is given below in Algorithm \ref{alg:Adapt loop}.

\begin{algorithm}
\caption{Mesh Adaption Procedure}
\begin{algorithmic}[1]

\Procedure{Adapt}{$\Omega_h^{(0)}$, $\eta_{TOL}$}       \Comment{Given initial mesh  and local error tolerance}
    \State Compute $u^{(i)}_h$ the approximate solution to the discrete problem on $\Omega_h^{(i)}$.
    \State Compute the error indicator $\eta_{\tau}$ over each element.
    \For{$\tau\in \Omega_h^{(i)}$}
    \If{$\eta_{\tau}\geq\eta_{TOL}$}
        \State Mark element $\tau$ for refinement.
        \EndIf
    \EndFor
     \If{$\eta_{\tau}<\eta_{TOL}$ for all $\tau$}
     \Return $(\Omega^{(i)}_h,u^{(i)}_h)$.
     \Else
       \State Refine marked elements and obtain the new mesh $\Omega^{(i+1)}_h$. 
       \State \Call{Adapt}{($\Omega_h^{(i+1)}$, $\eta_{TOL}$)}
        \EndIf
\EndProcedure

\end{algorithmic}
    \label{alg:Adapt loop}

\end{algorithm}

We now focus our attention on the performance of these different indicators.  To do so, we consider steady-state test problems. In the one dimensional case, we make use of explicit time stepping to steady-state.  When looping back from STEP 4 to STEP 1, we use the projection of $u^{(i)}_h$ onto $\Omega_h^{(i+1)}$ as the initialization to the discrete problem in STEP 1. This ensures that the whole procedure only requires the coarsest mesh's resolution of the initial data for initialization, rather than a re-initialization at ever increasing mesh resolutions. For the two-dimensional tests, we consider a problem that is a priori steady-state in the form of an elliptic problem. Hence, there is no modification of the refinement procedure for the two-dimensional case.

\subsection{Test Problems}
As a first PDE test case, consider the discrete problem of numerically approximating the solution to the forced viscous Burgers equation
\begin{equation}
u_t+\Big(\frac{u^2}{2}\Big)_x=\gamma u_{xx}+f(x,t)
\end{equation}
subject to appropriate boundary and initial conditions given below. We consider a modal RKDG discretization where the viscous term is treated by central fluxes and the convective term by a local Lax-Friedrich numerical flux (See \cite{Hes08W} Chapter 7. for the nodal version). In one dimension we consider the following test problems 
\begin{enumerate}
    \item Viscous Burgers subject to Dirichlet boundary conditions $u(\pm1,t)=\mp 1$ with $\gamma=0.02$, and initial data given by $u(x,0) = -\tanh(x/2\gamma)$, with $f=0$.
    \item Viscous Burgers subject to periodic boundary conditions with $\gamma=0.5$, and initial data given by $u(x,0) =\sin(2\pi x)$, with a forcing term of $f(x,t)=2\pi\sin(2\pi x)\cos(2\pi x)+4\pi^2\gamma\sin(2\pi x)$.
\end{enumerate}
The first test problem represents a function with a sharp gradient in the interior of the domain. The boundary conditions are satisfied to within an exponentially small term (see \cite{burgers_sol}), and so the initial condition is used as a proxy for the exact solution. As the data is non-periodic, the enhancement procedure is not applied in elements where the kernel support would overlap the domain boundary. The second test problem represents approximation of a smooth function, where the forcing term enforces convergence to the initial condition. In both cases, we enforce a 2:1 rule \cite{Giraldo} requiring elements be no larger than twice the size of their neighbors. The discrete problem given by the RKDG scheme is evolved to a steady state using a 4th order 2N-storage scheme given in \cite{Carpenter}. The time residual is computed as the $\ell_1$ difference of the modes at successive times steps and is required to be less than $10^{-12}$. 

In the two-dimensional tests, we consider an elliptic differential equation,
\begin{equation}
    \nabla \cdot (\nabla u)=f(\mathbf{x}),
\end{equation}
subject to Dirichlet boundary conditions. We consider test problems on the domain $\Omega=[0,1]^2$, where the forcing functions and boundary conditions are determined by the exact solutions given by:
\begin{enumerate}
    \item $u(\mathbf{x})=\sin(2\pi x)\sin(2\pi y)$,
    \item$u(\mathbf{x})=\text{exp}(-100\norm{\mathbf{x}-(1/2,1/2)}_{\ell^2})$.
\end{enumerate}
The first test problem represents a smooth test problem with no sharp features. For this type of function we expect the errors to not significantly differ in different parts of the domain.  As such, we expect the refinement procedure to either refine everywhere, or not at all. The second problem, while possessing a smooth solution, it  is very localized in space and as a result has a sharp gradient. Here we expect the refinement procedure to refine about the Gaussian spike in the interior and provide little or no refinement near the boundaries.

With respect to the DG method employed, the spatial discretization is performed by the Symmetric Interior Penalty (SIP) 
 DG method (see \cite{Arnold}) with penalty parameter chosen as in \cite{Alhawwary}. The solver used enforces a 2:1 rule which prohibits elements more than 1 level different from being neighbors (see \cite{Giraldo}). Note that despite the Dirichlet boundary conditions, we are applying the filter in a periodic fashion which for these test problems poses no issue. In future work, position-dependent SIAC filters may be employed for finite domains.
 
\subsection{Numerical Results: Viscous Burgers}

For the purposes of adaptivity indication in the one-dimensional case, we have observed that selecting the kernel scaling to be equal to the minimum over all element widths performs much better than the standard wisdom for obtaining superconvergence in DG data which usually elects to use the maximum element width. We employ this scaling choice to generate the results given below. For the first initial condition, an initial resolution of $2^3=8$ elements ($2^4=16$ for the $p=3$ case) is considered. A series of seven refinements (six for the $p=3$ case) is then performed. For this problem, $\eta_{TOL}$ has been selected to be the average element-wise $L^2$-error of the uniform $2^8=256$ element approximation. These values for $\eta_{TOL}$ for varying $p$ are given in Table \ref{tab:eta_tol}. 

\begin{figure}
    \centering
    \begin{tabular}{ccc}
      \multicolumn{1}{c}{Errors} &  \multicolumn{1}{c}{Effectivity Index} &  \multicolumn{1}{c}{Mesh Refinements}\\
      \multicolumn{3}{c}{$p=0$}\\
       \includegraphics[width=0.33\linewidth]{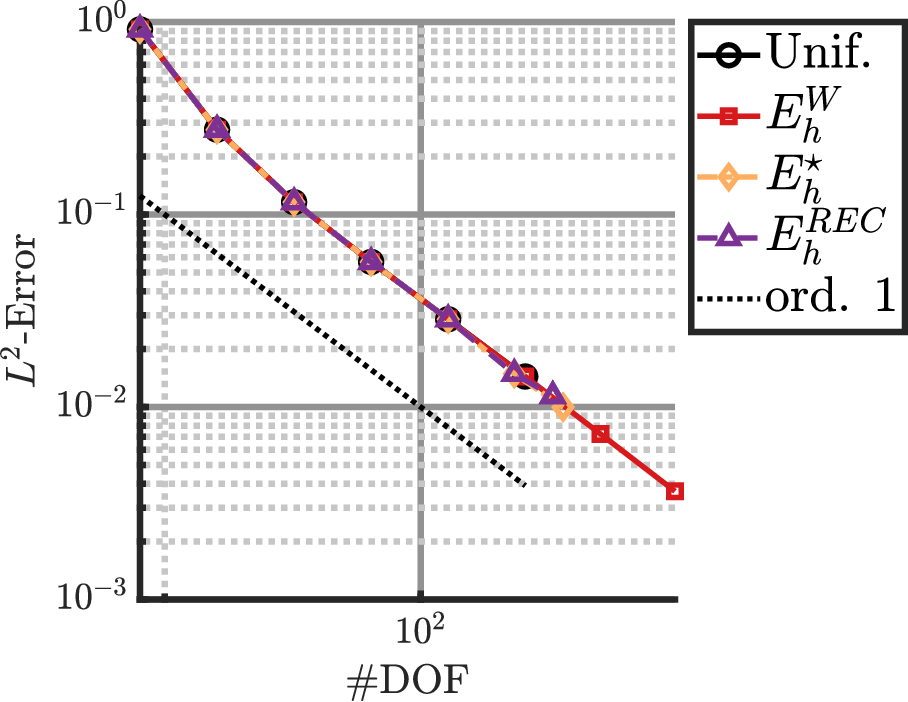} &  \includegraphics[width=0.33\linewidth]{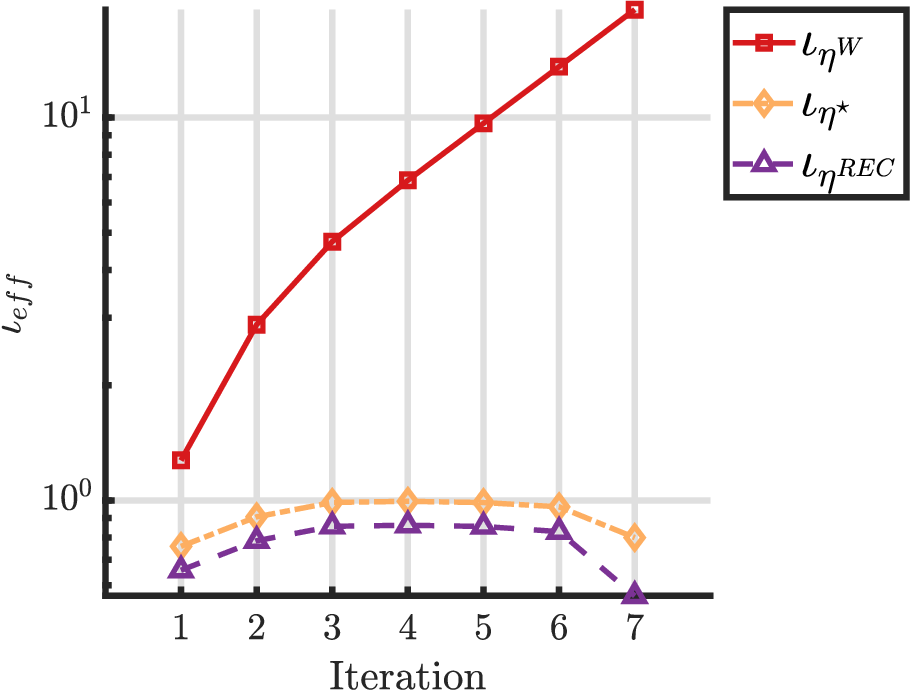} &

       \includegraphics[width=0.33\linewidth]{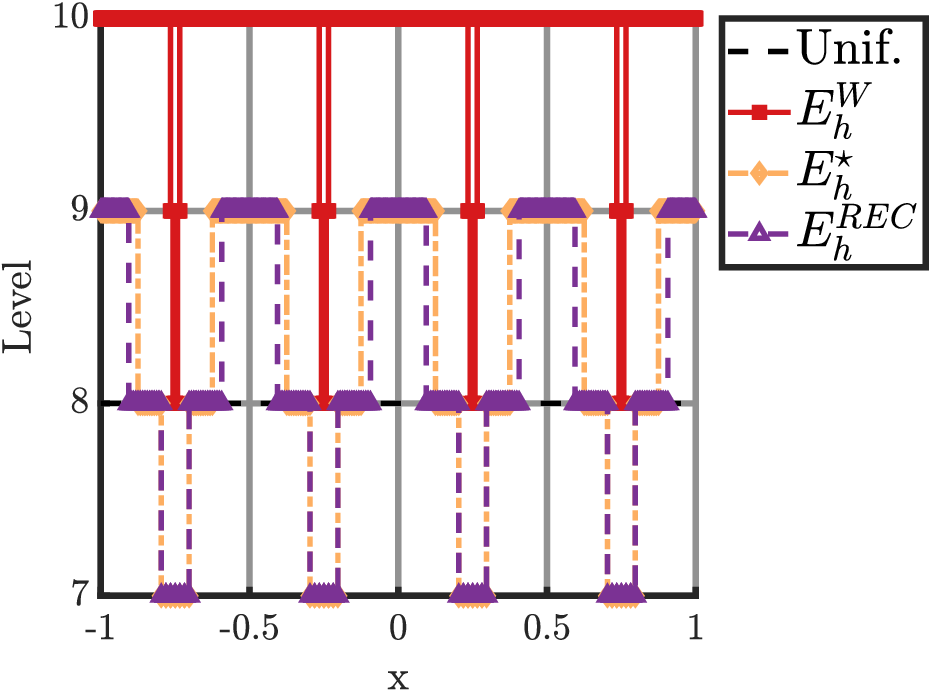}\\
             \multicolumn{3}{c}{$p=1$}\\
       \includegraphics[width=0.33\linewidth]{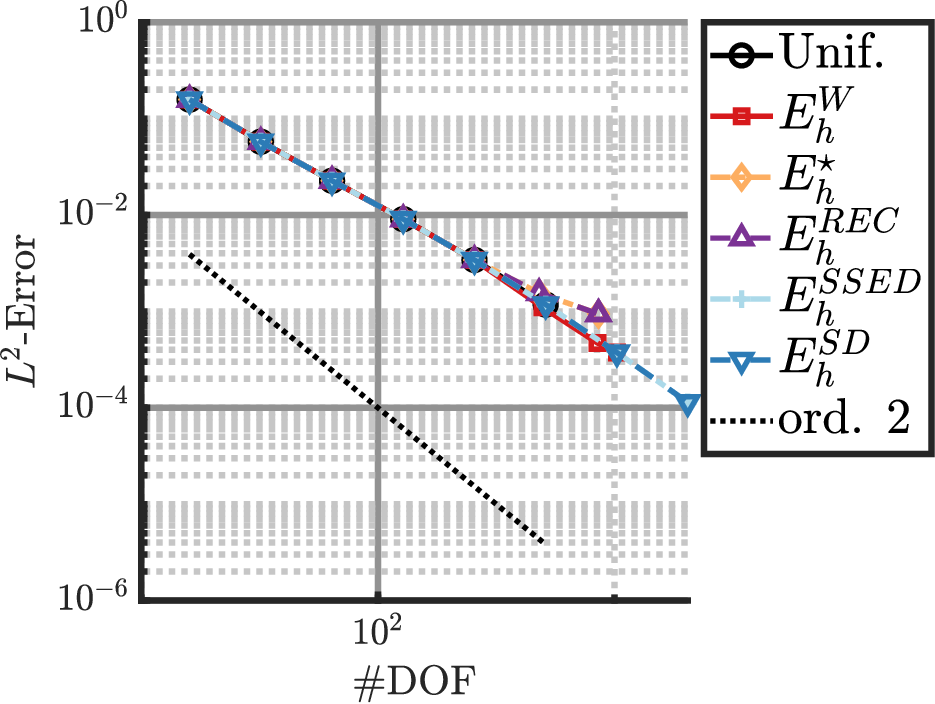}&  \includegraphics[width=0.33\linewidth]{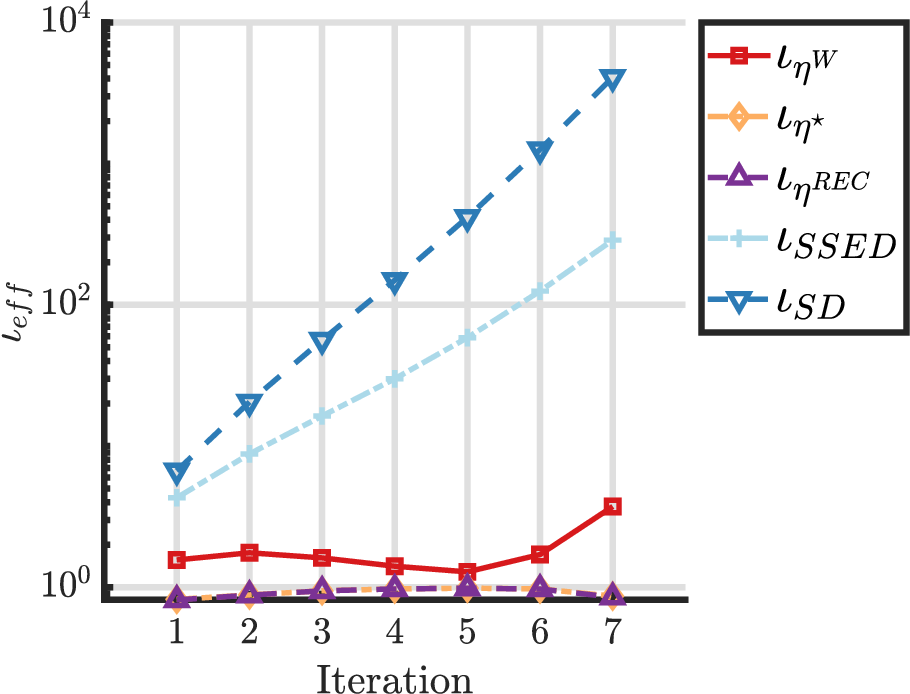} & \includegraphics[width=0.33\linewidth]{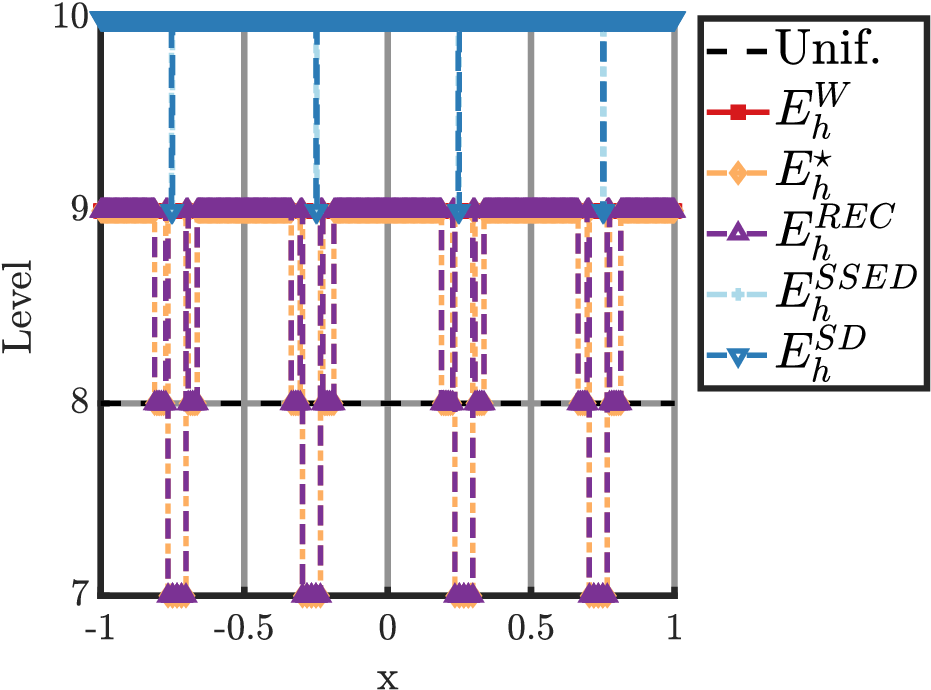}\\
        \multicolumn{3}{c}{$p=2$}\\
        \includegraphics[width=0.33\linewidth]{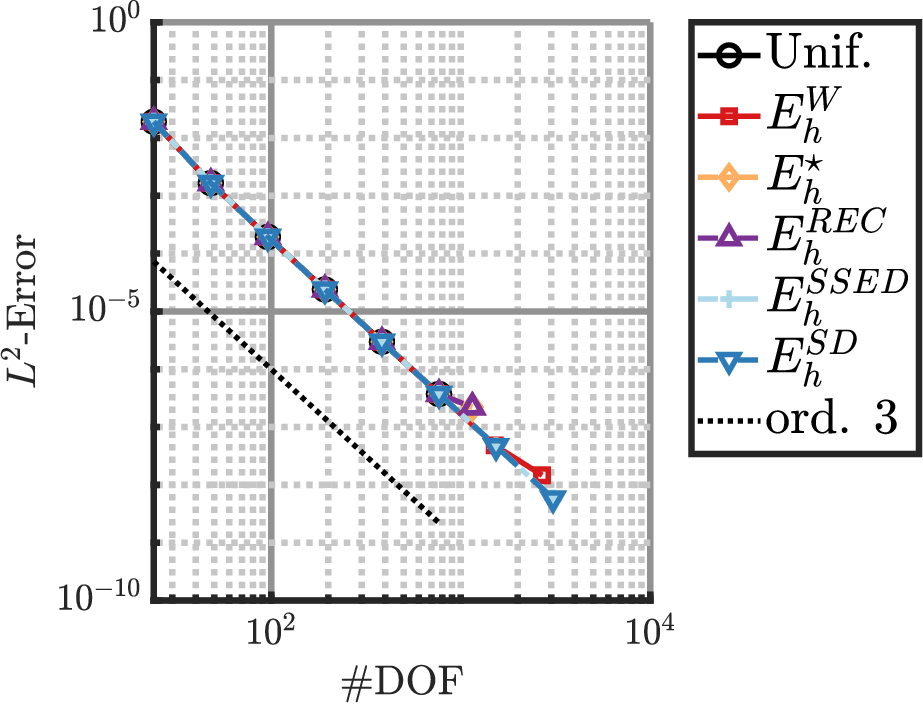} &  \includegraphics[width=0.33\linewidth]{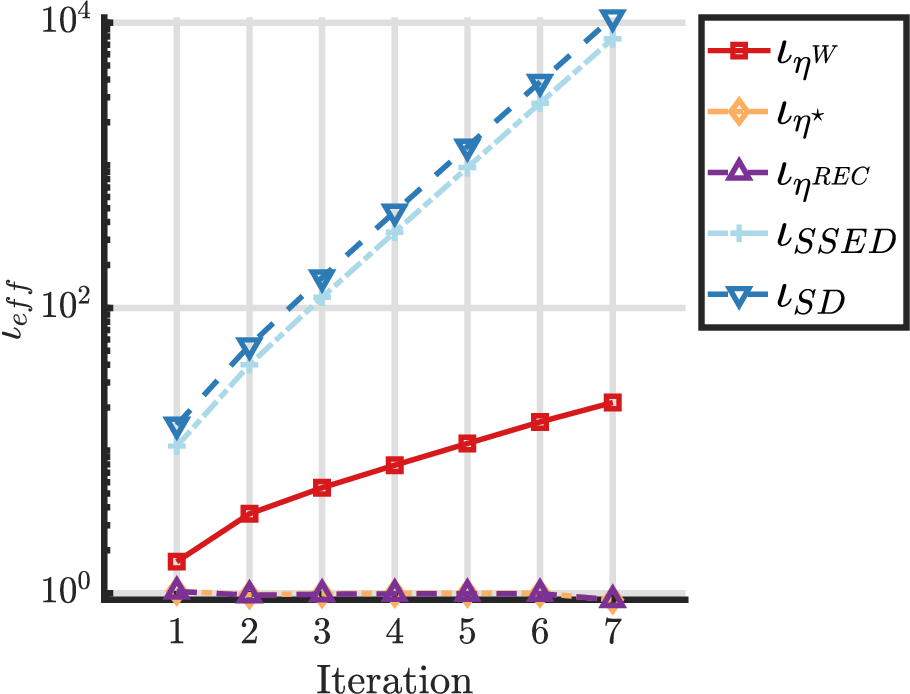} &  \includegraphics[width=0.33\linewidth]{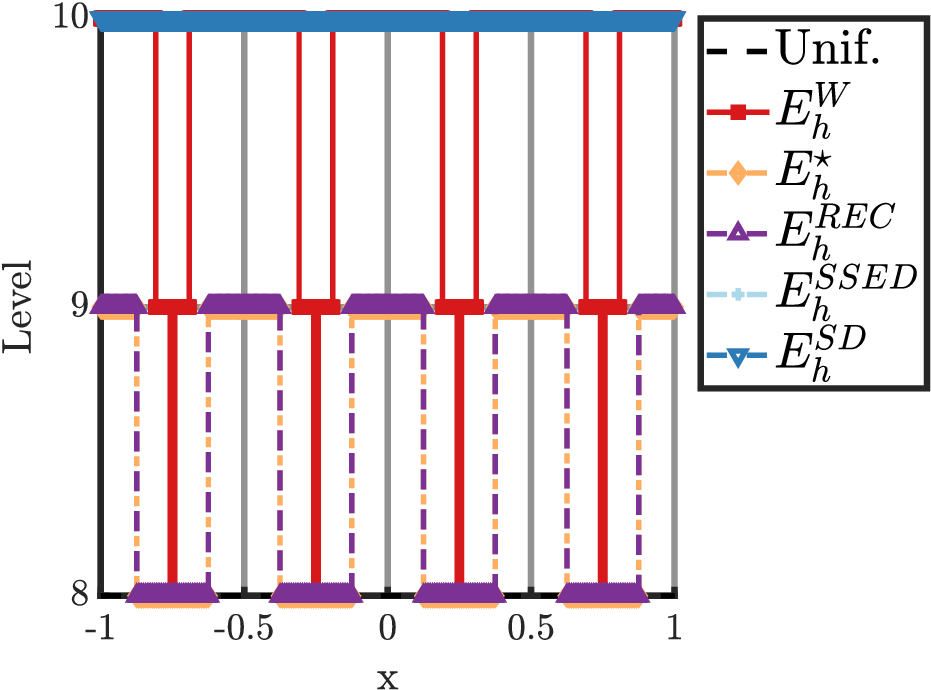}\\
       \multicolumn{3}{c}{$p=3$}\\
 \includegraphics[width=0.33\linewidth]{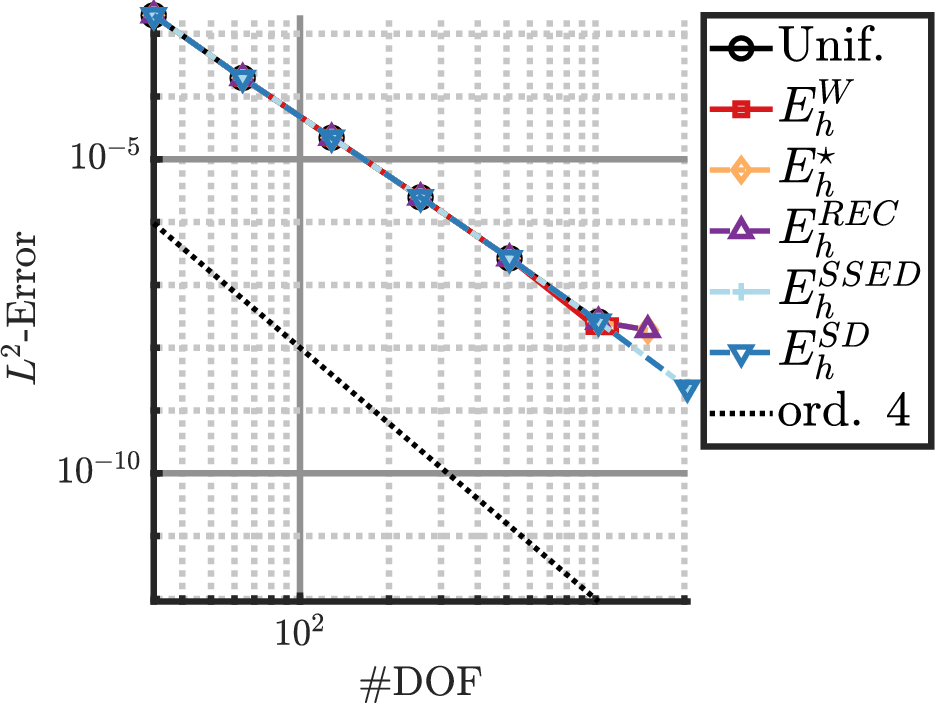} 
         & \includegraphics[width=0.33\linewidth]{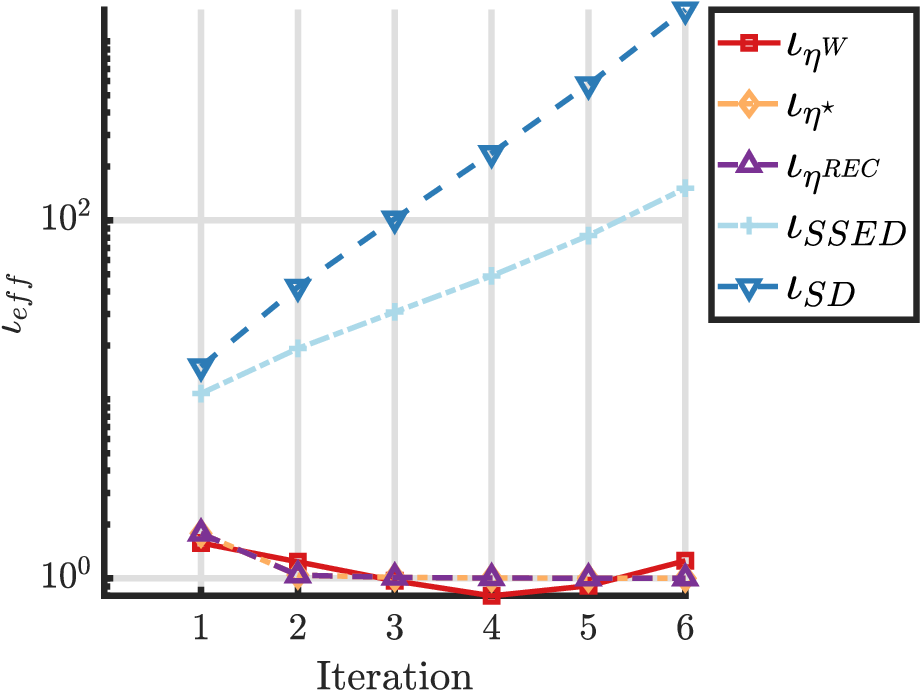}
       & \includegraphics[width=0.33\linewidth]{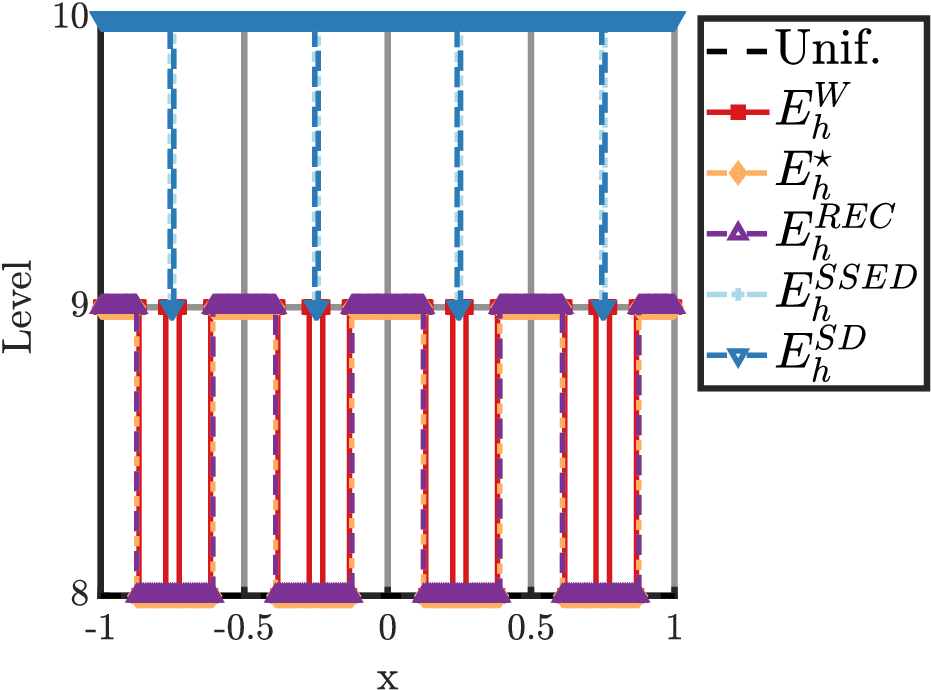} \\
    \end{tabular}
    \caption{$L^2$-Errors, effectivity indices, and final meshes achieved under different refinement indicators for the one-dimensional Sine IC. Here a level equal to $n$ indicates an element width of $2^{-n+1}$.}
    \label{fig:sine_varypoly}
    \end{figure}

 Consider first the smooth sine initial data. For this problem there is not a single feature or region of the mesh that has dramatically greater errors than any other region. As a result, the adaptivity procedure has to refine everywhere to obtain an approximation of similar quality to the uniformly refined standard. This is displayed in Column 1 of Figure \ref{fig:sine_varypoly} where it is observed that the adapted approximation closely tracks the uniformly refined approximations. This does not lead to any reduction in the number of degrees of freedom, and in fact a slight over-refinement occurs for the SIAC-based indicators. However, the SD and SSED indicators over-refine much more. The SIAC-based effectivity indices track well with the actual error, though the multiwavelet indicator consistently underestimates the actual error (see Column 2 of Figure \ref{fig:sine_varypoly}). Column 3 of \ref{fig:sine_varypoly} depicts the final meshes achieved by the varying error indicators. All the indicators refine the mesh everywhere in the domain, though the SD and SSED indicators over-refine the mesh relative to the uniform reference.

\begin{figure}
    \centering
    \begin{tabular}{c c c }
      \multicolumn{1}{c}{Errors} &  \multicolumn{1}{c}{Effectivity Index} &  \multicolumn{1}{c}{Mesh Refinements}\\
      \multicolumn{3}{c}{$p=0$}\\
       \includegraphics[width=0.33\linewidth]{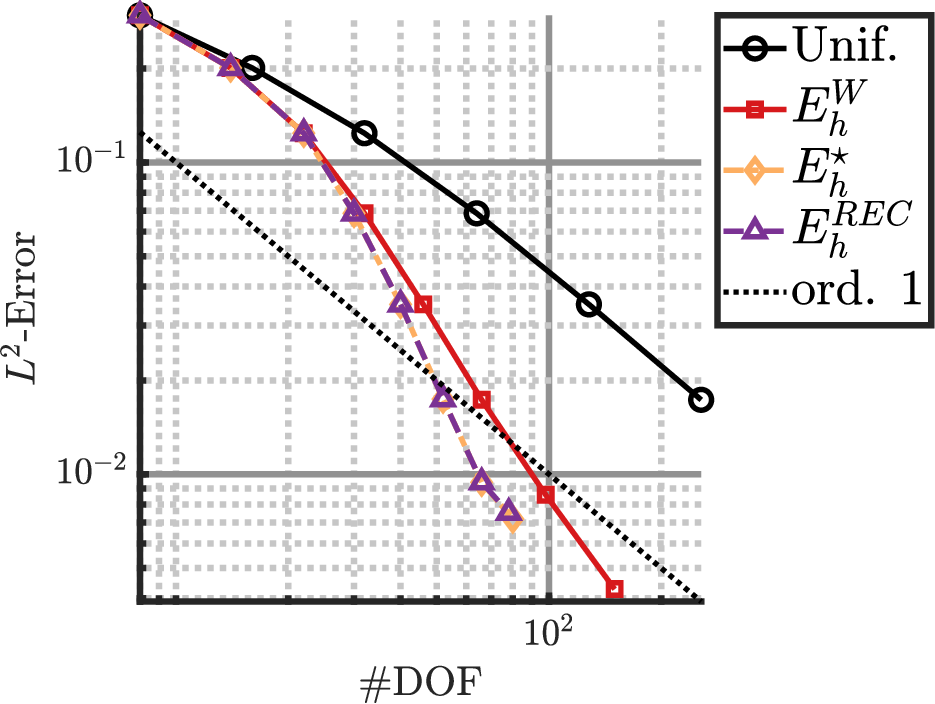} &  \includegraphics[width=0.33\linewidth]{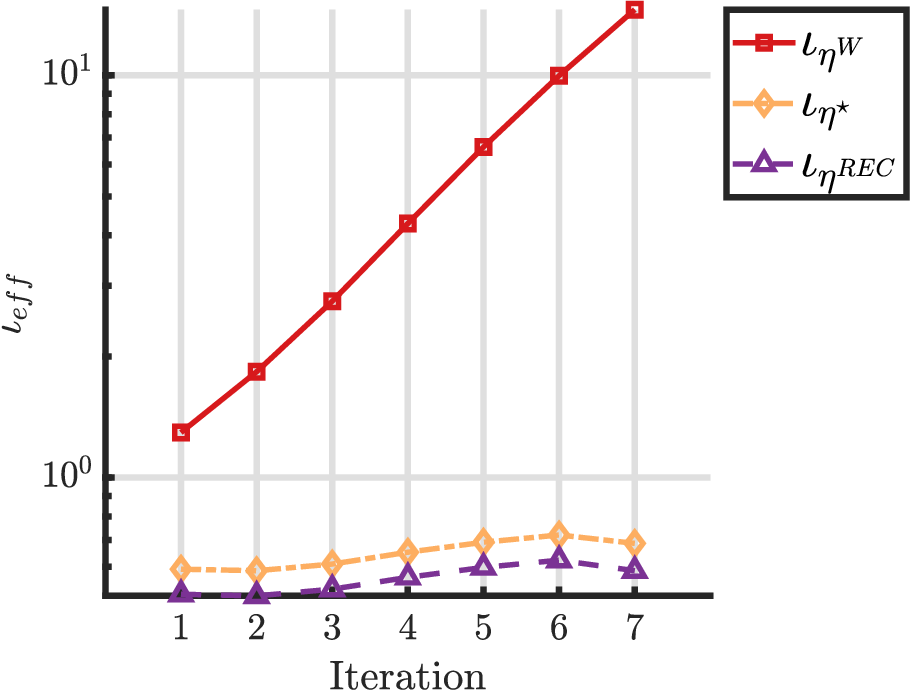} &

       \includegraphics[width=0.33\linewidth]{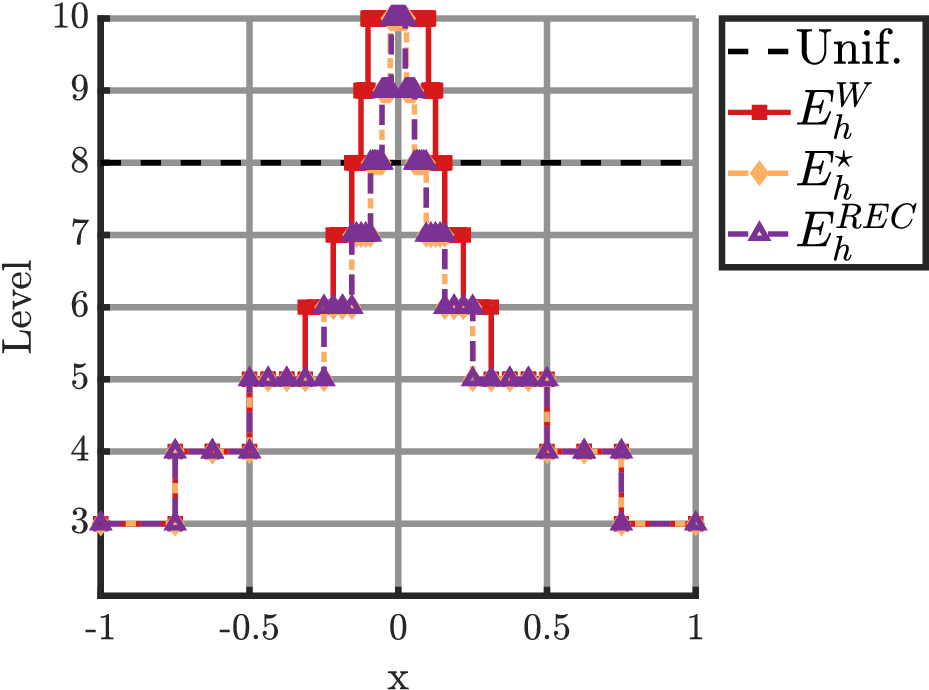}\\
             \multicolumn{3}{c}{$p=1$}\\
       \includegraphics[width=0.33\linewidth]{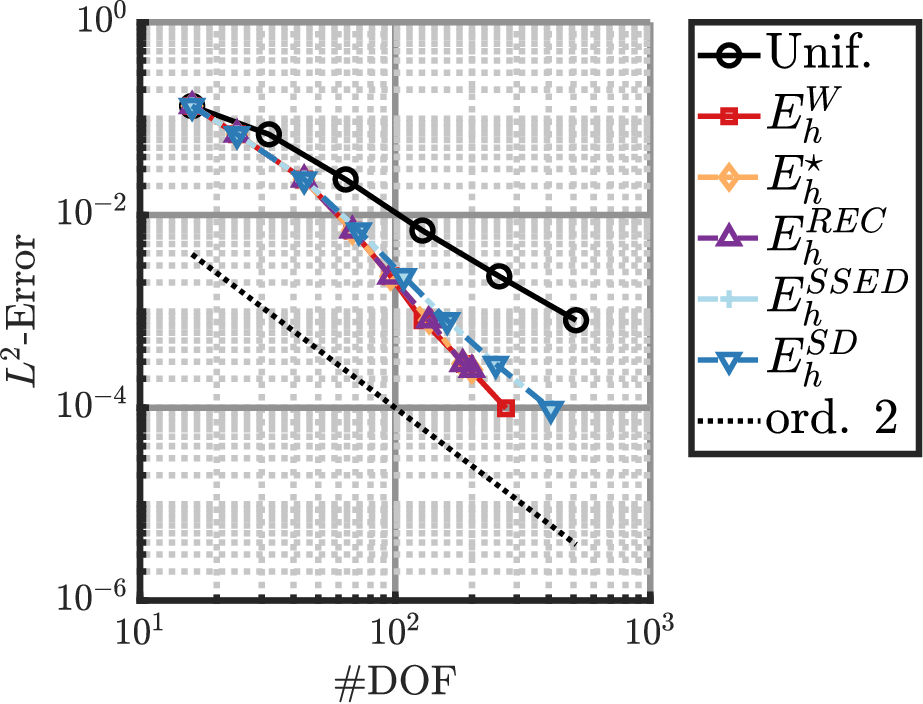}&  \includegraphics[width=0.33\linewidth]{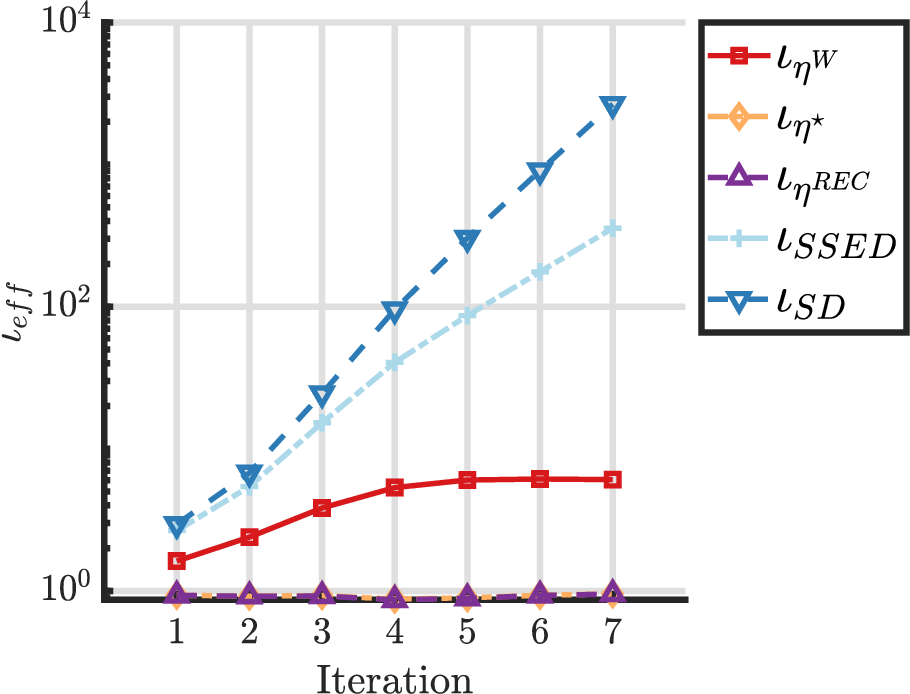} & \includegraphics[width=0.33\linewidth]{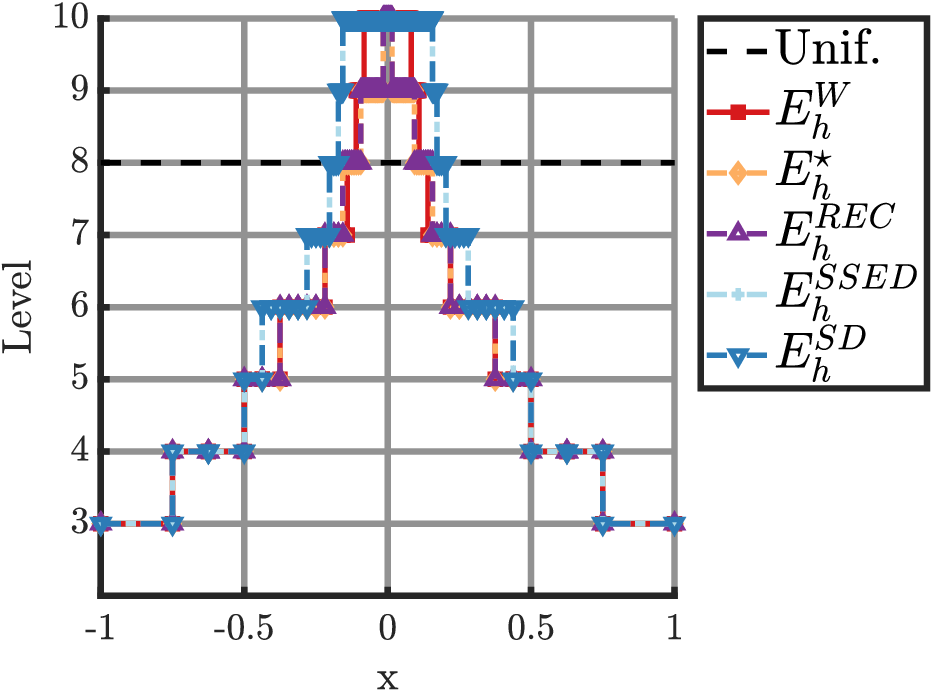}\\
        \multicolumn{3}{c}{$p=2$}\\
        \includegraphics[width=0.33\linewidth]{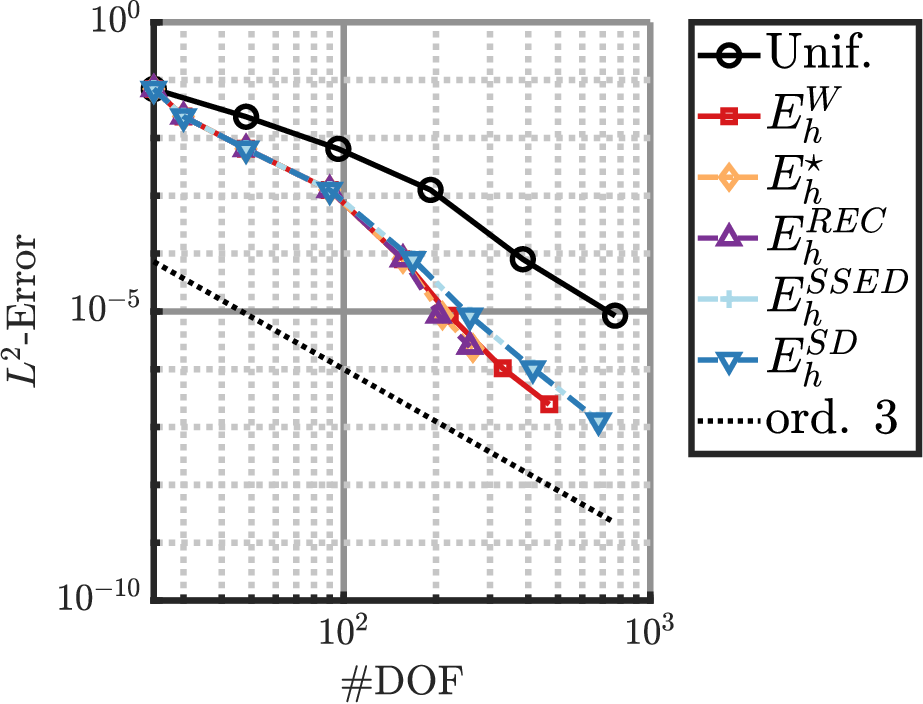} &  \includegraphics[width=0.33\linewidth]{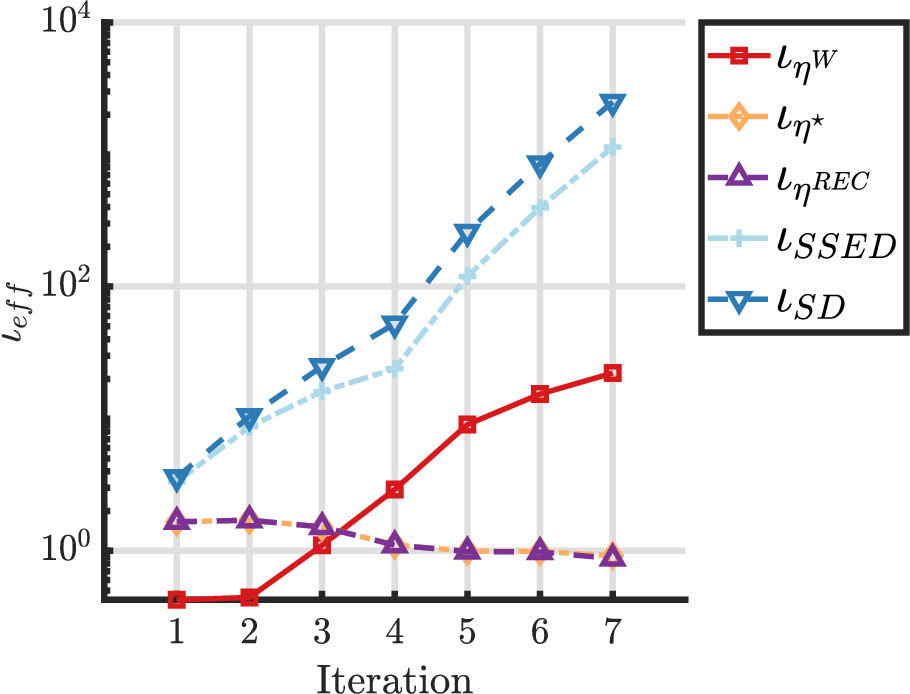} &  \includegraphics[width=0.33\linewidth]{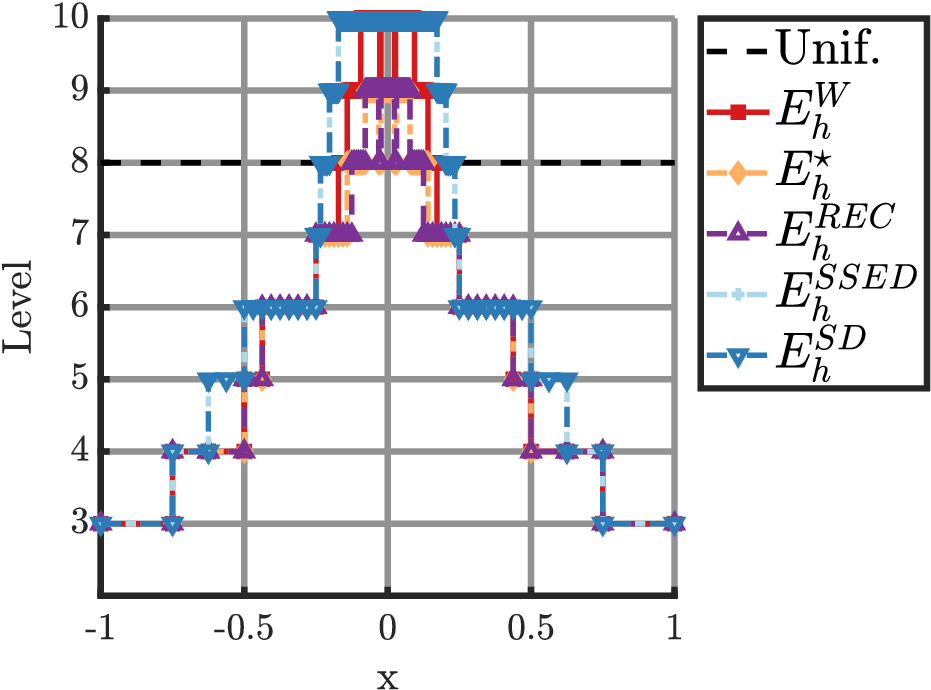}\\
       \multicolumn{3}{c}{$p=3$}\\
 \includegraphics[width=0.33\linewidth]{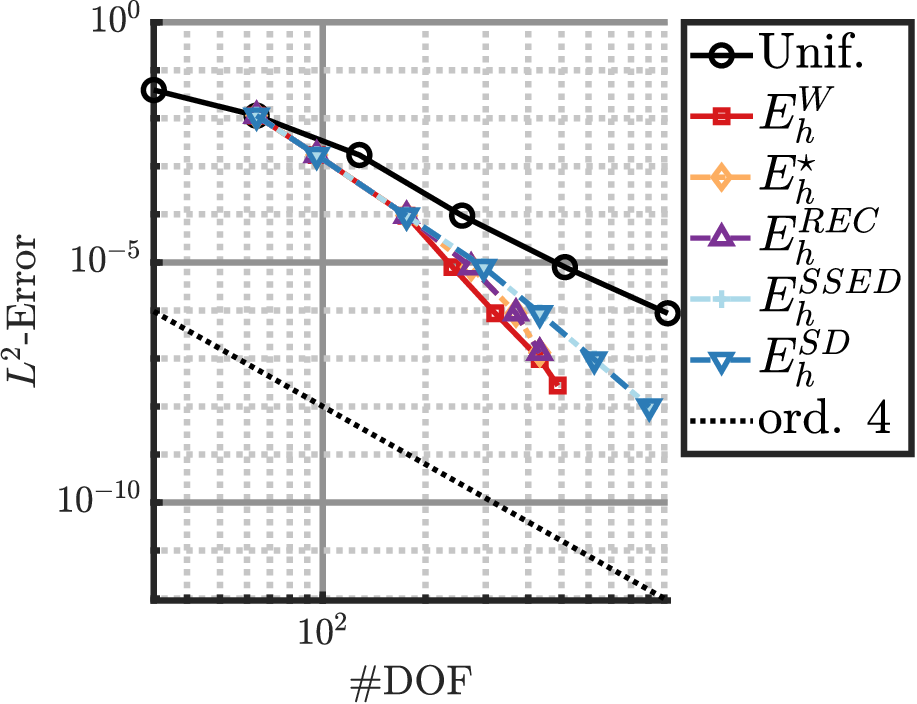} 
         & \includegraphics[width=0.33\linewidth]{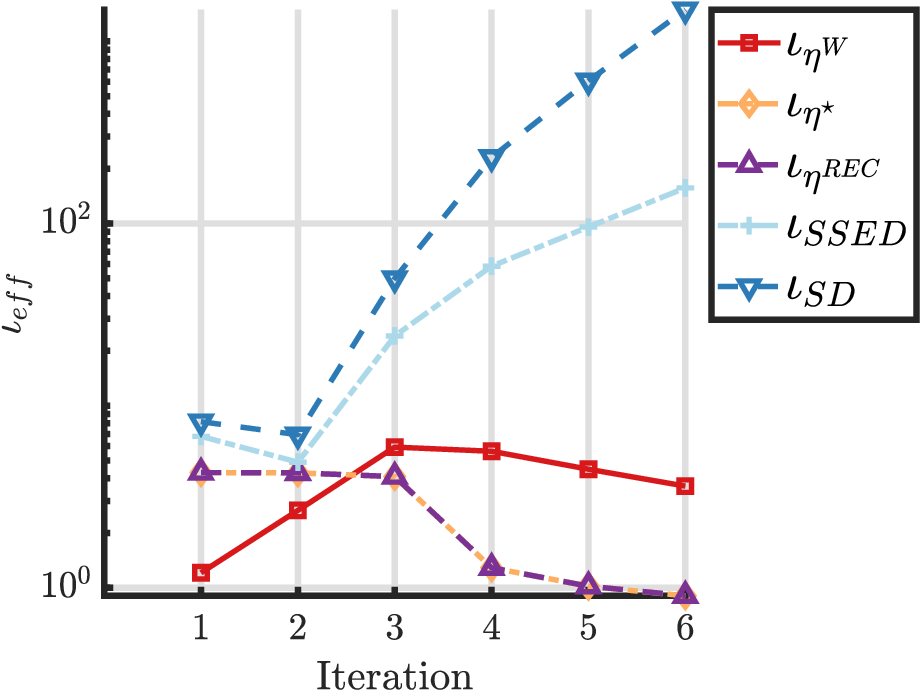}
       & \includegraphics[width=0.33\linewidth]{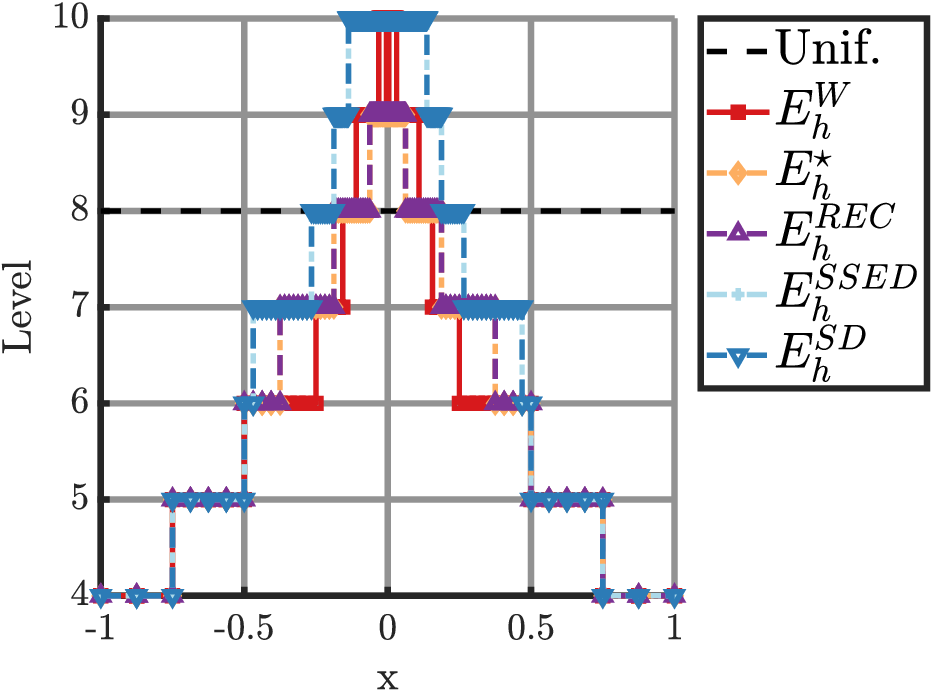} \\
    \end{tabular}
    \caption{$L^2$-Errors, effectivity indices, and final meshes achieved under different refinement indicators for the one-dimensional hyperbolic tangent IC. Here a level equal to $n$ indicates an element width of $2^{-n+1}$.}
    \label{fig:tanh_varypoly}
    \end{figure}

For the hyperbolic tangent initial data, accuracy comparable to the finest mesh of the uniformly refined solution is obtained with a significant reduction in numbers of degrees of freedom for the SIAC-based indicators (see Column 1 in Figure \ref{fig:tanh_varypoly}). Generally the filtered and reconstructed indicators perform similarly, as do the SD and SSED indicators. Effectivity indices are provided in  Column 2 of Figure \ref{fig:tanh_varypoly}. The error indicators track the actual error well for the filter-based indicators, though the multiwavelet indicator switches between over- and underestimating the error. The SSED and SD indicators overestimate the error more and as a result causes over-refinement and a large effectivity index. Looking at the final mesh resolutions obtained by these procedures (see column 3 \ref{fig:tanh_varypoly}), in all cases the greatest resolution is obtained about the sharp gradient near $x=0$. For $p>0$, the multiwavelet indicator results in refinement in between the filtered error indicators and the SD/SSED indicators.

\begin{table}[hp!]
    \centering
    \begin{tabular}{|c|c|c|c|c|c|}\hline
   $N_{elm}$& IC$/\eta_{TOL}$& $p=0$ & $p=1$ &  $p=2$ & $p=3$\\ \hline
   $256$  & 1D $\tanh$ IC & 1.1000e-03 &5.0603e-05&5.2650e-07&5.4431e-08  \\ \hline
     $256$  &  1D sine IC  &9.0045e-04 &7.3578e-05 & 2.3291e-08&1.6228e-09   \\ \hline
        $16,384$  & 2D Gaussian IC & 7.9961e-03 &3.0852e-04&2.4773e-06& 3.5193e-10\\ \hline
     $16,384$  &  2D sine IC  &2.9829e-02 &2.3229e-04 & 1.4290e-06&2.7362e-09   \\ \hline
    \end{tabular}
    \caption{Average element-wise $L^2$-errors of the $N_{elm}$ uniform element approximations. These are used as the user supplied error tolerances, $\eta_{TOL}$.}
    \label{tab:eta_tol}
\end{table}

\subsection{Numerical Results: Poisson Equation}

For the two-dimensional case, we again choose the kernel scaling equal to the minimum edge length over the mesh. An initial resolution of $16\times16$ elements is considered, and a series of four refinements is performed. Again, $\eta_{TOL}$ has been selected as given in Table \ref{tab:eta_tol} and corresponds to the average element-wise $L^2$-error of a uniform $128\times 128$ mesh approximation.

\begin{figure}[hp!]
    \centering
    \begin{tabular}{c c  }
      \multicolumn{1}{c}{Errors} &  \multicolumn{1}{c}{Effectivity Index}  \\
      \multicolumn{2}{c}{$p=0$}\\
       \includegraphics[width=0.33\linewidth]{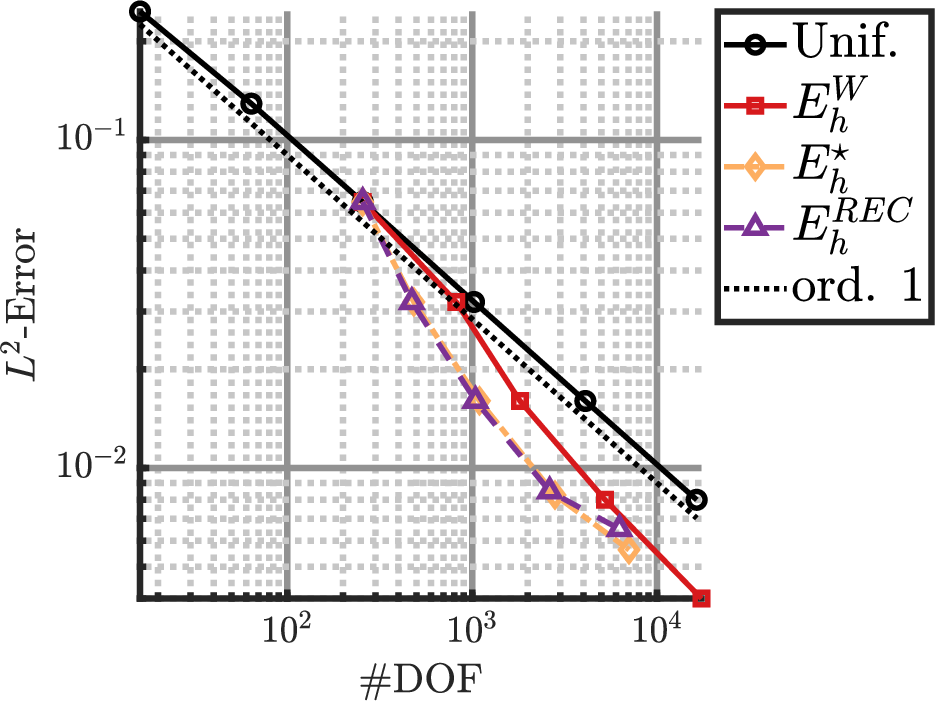}&  \includegraphics[width=0.33\linewidth]{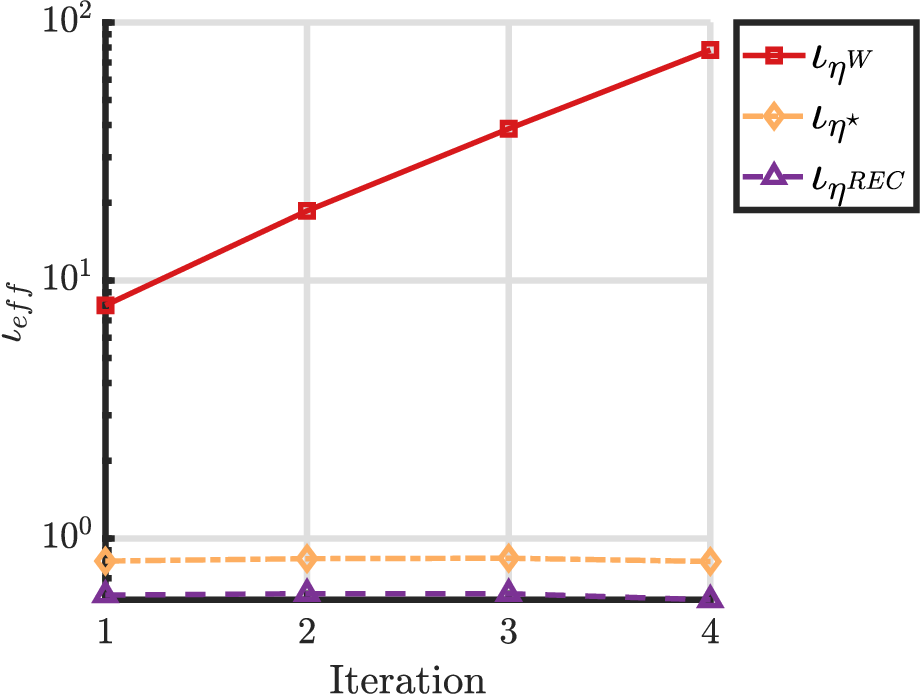} \\
             \multicolumn{2}{c}{$p=1$}\\
        \includegraphics[width=0.33\linewidth]{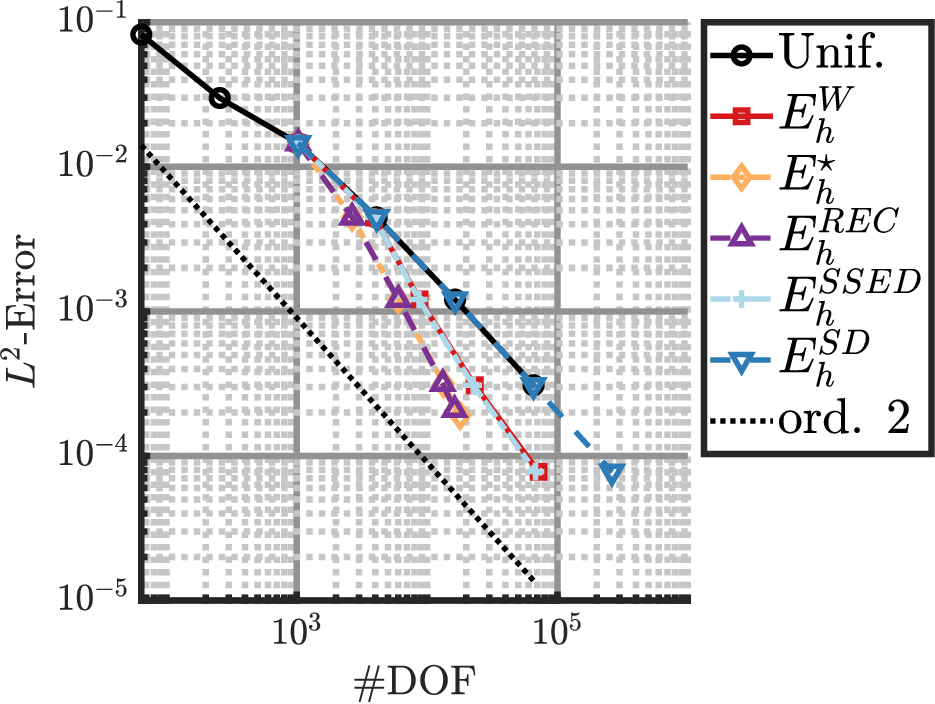} &  \includegraphics[width=0.33\linewidth]{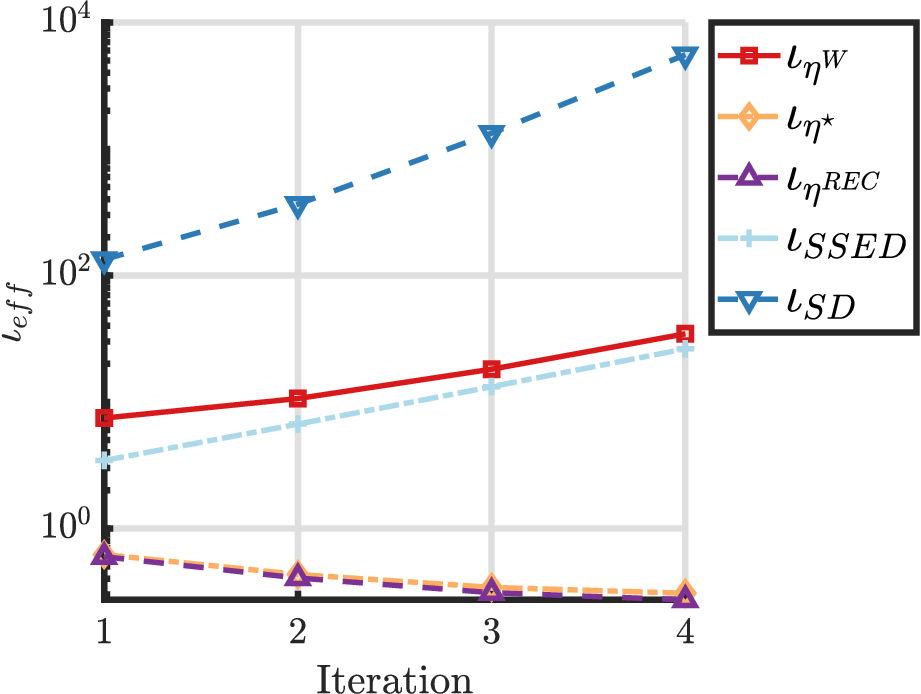} \\
        \multicolumn{2}{c}{$p=2$}\\
           \includegraphics[width=0.33\linewidth]{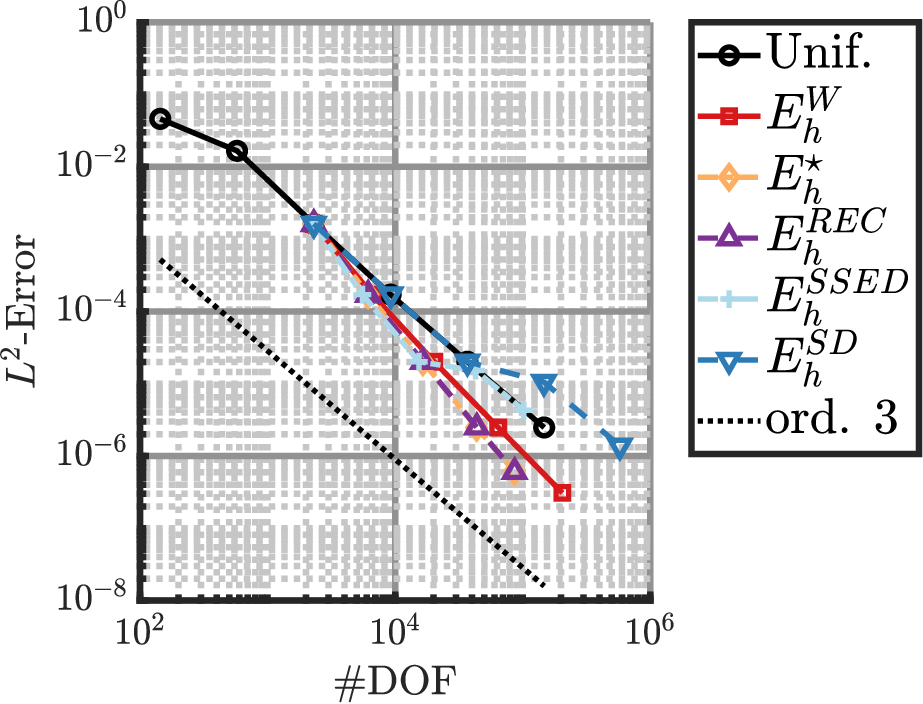} &  \includegraphics[width=0.33\linewidth]{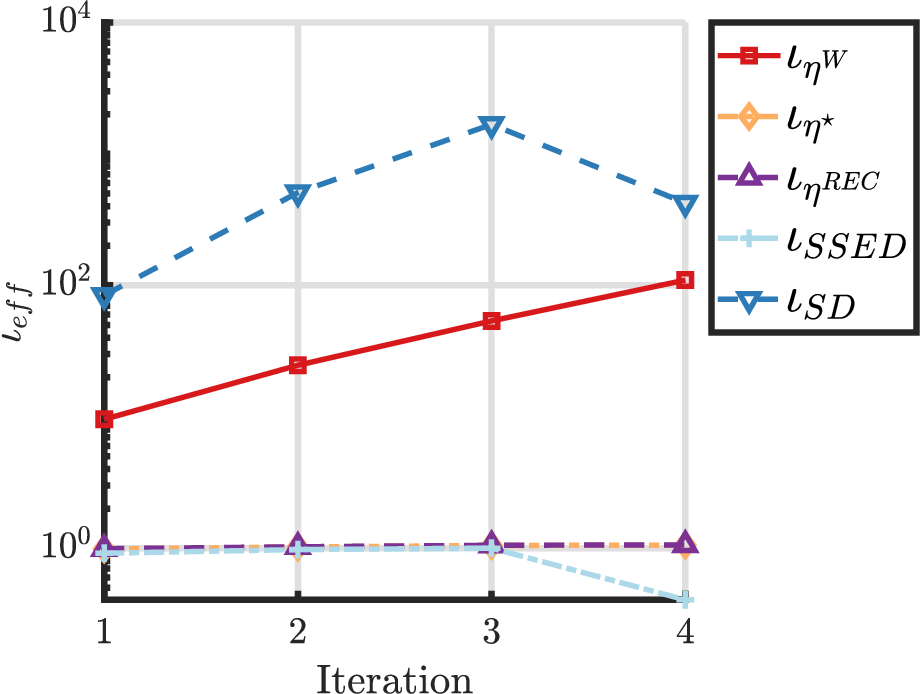}     \\
        \multicolumn{2}{c}{$p=3$}\\
           \includegraphics[width=0.33\linewidth]{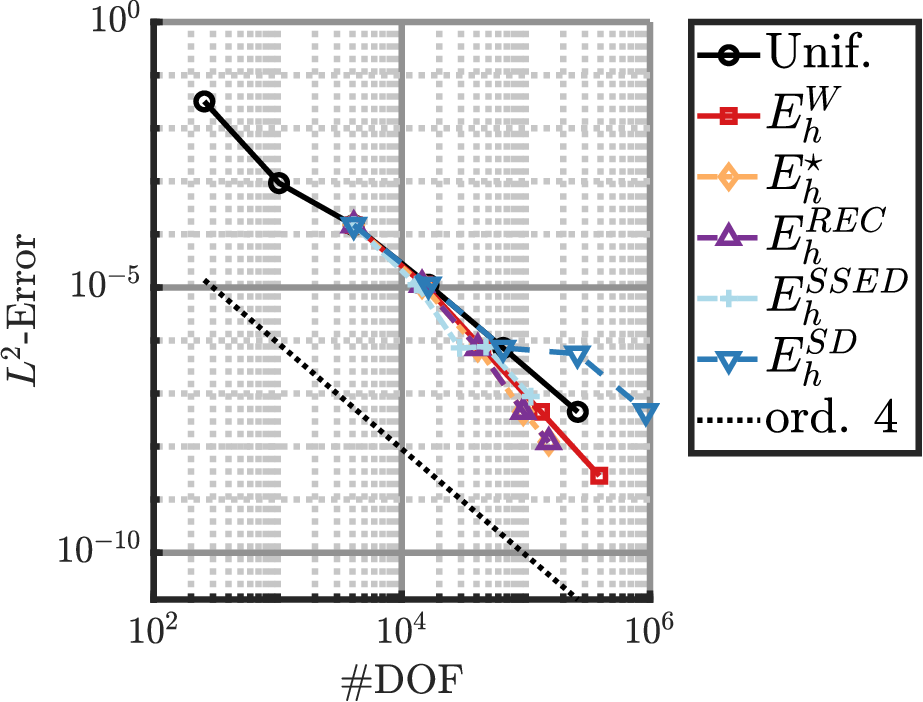} &  \includegraphics[width=0.33\linewidth]{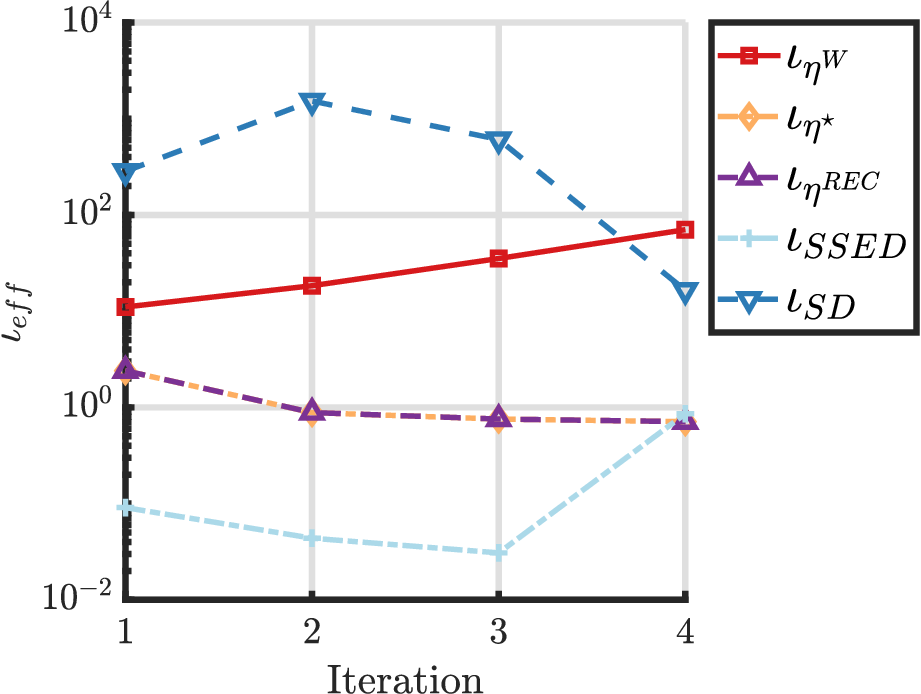}    
    \end{tabular}
    \caption{ The first column depicts the $L^2$-Errors and second the effectivity indices achieved under different refinement indicators for the two-dimensional Gaussian IC.}
    \label{fig:gauss_2D_varypoly}
    \end{figure}

With respect to the Gaussian problem, we observe lower errors in most cases with reduced numbers of degrees of freedom (see Column 1 Figure \ref{fig:gauss_2D_varypoly}), except for the spectral decay and multiwavelet indicators, for which over-refinement occurs. Looking at the effectivity indices (see Column 2 Figure \ref{fig:gauss_2D_varypoly}), we observe that the multiwavelet indicator overestimates the error at a given level, though typically not as bad as the SD indicator for lower orders. The other SIAC-based indicators generally appear to track the actual error well, if slightly underestimating it, though this does lead to reduced error and reduced numbers of degrees of freedom. The SSED and SD indicators do not perform well, either overestimating the error dramatically in the case of the SD indicator, or underestimating the error too much for $p>1$ in the SSED case. This prevents a reduced error compared to the reference uniform approximation. 

In Figure \ref{fig:wavelet_meshes}, we plot for illustrative purposes the sequence of meshes obtained after four enhancements using the multiwavelet indicator for the $p=2$ case. This shows that the procedure increases the refinement about the sharp gradient of the Gaussian spike, and provides less refinement elsewhere in the domain. For comparison between different error indicators, we also depict the final meshes achieved under different indication schemes for the $p=2$ case in \ref{fig:select_final_meshes}. Notice that all but the SD indicator refine about the center of the domain, whereas the SD indicator refines practically everywhere.

\input{Figures/meshes/select_final_meshes}

\begin{figure}[hp!]
    \centering
    \begin{tabular}{c c }
      \multicolumn{1}{c}{Errors} &  \multicolumn{1}{c}{Effectivity Index} \\
      \multicolumn{2}{c}{$p=0$}\\
       \includegraphics[width=0.33\linewidth]{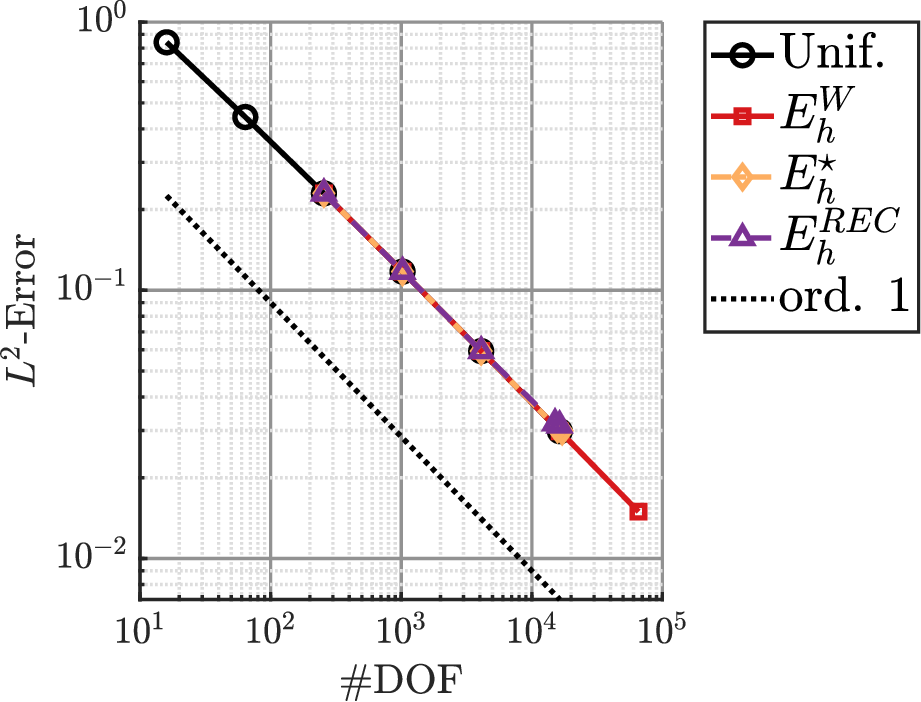} &  \includegraphics[width=0.33\linewidth]{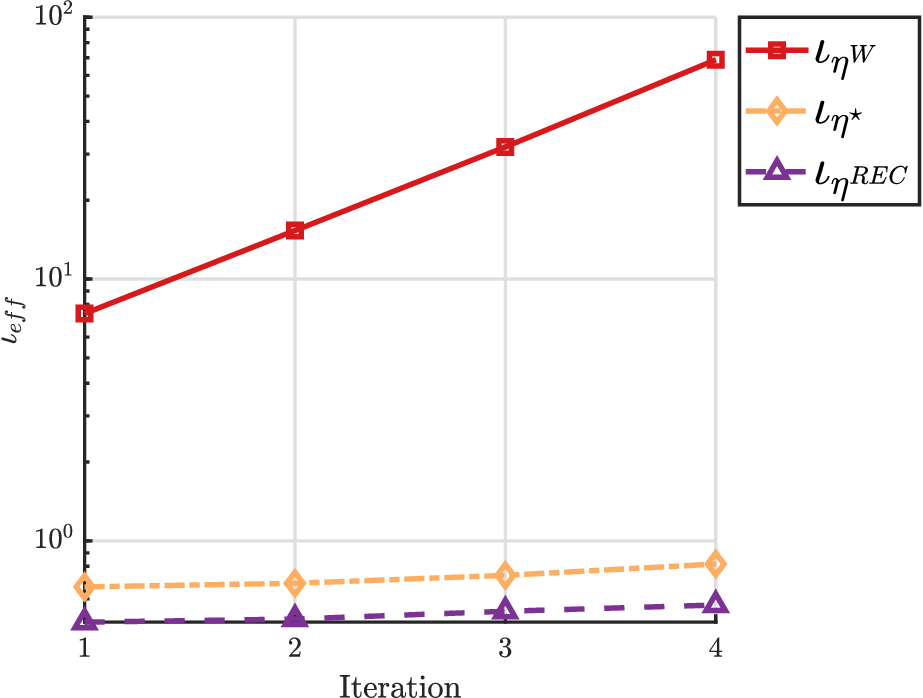} \\
             \multicolumn{2}{c}{$p=1$}\\
        \includegraphics[width=0.33\linewidth]{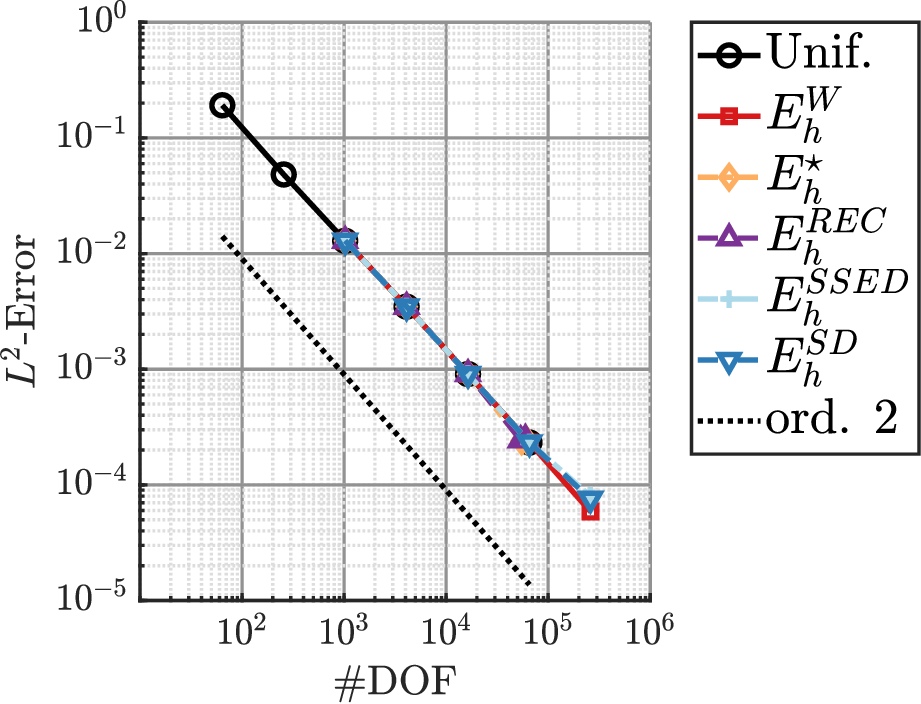} &  \includegraphics[width=0.33\linewidth]{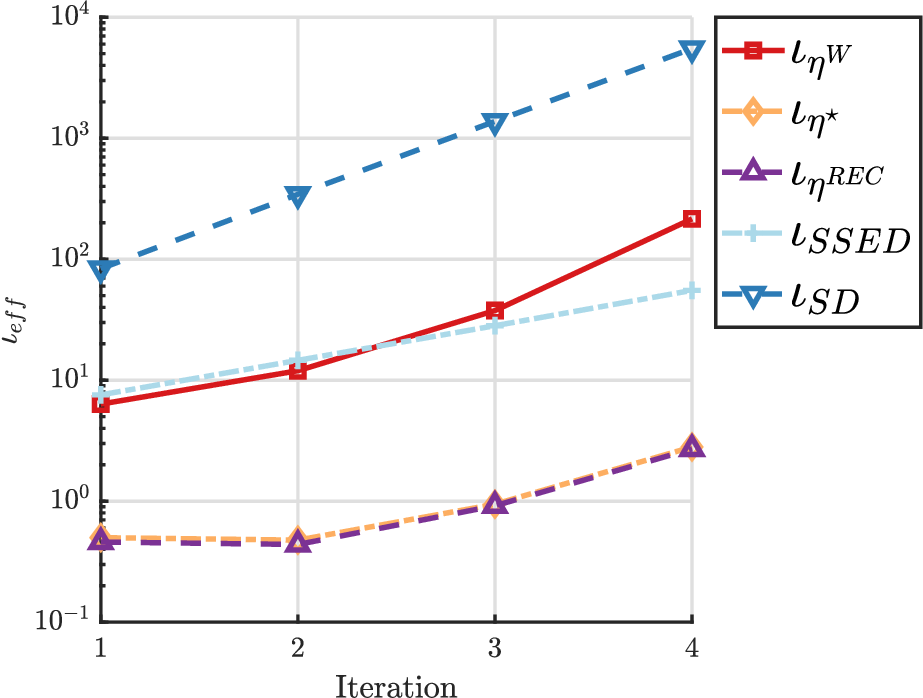} \\
        \multicolumn{2}{c}{$p=2$}\\
           \includegraphics[width=0.33\linewidth]{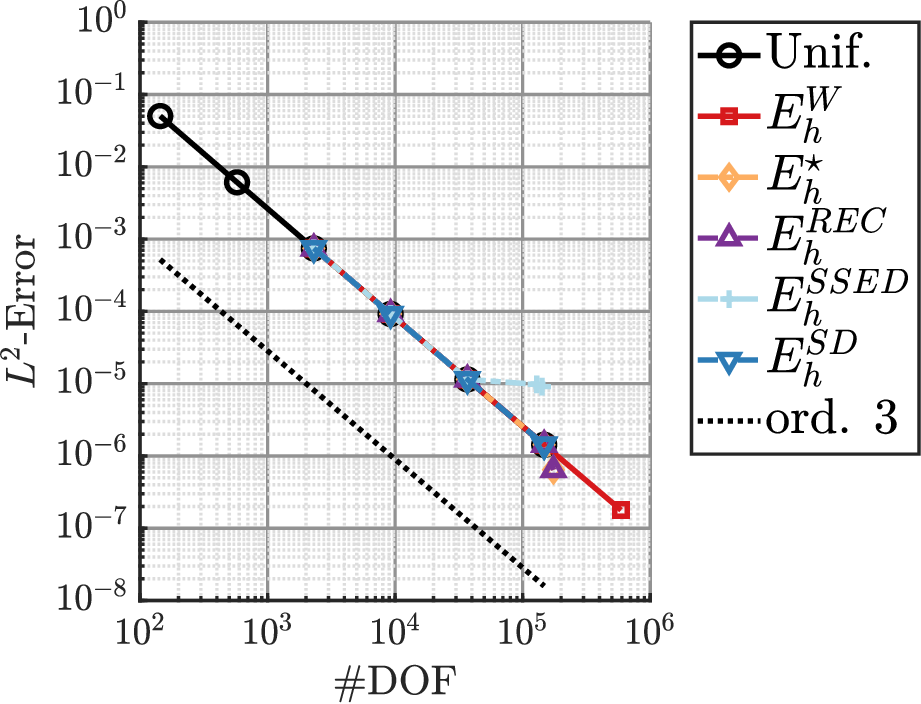} &  \includegraphics[width=0.33\linewidth]{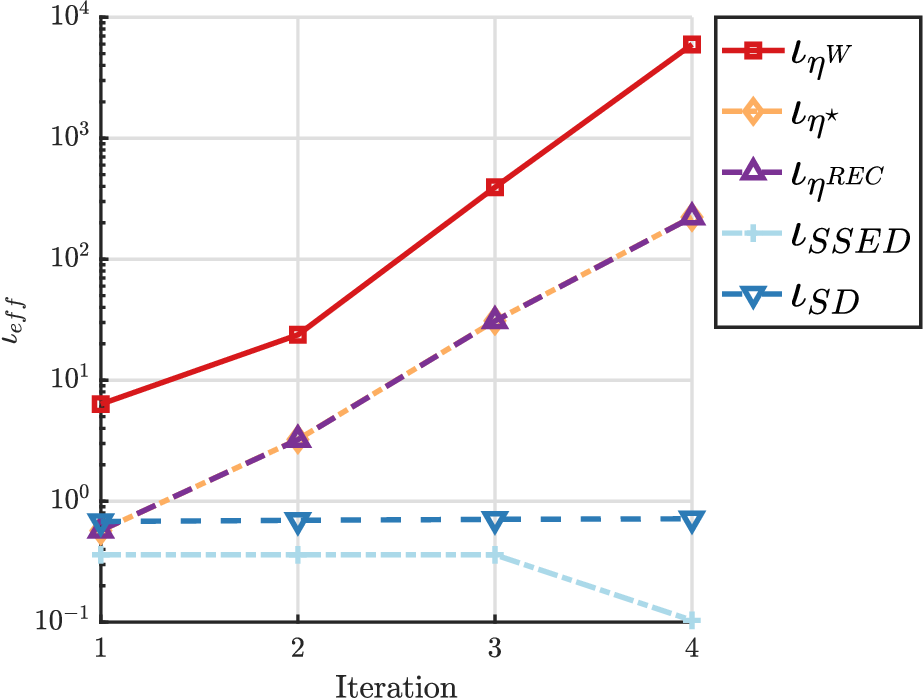}\\
                   \multicolumn{2}{c}{$p=3$}\\
           \includegraphics[width=0.33\linewidth]{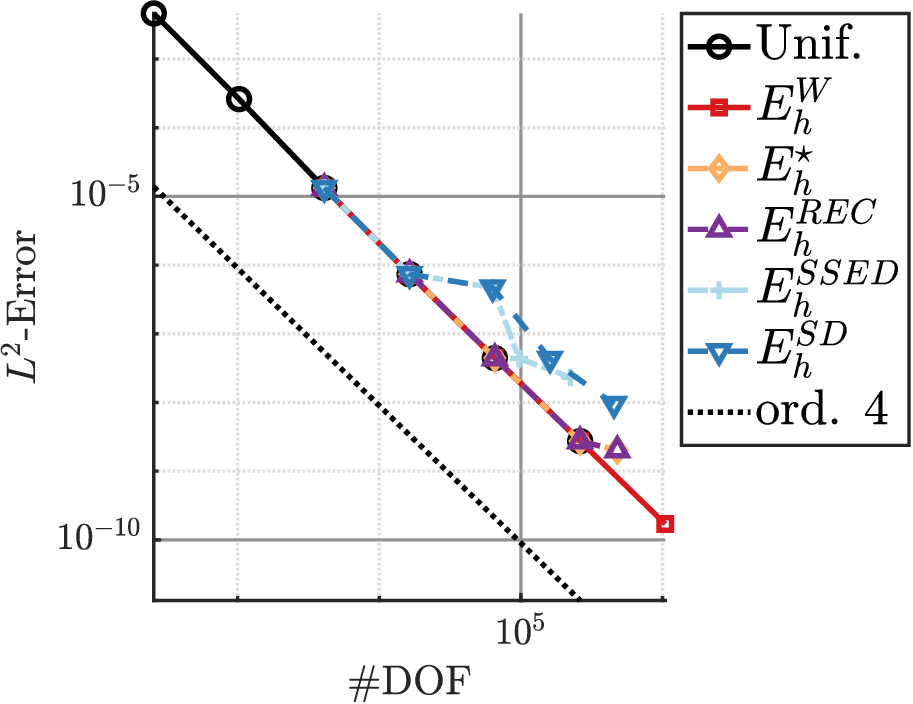} &  \includegraphics[width=0.33\linewidth]{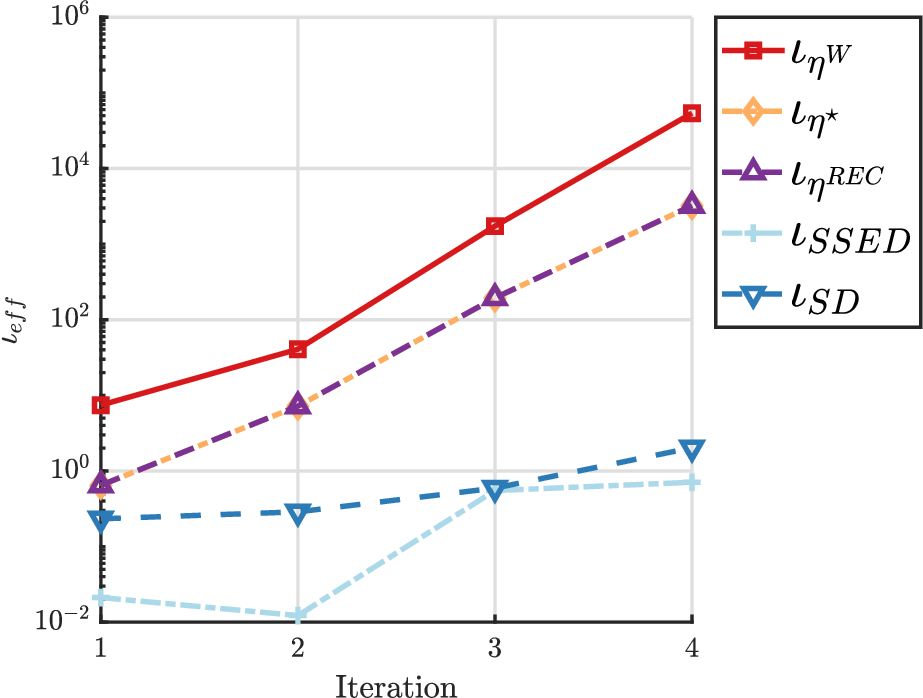}
    \end{tabular}
    \caption{The first and second column depict the $L^2$-Errors and effectivity indices achieved under different refinement indicators for the two-dimensional tensor sine IC.}
    \label{fig:sine_2D_varypoly}
    \end{figure}

Lastly, we consider the smooth sine problem (see Figure \ref{fig:sine_2D_varypoly}). Owing to the lack of sharp gradients, the different refinement indicators refine everywhere trying to reach the desired error tolerance. Generally the SIAC-based indicators perform well and reach a similar error level as the uniform reference,  though they start to over-refine for $p>1$. The multiwavelet indicator again tends to overestimate the error and as a result over-refines. The SD and SSED indicators flip between overestimation ($p=1$) and underestimation $(p=2,3)$. Outside of the $p=2$ case for the SD indicator, these indicators either over-refine and overshoot the desired error target, or refine in the wrong places and increase the number of degrees of freedom without substantially affecting the global error.

\subsection{Summary of Adaptivity results}
Applying the adaptivity procedure in one and two dimensions, we observe a few trends. The SIAC-based indicators $E^{\star}_h$ and $E^{REC}_h$ perform well in all  one-dimensional cases considered, tracking the uniform error well for smooth problems, and refining about areas of higher error for problems with sharp features, resulting in a reduction of the total number of degrees of freedom. In higher dimensions they begin to overestimate the error, perhaps suggesting that the one-dimensional support of the LSIAC filter may be too restrictive. The SD and SSED do not perform as well in one dimension, where they over-refine for smooth problems, and overestimate the error for problems with sharp features. In two dimensions, they perform erratically, either underestimating, or overestimating the actual error leading to varied performance. 
Lastly, the multiwavelet indicator has mixed performance in one dimension, overestimating the actual error, but achieving reductions in total numbers of degrees of freedom for problems possessing sharp features. In two dimensions the multiwavelet indicator consistently overestimates the error, which while less impactful for the Gaussian problem, leads to over-refinement and consequently poor performance for the sine problem. It is possible that the line support filter is too minimal in its averaging support, and that a tensor-product style filter may be prove better in the context of the multiwavelet indicator. In total, of the methods considered, the SIAC-based indicators appear most consistent in attaining error-reductions at reduced numbers of degrees of freedom.

%% file: Figures/meshes/select_final_meshes.tex
\begin{figure}
\centering
    \begin{tabular}{c c c  }
    \multicolumn{1}{c}{Initial Mesh} & \multicolumn{1}{c}{One Enhancement} & \multicolumn{1}{c}{Two Enhancements}\\
  \includegraphics[width=0.2\linewidth]{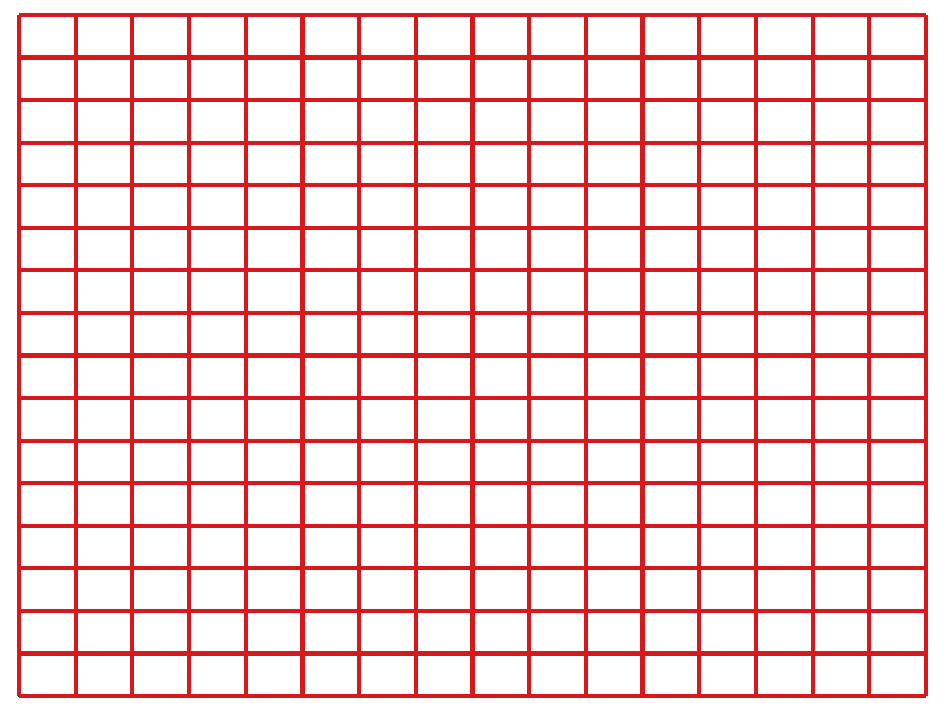}        &  
    \includegraphics[width=0.2\linewidth]{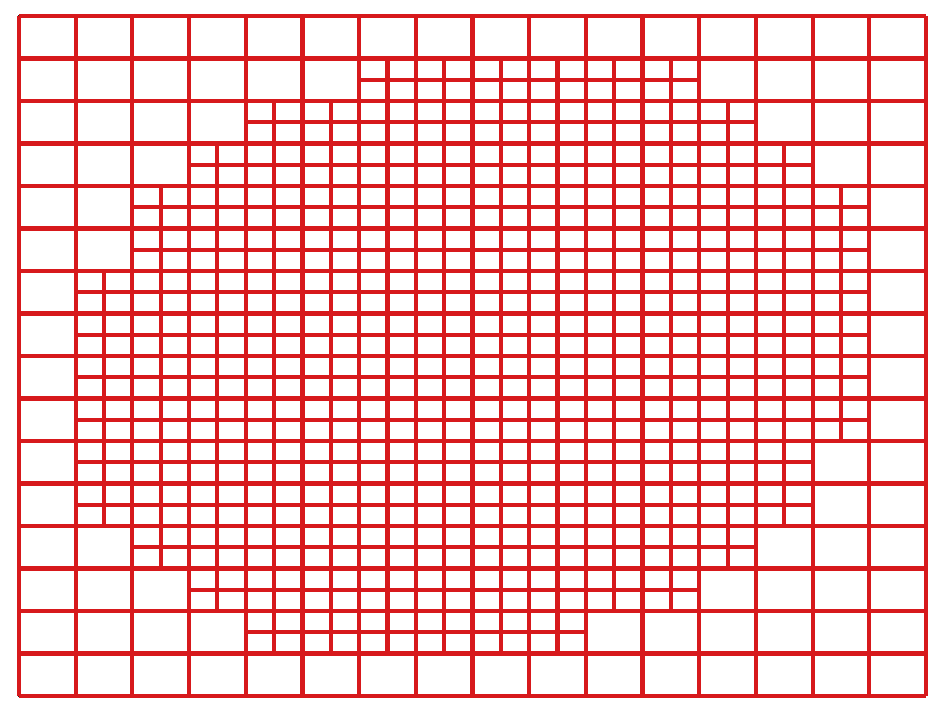}       
     & \includegraphics[width=0.2\linewidth]{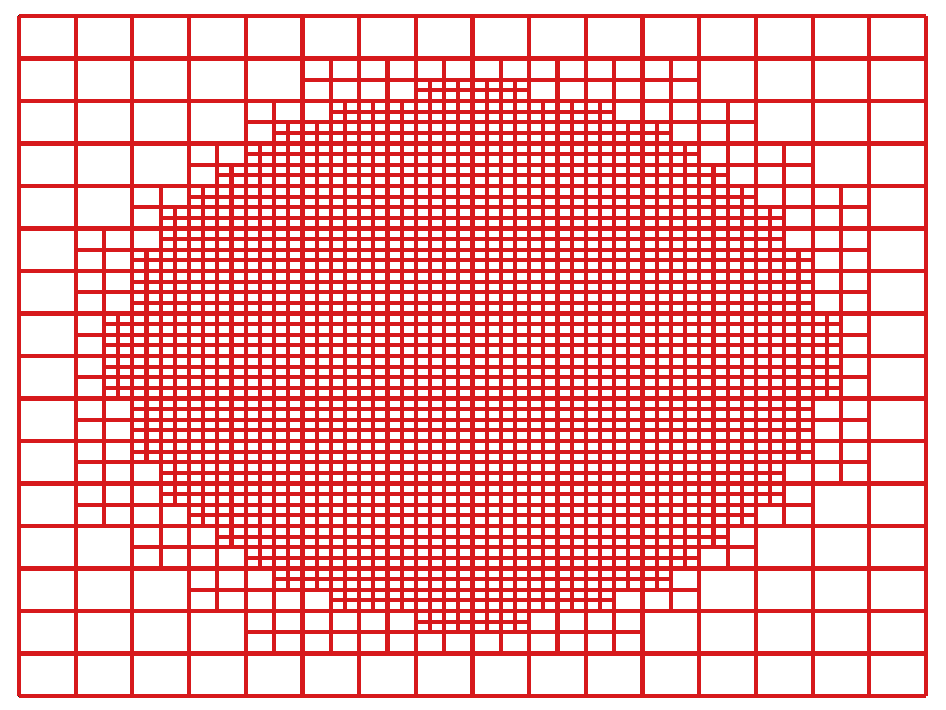}  
         \end{tabular}\\
\begin{tabular}{c c}
\multicolumn{1}{c}{Three Enhancements} & \multicolumn{1}{c}{Four Enhancements}\\
      \includegraphics[width=0.2\linewidth]{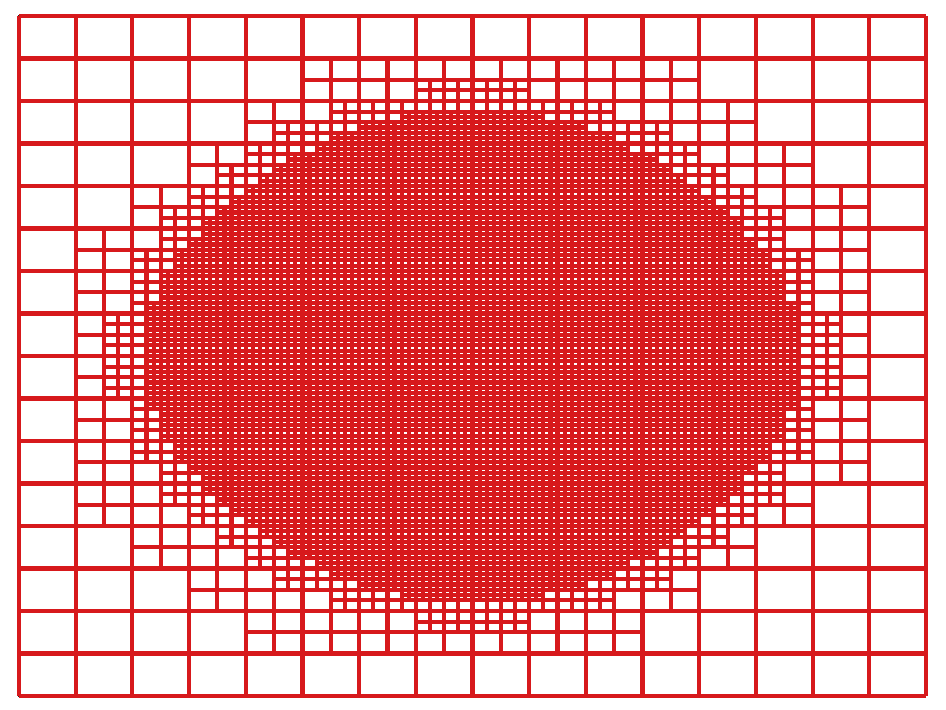}& 
          \includegraphics[width=0.2\linewidth]{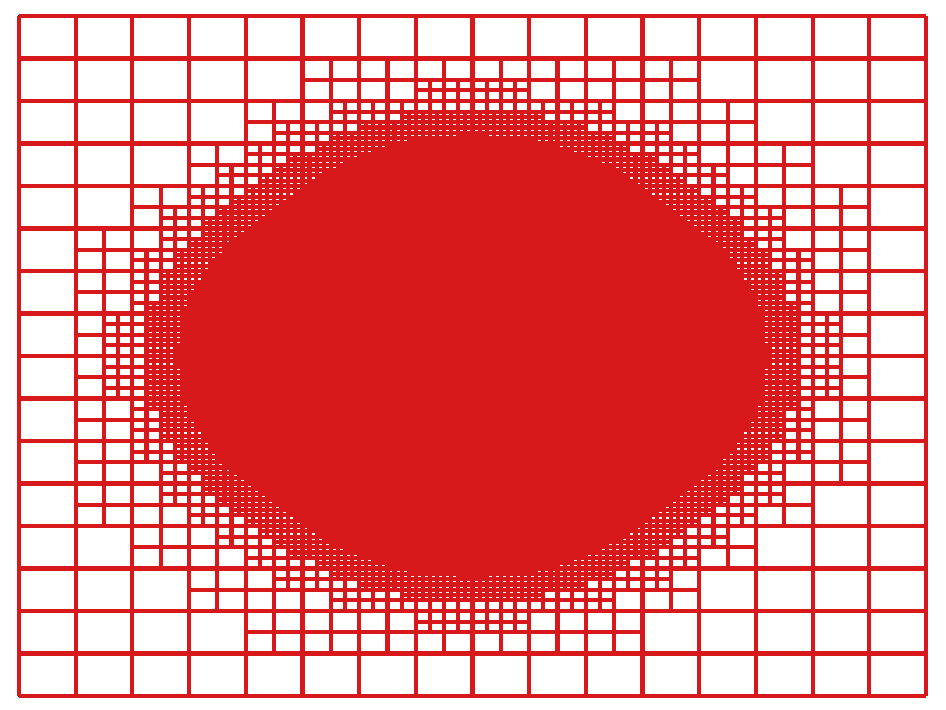} 
\end{tabular}
  
    \caption{For illustrative purposes we depict the meshes produced by the multiwavelet-based error indicator for $p=2$. Over repeated enhancements the resolution about the sharp Gaussian spike in the interior increases while the resolution near the boundary remains coarse. }
    \label{fig:wavelet_meshes}
\end{figure}

\begin{figure}
\centering
    \begin{tabular}{c c c }
    \multicolumn{1}{c}{$\eta^{W}$ Indicator} & \multicolumn{1}{c}{$\eta^{\star}$ Indicator} & \multicolumn{1}{c}{$\eta^{REC}$ Indicator}\\
  \includegraphics[width=0.2\linewidth]{Figures/meshes/fxd_w_gauss_periodic_mesh_p=2_enh=4.eps}        &  
    \includegraphics[width=0.2\linewidth]{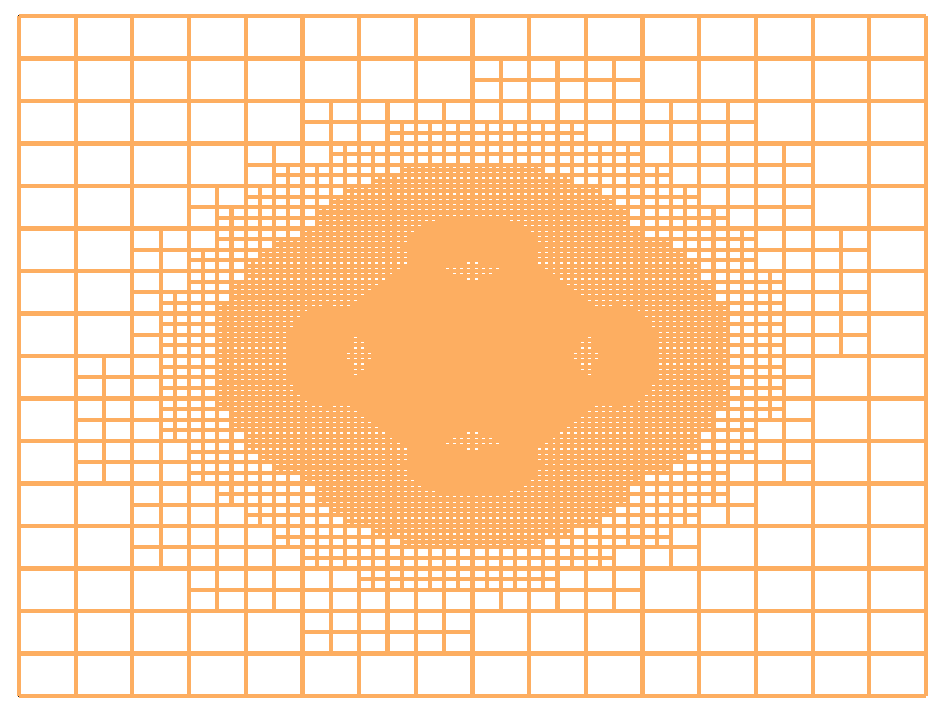}       
     & \includegraphics[width=0.2\linewidth]{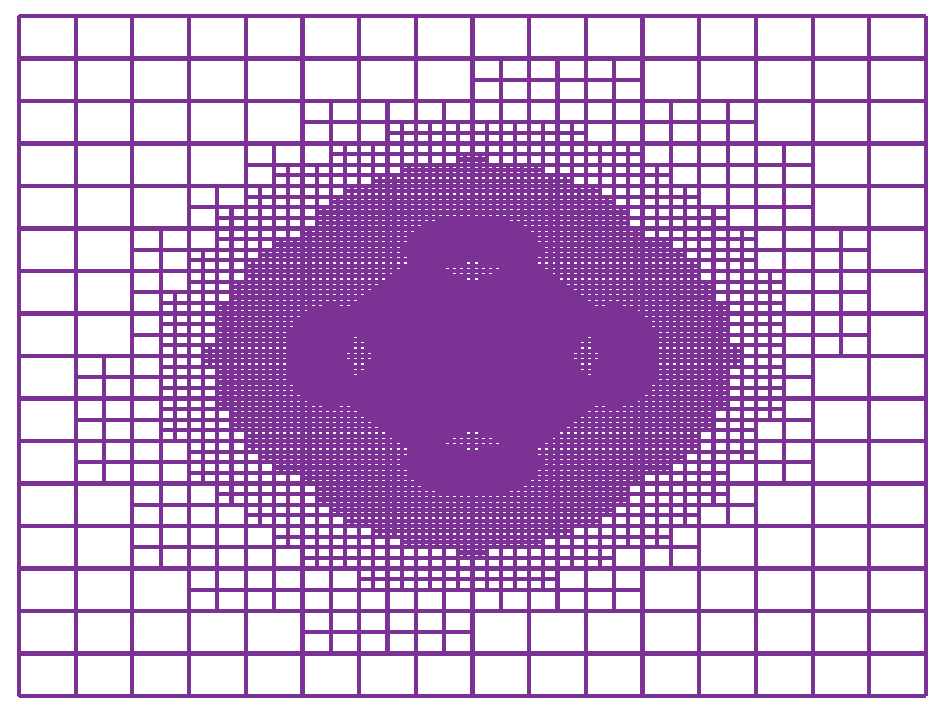}  
         \end{tabular}\\
\begin{tabular}{c c}
\multicolumn{1}{c}{$\eta^{SD}$ Indicator} & \multicolumn{1}{c}{$\eta^{SSED}$}\\
      \includegraphics[width=0.2\linewidth]{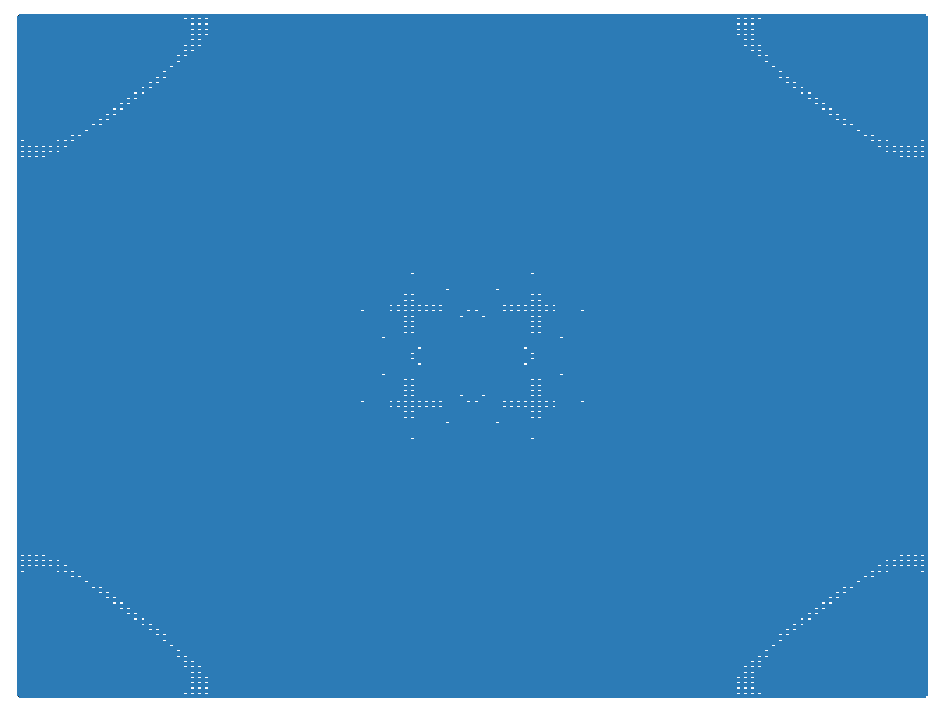} & 
          \includegraphics[width=0.2\linewidth]{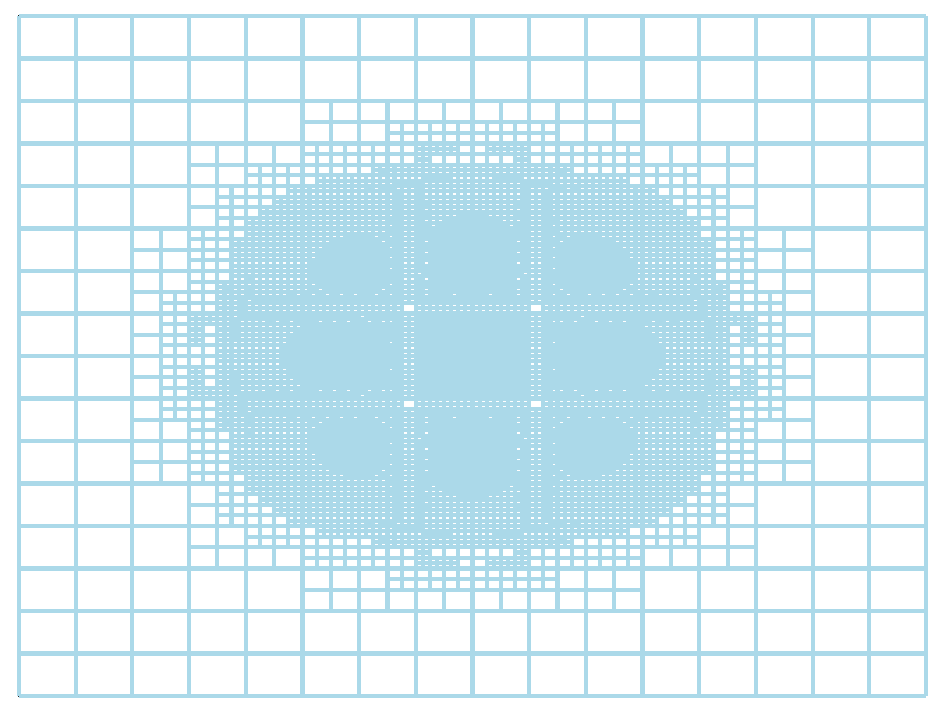}
\end{tabular}
  
    \caption{Comparison of final meshes generated after the last refinement for the Gaussian test problem with $p=2$ for different error indicators.}
    \label{fig:select_final_meshes}
\end{figure}

%% file: sections/conclusions.tex
\section{Conclusions}\label{sec:conclusions}
In this paper we have demonstrated the applicability of (L)SIAC-MRA to nonuniform meshes including perturbed quadrilateral meshes and Delaunay triangulations. We have also confirmed that the inclusion of the adaptive kernel scaling of \cite{Jallepalli2019B} serves to improve error reduction for meshes with widely varying mesh sizes as compared to the standard choice of the maximum element edge length. With respect to adaptivity, the enhanced reconstructions from (L)SIAC-MRA can be used to inform adaptation in DG schemes in a manner similar to \cite{Bautista_paper}. Furthermore, adapted solutions based off SIAC-based indicators can be constructed to have fewer degrees of freedom while maintaining similar errors to the uniformly refined solutions. In the future, we anticipate that the (L)SIAC-MRA techniques will prove useful in the resolution of sub-grid dynamics within multi-scale turbulence modelling. An investigation of this application and of the effect of tensor-product filters on the multiwavelet indicator will be pursued in future work.